\input amstex
\documentstyle{amsppt}\nologo\footline={}\subjclassyear{2000}
\input epsf

\def\PU{\mathop{\text{\rm PU}}}
\def\T{\mathop{\text{\rm T}}}
\def\SU{\mathop{\text{\rm SU}}}
\def\ta{\mathop{\text{\rm ta}}}
\def\Lin{\mathop{\text{\rm Lin}}}
\def\B{\mathop{\text{\rm B}}}
\def\R{\mathop{\text{\rm R}}}
\def\Im{\mathop{\text{\rm Im}}}
\def\Arg{\mathop{\text{\rm Arg}}}
\def\G{\mathop{\text{\rm G}}}
\def\Area{\mathop{\text{\rm Area}}}
\def\Re{\mathop{\text{\rm Re}}}
\def\tr{\mathop{\text{\rm tr}}}
\def\dist{\mathop{\text{\rm dist}}}
\def\L{\mathop{\text{\rm L}}}
\def\U{\mathop{\text{\rm U}}}
\def\Tn{\mathop{\text{\rm Tn}}}

\hsize450pt\topmatter\title\vskip-40pt Complex hyperbolic structures on
disc bundles over surfaces\vskip-15pt\endtitle\author Sasha Anan$'$in,
Carlos H.~Grossi, Nikolay Gusevskii\vskip-5pt\endauthor\address
Departamento de Matem\'atica, IMECC, Universidade Estadual de
Campinas,\newline13083-970--Campinas--SP, Brasil\endaddress\address
Max-Planck-Institut f\"ur Mathematik, Vivatsgasse 7, 53111 Bonn,
Germany\endaddress\email grossi$_-$ferreira\@yahoo.com\endemail\address
Departamento de Matem\'atica, ICEX, Universidade Federal de Minas
Gerais,\newline31161-970--Belo Horizonte--MG,
Brasil\endaddress\subjclass 57S30 (30F35, 51M10,
57M50)\endsubjclass\abstract We study complex hyperbolic disc bundles
over closed orientable surfaces that arise from discrete and faithful
representations $H_n\to\PU(2,1)$, where $H_n$ is the fundamental group
of the orbifold $\Bbb S^2(2,\dots,2)$ and thus contains a surface group
as a subgroup of index $2$ or $4$. The results obtained provide the
first complex hyperbolic disc bundles $M\to\Sigma$ that

$\bullet$ admit both real and complex hyperbolic structures,

$\bullet$ satisfy the equality $2(\chi+e)=3\tau$,

$\bullet$ satisfy the inequality $\frac12\chi<e$, and

$\bullet$ induce discrete and faithful representations
$\pi_1\Sigma\to\PU(2,1)$ with fractional Toledo invariant,

\noindent
where $\chi$ is the Euler characteristic of $\Sigma$, $e$ denotes the
Euler number of $M$, and $\tau$ stands for the Toledo invariant of $M$.
To get a satisfactory explanation of the equality $2(\chi+e)=3\tau$, we
conjecture that there exists a holomorphic section in all our examples.

In order to reduce the amount of calculations, we systematically
explore coordinate-free
methods.\endabstract\endtopmatter\document\rightheadtext{Complex
Hyperbolic Structures on Disc Bundles over Surfaces}\leftheadtext{Sasha
Anan$'$in, Carlos H.~Grossi, Nikolay Gusevskii}

\vskip-70pt

{\eightpoint

\centerline{CONTENTS}

\qquad\qquad\qquad\qquad1.~Introduction\hfill1\phantom{aaaaaa}

\qquad\qquad\qquad\qquad\quad Outline of the general
construction\hfill3\phantom{aaaaaa}

\qquad\qquad\qquad\qquad\quad Acknowledgements\hfill4\phantom{aaaaaa}

\qquad\qquad\qquad\qquad\quad Some notation\hfill4\phantom{aaaaaa}

\qquad\qquad\qquad\qquad2.~General construction\hfill4\phantom{aaaaaa}

\qquad\qquad\qquad\qquad\quad2.1.~Preliminaries. Cycle of bisectors and
the Toledo invariant\hfill4\phantom{aaaaaa}

\qquad\qquad\qquad\qquad\quad2.2.~Discreteness\hfill9\phantom{aaaaaa}

\qquad\qquad\qquad\qquad\quad2.3.~Simplicity and
transversality\hfill10\phantom{aaaaaa}

\qquad\qquad\qquad\qquad\quad2.4.~Fibred polyhedra. Euler
number\hfill12\phantom{aaaaaa}

\qquad\qquad\qquad\qquad\quad2.5.~Transversal
triangles\hfill14\phantom{aaaaaa}

\qquad\qquad\qquad\qquad3.~Series of explicit
examples\hfill16\phantom{aaaaaa}

\qquad\qquad\qquad\qquad\quad3.1.~A couple of transversal
triangles\hfill16\phantom{aaaaaa}

\qquad\qquad\qquad\qquad\quad3.2.~Examples of disc
bundles\hfill20\phantom{aaaaaa}

\qquad\qquad\qquad\qquad\quad3.3.~Some interesting
examples\hfill25\phantom{aaaaaa}

\qquad\qquad\qquad\qquad4.~Appendix: technical tools, elementary
algebraic and geometric background\hfill28\phantom{aaaaaa}

\qquad\qquad\qquad\qquad\quad4.1.~Preliminaries: hermitian metric,
complex geodesics, geodesics,\hfill\phantom{aaaaaa}

\qquad\qquad\qquad\qquad\quad\quad\quad $\Bbb R$-planes, and
bisectors\hfill29\phantom{aaaaaa}

\qquad\qquad\qquad\qquad\quad4.2.~Bisectors\hfill34\phantom{aaaaaa}

\qquad\qquad\qquad\qquad\quad4.3.~Bisectors with common
slice\hfill38\phantom{aaaaaa}

\qquad\qquad\qquad\qquad\quad4.4.~Triangles of
bisectors\hfill43\phantom{aaaaaa}

\qquad\qquad\qquad\qquad\quad4.5.~K\"ahler
primitive\hfill49\phantom{aaaaaa}

\qquad\qquad\qquad\qquad5.~References\hfill51\phantom{aaaaaa}}

\bigskip

\centerline{\bf1.~Introduction}

\medskip

This paper is devoted to a particular case of the traditional question
concerning the interrelation between topology and geometry: Does a
topological $4$-manifold admit a specific geometric structure (real
hyperbolic, complex hyperbolic, quaternionic, etc.)? Henceforth, $M$
denotes an oriented disc bundle over a closed orientable surface
$\Sigma$ of Euler characteristic $\chi$ and $e$ stands for the Euler
number of the bundle.

For the existence of a real hyperbolic structure on $M$, there are
various conditions in terms of $\chi$ and $e$ (see [GLT], [Kui], [Luo],
and also [Kap1, Kap2]). We study the case of complex hyperbolic
geometry. In~this case, there is one more discrete invariant, the
Toledo invariant of the representation $\pi_1\Sigma\to\PU(2,1)$
provided by $M$. It is related to the complex ($=$ Riemannian)
structure on $M$, takes values in $\frac23\Bbb Z$, and satisfies the
inequality $|\tau|\le|\chi|$ [Tol]. The Toledo invariant is the only
discrete invariant of a (not necessarily discrete) representation
$\pi_1\Sigma\to\PU(2,1)$ [Xia].

There are not so many known complex hyperbolic disc bundles. Simple
types of such bundles are the $\Bbb R$-Fuchsian bundles (they satisfy
$e=\chi$ and $\tau=0$) and the $\Bbb C$-Fuchsian bundles (they are
characterized by $\chi=\tau$ [Tol] and satisfy $e=\frac12\chi$). The
first complex hyperbolic bundles different from $\Bbb R$- and
$\Bbb C$-Fuchsian ones were constructed in [GKL]. They satisfy the
relations $e=\chi+|\tau/2|$ and $\chi\le e\le\frac12\chi$. Thus,
$\Bbb R$- and $\Bbb C$-Fuchsian bundles provide the extreme values of
$e$ in all examples known before the present paper.

Here, we consider a class of discrete and faithful representations that
is slightly simpler than that of general representations
$\pi_1\Sigma\to\PU(2,1)$. Let $H_n$ denote the group generated by
$r_1,\dots,r_n$ with the defining relations $r_n\dots r_1=1$ and
$r_i^2=1$, $i=1,\dots,n$. For a hyperelliptic $\Sigma$, the group $H_n$
is nothing but the index $2$ extension of $\pi_1\Sigma$ by the
hyperelliptic involution of $\Sigma$ and $n$ is even (for~odd~$n$, a
surface group is a subgroup of index $4$ in $H_n$). We study discrete
and faithful representations $\varrho:H_n\to\PU(2,1)$ such that the
$\varrho r_i$'s are reflections in ultraparallel complex geodesics. In
order to prove discreteness, we~construct fundamental polyhedra bounded
by cycles of segments of bisectors (see Outline at the end of
Introduction for more details). The main difficulties at these stages
are to check that segments of bisectors intersect properly and the fact
that adjacent segments have nonconstant angle along their intersection.

The oriented transversal triangles of bisectors studied in Subsections
2.5 and 4.4 serve as building blocks for fundamental polyhedra. The
vertices of such a triangle are complex geodesics ($=$ Poincar\'e
discs) naturally equipped with an isometry called the holonomy of the
triangle. Every oriented transversal triangle of bisectors bounds a
(closed) fibred polyhedron, i.e., a disc bundle over a disc that
extends the slice bundle structure of the triangle itself. This fact is
established by showing that the space of oriented transversal triangles
of bisectors is path-connected and that a triangle of bisectors with
common complex spine bounds a fibred polyhedron. It is crucial that the
holonomy of a triangle cannot be R-parabolic (see 2.5.1 for the
definition) or equal to the identity. Moreover, every oriented
transversal triangle of bisectors possesses a `fractional Euler
quantity' that can actually be interpreted as its Euler number. The
fractional Euler quantity is an arc on the ideal boundary of a vertex
of the triangle. It is determined by the holonomy of the triangle and
is additive with respect to gluing triangles.

Gluing a couple of transversal triangles of bisectors, we obtain in
Section 3 a quadrangle of bisectors bounding a fundamental polyhedron
for the group generated by $U,W\in\PU(2,1)$ with the relations
$U^n=W^n=(U^{-1}W)^2=1$. The quadrangle furnishes a complex hyperbolic
disc bundle over the turnover orbifold $\Bbb S^2(n,n,2)$ whose rational
Euler number (see [BSi] for the definition) can be inferred from the
fractional Euler quantities of the involved triangles. Since $H_n$ is a
subgroup of index $n$ in the turnover group, we get a complex
hyperbolic disc bundle $M$ over a closed surface. Clearly,
we~simultaneously obtain a compact $3$-manifold (a circle bundle over a
closed surface) admitting a spherical CR structure.

In the above way, we construct a huge family of explicit examples of
disc bundles $M$. All of them satisfy the inequality $\frac12\chi<e$
and the equality $2(\chi+e)=3\tau$ (with negative $\chi,e,\tau$). The
inequality was never valid for previously known examples and the
equality was valid only in the $\Bbb C$-Fuchsian case. Since the
equality is a necessary condition for the existence of a holomorphic
section of $M$ (in the sense that there are a disc bundle structure on
$M$ and a smooth holomorphic surface $\Sigma\subset M$ that intersects
every fibre exactly once), we conjecture that there exists a
holomorphic section in all our examples. For many examples, $\tau$ is
not integer, which implies in particular that the corresponding
representation $\pi_1\Sigma\to\PU(2,1)$ cannot be lifted to $\SU(2,1)$.
We obtain the first disc bundles $M$ admitting both real and complex
hyperbolic structures. Passing to the corresponding circle bundles, we
see that there exist circle bundles over closed orientable surfaces
admitting simultaneously conformally flat and spherical CR structures.
For a more detailed discussion about the features of our examples and
related conjectures, see Subsection 3.3.

Whenever reasonable, we work without coordinates. This concerns
especially the appendix (Section 4) which is devoted to an exposition
of algebraic and geometric background and contains the proofs absent in
Sections 2--3. We believe that some facts in the appendix are
interesting {\it per se.} This includes: explicit formulae for the
normal vector to an oriented bisector (Proposition 4.2.11) and for the
angle between two cotranchal bisectors (Lemma 4.3.1); a numerical
transversality criterion for a pair of cotranchal bisectors (Criterion
4.3.3); a sort of metric separability of transversal cotranchal
bisectors (Lemma 4.3.6) that is a key point allowing us to apply
Poincar\'e's Polyhedron Theorem; and an explicit formula for a K\"ahler
primitive in 4.5.1.

The reader unfamiliar with complex hyperbolic geometry may begin
reading with Section 4. The reader familiar with complex hyperbolic
geometry may first read Sections 2--3 (the geometric core of the
paper).

\smallskip

{\bf Outline of the general construction.} Roughly speaking, the
general construction of a complex hyperbolic disc bundle $M$ deals with
a certain fundamental polyhedron $F$ for the group
$\pi_1M\subset\PU(2,1)$. The polyhedron $F$ is bounded by a cycle of
oriented segments of bisectors $(\B_1,\dots,\B_m)$ such that the final
slice of $\B_i$ is the initial slice of $\B_{i+1}$ for every $i$ (the
indices are modulo $m$). Suppose that the obvious slice bundle
structure of the boundary $\partial_0F:=\B_1\cup\dots\cup\B_m$ of $F$
is extendable to $F$. Then, taking into account that the
identifications of the segments by means of elements in $\pi_1M$
preserve the slice bundle structure of $\partial_0F$, we obtain a
complex hyperbolic disc bundle $M$.

For technical reasons, the polyhedron $F$ is glued from two or four
polyhedra $P$, each bounded by a cycle $\Cal C$ of segments
$(\B_1,\dots,\B_n)$. The identification of faces of $P$ is given by the
reflections $R_i$ in the middle slices $M_i$ of the $\B_i$'s which
satisfy the relation $R_n\dots R_1=1$. If we show that $P$ is a
fundamental polyhedron for the group $H$ generated by the $R_i$'s and
that the slice bundle structure of $\partial_0P$ is extendable to $P$,
then we arrive at the desired complex hyperbolic disc bundle $M$.

The solid torus $\partial_0P$ can be readily fibred into circles called
meridional curves of $\Cal C$ : in its segment $b_i$ belonging to
$\B_i$, every meridional curve $b$ is equidistant from the real spine
of $\B_i$ and lies in a suitable meridian of $\B_i$. Moreover,
$R_ib_i=b_i$ and $R_i$ acts on $b_i$ as a `reflection' in the middle
point of $b_i$, i.e., in~the point $b_i\cap M_i$. In other words, the
reflections $R_i$'s preserve the fibering of $\partial_0P$ by
meridional curves. Every meridional curve of the cycle $\Cal C$ is
intended to generate a section of the disc bundle, i.e., to be the
boundary of a disc that provides a section.

Considering those `formally' adjacent polyhedra congruent to $P$ that
are intended to tessellate a neighbourhood of the slice $S_{i+1}$
common to the segments $\B_i$ and $\B_{i+1}$, we can see that the sum
of the interior angles of these adjacent polyhedra at a point
$s\in S_{i+1}$ (called the total angle at $s$) is an integer multiple
of $2\pi$ and that the total angle at $s$ is independent of the choice
of $s$ in $S_{i+1}$. Requiring that the total angle is $2\pi$ and that
the full bisectors ${\prec}\B_i{\succ}$ and ${\prec}\B_{i+1}{\succ}$ of
the segments $\B_i$ and $\B_{i+1}$ are transversal along their common
slice $S_{i+1}$ for all $i$ (such a cycle $\Cal C$ is said to be simple
transversal with total angle $2\pi$), we show that $P$ is indeed a
fundamental polyhedron for $H$.

It is not so easy to construct a simple transversal cycle of bisectors.
The problem is that the segments $\B_i$ and $\B_j$ can unexpectedly
intersect, that is, the cycle $\Cal C$ can be nonsimple and, thus, can
provide no polyhedron $P$. Under some additional requirement (appealing
to some kind of simple Poincar\'e duality), we show that the cycle is
simple.

In general, we cannot expect the slice bundle structure of
$\partial_0P$ to be extendable to $P$. It is quite possible that the
polyhedron $P$ is `knotted.' Although this is a topological problem, it
looks even more difficult than the problem of simplicity because it
requires to dwell into the mutual position of all the~$\B_i$'s (not
only to verify that the $\B_i$'s intersect properly). An easy way to
avoid such a situation is to show that the polyhedron $P$ can be cut
into some simple polyhedra fibred in a way compatible with gluing them
back. Besides, this helps to solve the problem of the simplicity of a
cycle. Transversal triangles of bisectors are good candidates for such
simple fibred polyhedra.

\smallskip

\newpage

{\bf Acknowledgements.} The authors are grateful to Victor Gerasimov
for the discussions in the beginning of this work. We are indebted to
Misha Kapovich for his interest in our work and many important
suggestions that considerably improved the exposition. We express our
gratitude to the referee for many helpful comments. We thank Israel
Vainsencher for a valuable remark. This article is dedicated to the
memory of Igor$'$ Vladimirovich L$'$vov.

\medskip

\centerline{\bf Some notation}

\medskip

\noindent$\Bbb B$\ \ {\bf2.1.1} 4

\noindent$\overline{\Bbb B}$\ \ {\bf2.1.1} 4

\noindent$\B^-$\ \ {\bf2.2.3} 9

\noindent$\B{\succ}$\ \ {\bf2.3.1} 10

\noindent${\prec}\B{\succ}$\ \ {\bf2.3.1} 10

\leftskip80pt

\vskip-61pt

\noindent$\B[S,S']$\ \ {\bf2.1.1} 5

\noindent$\B[p_1,p_2]$\ \ {\bf2.3.1} 10

\noindent$\B{\wr}g_1,g_2{\wr}$\ \ {\bf4.1.20} 33

\noindent$\partial_0$\ \ {\bf2.2.2} 9

\noindent$\partial_\infty$\ \ {\bf2.1.1} 4

\leftskip180pt

\vskip-62pt

\noindent$e_P$\ \ {\bf2.4.6} 13

\noindent$\G S$\ \ {\bf4.1.9} 31

\noindent$\G[p_1,p_2]$\ \ {\bf2.3.1} 10

\noindent$\G{\wr}p_1,p_2{\wr}$\ \ {\bf4.1.11} 31

\noindent$K^+$\ \ {\bf2.3.1} 10

\leftskip280pt

\vskip-62pt

\noindent$K^-$\ \ {\bf2.3.1} 10

\noindent$\L\G$\ \ {\bf4.1.9} 31

\noindent$\L{\wr}p_1,p_2{\wr}$\ \ {\bf4.1.6} 30

\noindent$\pi'[p]v$\ \ {\bf4.1.1} 29

\noindent$\pi[p]v$\ \ {\bf4.1.1} 29

\leftskip374pt

\vskip-62pt

\noindent$R(p)$\ \ {\bf2.1.3} 5

\noindent$\R S$\ \ {\bf4.1.9} 31

\noindent$\ta(p_1,p_2)$\ \ {\bf2.1.2} 4

\noindent$t_\varphi$\ \ {\bf2.1.11} 7

\noindent$v_p$\ \ {\bf2.1.11} 7

\leftskip0pt

\bigskip

\centerline{\bf2.~General construction}

\medskip

In this section, we develop some general tools applicable to
constructing complex hyperbolic disc bundles over closed orientable
surfaces.

\medskip

\centerline{\bf2.1.~Preliminaries. Cycle of bisectors and the Toledo
invariant}

\medskip

In this subsection, we introduce some basic notation, conventions, and
simple facts. Then we define configurations and cycles of bisectors,
meridional curves, and related concepts. Finally, we calculate the
Toledo invariant of representations defined by cycles of bisectors.

\smallskip

{\bf2.1.1.}~Fix, once and for all, a three-dimensional $\Bbb C$-vector
space $V$ equipped with a hermitian form $\langle-,-\rangle$ of
signature $++-$. The corresponding unitary, special unitary, and
projective unitary groups are respectively denoted by $\U$, $\SU$, and
$\PU$. Depending on the context, we will use the elements of $V$ to
denote the points in the complex projectivization $\Bbb P$ of $V$. In
general, $\Bbb PW$ denotes the complex projectivization of
$W\subset V$.

\smallskip

It is well known that the complex hyperbolic plane can be identified
with the open $4$-ball\footnote{Here and in what follows, the symbol
$:=$ stands for `equals by definition.'}
$$\Bbb B:=\big\{p\in\Bbb P\mid\langle p,p\rangle<0\big\}.$$
The complex hyperbolic distance between two points $p_1,p_2\in\Bbb B$
can be expressed in terms of the {\it tance}
$$\ta(p_1,p_2):=\frac{\langle p_1,p_2\rangle\langle
p_2,p_1\rangle}{\langle p_1,p_1\rangle\langle
p_2,p_2\rangle}\eqno{\bold{(2.1.2)}}$$
as explained in [Gol, p.~77]. Sometimes, it is convenient to use the
tance instead of the distance: the~latter is a monotonic function of
the former, whereas the tance has an algebraic nature and a simpler
form. Furthermore, the tance is involved in many complex hyperbolic
concepts.

The ideal boundary of the complex hyperbolic plane is the $3$-sphere
$$\partial_\infty\Bbb B:=\big\{p\in\Bbb P\mid\langle
p,p\rangle=0\big\}$$
formed by all isotropic points. We denote
$\overline{\Bbb B}:=\Bbb B\cup\partial_\infty\Bbb B$. Preferentially,
we consider geometrical objects such as geodesics, complex geodesics,
$\Bbb R$-planes, bisectors, etc.~as being closed in
$\overline{\Bbb B}$. The exceptions occur when explicitly stated or
when a concept or formula is clearly inapplicable to isotropic points.
In~general, $\partial_\infty X$ stands for
$\partial_\infty\Bbb B\cap X$.

For every nonisotropic $p\in\Bbb P$, define the linear map
$$R(p):x\mapsto2\frac{\langle x,p\rangle}{\langle
p,p\rangle}p-x.\eqno{\bold{(2.1.3)}}$$
Such a definition provides the reflection $R(p)\in\SU$. If $p$ is
positive, that is, if $p\notin\overline{\Bbb B}$, then $R(p)$ is the
reflection in the complex geodesic
$\Bbb Pp^\perp\cap\overline{\Bbb B}$, where
$p^\perp:=\big\{x\in V\mid\langle x,p\rangle=0\big\}$. The point $p$ is
said to be the {\it polar point\/} to $\Bbb Pp^\perp$. For brevity,
$\Bbb Pp^\perp$ will denote both the complex geodesic
$\Bbb Pp^\perp\cap\overline{\Bbb B}$ and its projective line, depending
on the context.

It is well known that every bisector in $\overline{\Bbb B}$ admits two
decompositions into totally geodesic $2$-planes: the slice
decomposition into complex geodesics and the meridional decomposition
into $\Bbb R$-planes (see [Gol, p.~152]). We treat every bisector as
being oriented by an orientation of its real spine and the natural
(complex) orientation of its slices. Let $S$ and $S'$ be two complex
geodesics. The inequality $\ta(p,p')>1$ means that $S$ and $S'$ are
ultraparallel [Gol, p.~100], where $p$ and $p'$ stand for the polar
points to $S$ and $S'$. In this case, there exists a unique bisector
$\Cal B$ containing $S$ and $S'$ as slices [Gol, p.~165,
Theorem~5.2.4]. We denote by $\B[S,S']\subset\Cal B$ the closed
oriented {\it segment of bisector\/} formed by the slices from $S$ to
$S'$ and call $\Cal B$ the {\it full bisector\/} of the segment
$\B[S,S']$. It is easy to see that there exists a unique slice $M$ of
$\B[S,S']$ such that the reflection in $M$ exchanges $S$ and $S'$. We
call $M$ the {\it middle slice\/} of $\B[S,S']$.

Let $\B[S,S']$ be an oriented segment of bisector and let $s\in S$. We
take a meridian $\R$ of $\B[S,S']$ such that $s\in\R$. Let $b$ denote
the curve contained in $\R$ that begins at $s\in S$, ends at some point
$s'\in S'$, and is equidistant from the real spine of the bisector ($b$
is a segment of a hypercycle in $\R$). We call $b$ the {\it meridional
curve\/} of $\B[S,S']$ generated by $s\in S$. When
$s\in\partial_\infty S$, the {\it ideal meridional curve\/} $b$ is
defined as the segment of the $\Bbb R$-circle in $\B[S,S']$ that
contains $s$.

\medskip

{\bf Lemma 2.1.4.} {\sl Let\/ $\B[S,S']$ be an oriented segment of
bisector and let\/ $b$ be an\/ {\rm(}ideal\/{\rm)} meridional curve
of\/ $\B[S,S']$ that begins at\/ $s\in S$ and ends at\/ $s'\in S'$.
Then\/ $Rb=b$ and\/ $Rs=s'$, where\/ $R$ stands for the reflection in
the middle slice of\/ $\B[S,S']$.}

\medskip

{\bf Proof.} There is a meridian $\R$ of $\B[S,S']$ such that
$b\subset\R$. The reflection in the middle slice $M$ of $\B[S,S']$,
being restricted to $\R$, is the reflection in the geodesic $\R\cap M$
that is orthogonal to the real spine of
$\B[S,S']$.\hfill$_\blacksquare$

\medskip

{\bf Lemma 2.1.5 {\rm[Appendix, Lemma 4.2.16]}.} {\sl An\/
{\rm(}ideal\/{\rm)} meridional curve depends continuously on the
segment of bisector and on the initial point of the curve.}

\medskip

{\bf2.1.6.}~Let $S_1,\dots, S_n$ be complex geodesics such that $S_i$
and $S_{i+1}$ are ultraparallel for all $i$ (the indices are modulo
$n$). Denote $\B_i:=\B[S_i,S_{i+1}]$. We call $(\B_1,\dots,\B_n)$ a
{\it configuration of bisectors.} For~every~$\B_i$, we denote by $M_i$
the middle slice of $\B_i$ and by $R_i$, the reflection in $M_i$.

A configuration of bisectors $(\B_1,\dots,\B_n)$ is called {\it
simple\/} if $\B_i\cap\B_j\ne\varnothing$ implies that
$j\in\{i-1,i,i+1\}$ and if $\B_{i-1}\cap\B_i=S_i$ for all $i$. A
configuration is {\it transversal\/} if, for all $i$,
$\B_{i-1}\cap\B_i=S_i$ and the full bisectors of $\B_{i-1}$ and $\B_i$
are transversal along $S_i$ (this includes transversality along
$\partial_\infty S_i$).

\medskip

{\bf Definition 2.1.7.} A {\it cycle of bisectors\/} is a configuration
of bisectors $(\B_1,\dots,\B_n)$ such that the reflections $R_i$'s in
the middle slices $M_i$'s satisfy the relation $R_n\dots R_1=1$ in
$\PU$. In terms of $\SU$, identifying the centre of $\U$ with the unit
complex numbers, we have $R_n\dots R_1=\delta$, where $\delta\in\Bbb C$
and $\delta^3=1$. We denote by $H:=\langle R_i\mid1\le i\le n\rangle$
the subgroup in $\PU$ generated by the $R_i$'s.

\medskip

Given reflections $R_1,\dots,R_n$ in complex geodesics $M_1,\dots,M_n$
subject to the relation $R_n\dots R_1=1$ in~$\PU$, a cycle of bisectors
can be constructed as follows. Let $S_1$ be a complex geodesic. Define
$S_{i+1}:=\nomathbreak R_iS_i$. Requiring that $M_i$ and $S_i$ are
ultraparallel for all $i$, we obtain the segments
$\B_i:=\B[S_i,S_{i+1}]$ forming the cycle $(\B_1,\dots,\B_n)$.
Constructing a cycle in this way does depend on the choice of $S_1$.
There are many suitable choices for $S_1$. Indeed, for some positive
points $m_i$\ and $p_1$, we have $R_i=R(m_i)$, $M_i=\Bbb Pm_i^\perp$,
$S_1=\Bbb Pp_1^\perp$, and $S_i:=\Bbb Pp_i^\perp$, where
$p_{i+1}:=R_ip_i$. In these terms, the inequality $\ta(m_i,p_i)>1$ says
that $M_i$ and $S_i$ are ultraparallel. We have
$\ta(m_i,p_i)=\ta(R_1\dots R_{i-2}R_{i-1}m_i,p_1)$. All these tances
tend to infinity while $p_1\notin\overline{\Bbb B}$ tends to a generic
point in $\partial_\infty\Bbb B$.

\smallskip

{\bf2.1.8.~Meridional curve of a cycle.} Let $(\B_1,\dots,\B_n)$ be a
configuration of bisectors. Take $s_1\in S_1$. The point $s_1$
generates a meridional curve $b_1$ of $\B_1$ which ends at some
$s_2\in S_2$. By Lemma 2.1.4, $R_1b_1=b_1$ and $R_1$ exchanges $s_1$
and $s_2$. Inductively, $s_i$ generates the meridional curve $b_i$ of
$\B_i$ which ends at some $s_{i+1}\in S_{i+1}$. Again, $R_ib_i=b_i$,
$R_is_i=s_{i+1}$, and $R_is_{i+1}=s_i$. We call\footnote{For two curves
$c_1$ and $c_2$ such that $c_2$ begins at the final point of $c_1$, we
denote by $c_1\cup c_2$ their concatenation.}
the curve $b:=b_1\cup\dots\cup b_n$ the {\it meridional curve\/} of the
configuration generated by $s_1$. Similarly,
$s_1\in\partial_\infty S_1$ generates an {\it ideal meridional curve.}
The arcs $b_1,\dots,b_n$ are {\it edges\/} and the points
$s_1,\dots,s_n$ are {\it vertices\/} of the (ideal) meridional curve.
In the case of a cycle, the meridional curves are closed due to the
relation $R_n\dots R_1=1$, as shown in the picture.

\noindent
$\hskip111pt\vcenter
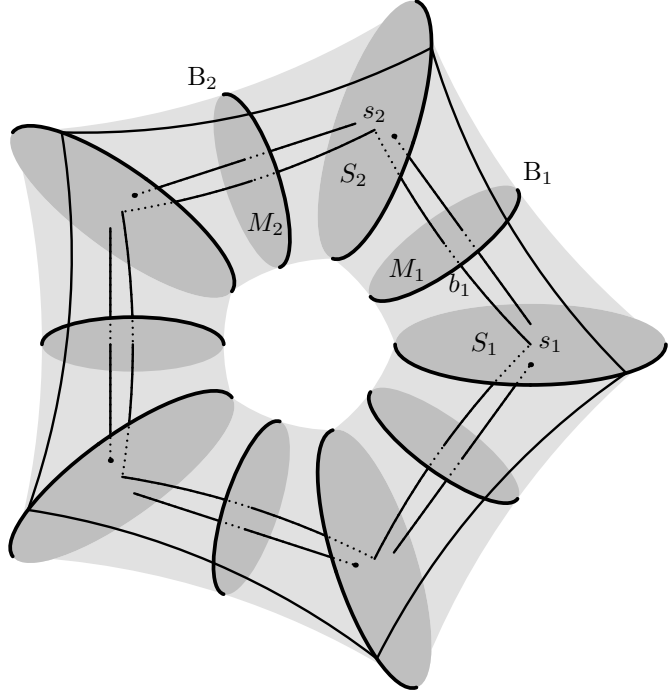
{\hbox{\epsfbox{Pictures.14}}}$

\vskip5pt

\smallskip

{\bf2.1.9.~Planar example.} The concept of a cycle of bisectors in
$\Bbb B$ originates from the following planar example.

\smallskip

Fix an integer $n\ge5$. Let $P$ be a simply connected geodesic $n$-gon
in the real hyperbolic plane $\Bbb H_\Bbb R^2$ with vertices
$v_1,\dots,v_n$ and interior angles $\alpha_1,\dots,\alpha_n$ such that
$\alpha_1+\dots+\alpha_n=2\pi$. We denote the edges of $P$ by
$e_1,\dots,e_n$, where $e_i$ connects $v_i$ and $v_{i+1}$ (the indices
are modulo $n$). Let $r_i$ stand for the reflection in the middle point
of $e_i$. By Poincar\'e's Polyhedron Theorem, $P$ provides a
fundamental domain for the group
$$H_n:=\langle r_1,\dots,r_n\mid r_n\dots r_1=1{\text{ \rm and
}}r_i^2=1\rangle\eqno{\bold{(2.1.10)}}$$
generated by the $r_i$'s subject to the defining relations
$r_n\dots r_1=1$ and $r_i^2=1$. The group $H_n$ is the fundamental
group of the sphere with $n$ elliptic points of order $2$, i.e., of the
orbifold $\Bbb S^2(2,\dots,2)$. Passing to a torsion-free subgroup of
finite index, we arrive at surface groups as follows.

For even $n$, $P\cup r_1P$ provides a fundamental region for the
subgroup
$$G_n:=\langle r_1r_i\mid2\le i\le n\rangle$$
of index $2$, implying that $G_n$ is the fundamental group of a closed
orientable surface of genus $\frac n2-1$ (the~polygon $P\cup r_1P$ has
two cycles of vertices and $n-1$ pairs of edges to identify).

For odd $n$, the polygon $P\cup r_1P\cup r_2P\cup r_2r_1P$ provides a
fundamental region for the group $T_n$ of index $4$ generated by
$$x:=r_2r_1r_n,\qquad y:=r_2r_nr_1,\qquad z:=r_2r_1r_2r_1,\qquad
x_i:=r_1r_i,\qquad y_i:=r_2r_1r_ir_2,$$
where $3\le i\le n-1$,
$$T_n:=\langle x,y,z,x_i,y_i\mid3\le i\le n-1\rangle,$$
implying that $T_n$ is the fundamental group of a closed orientable
surface of genus $n-3$.

\smallskip

{\bf2.1.11.~Toledo invariant.} Following Toledo [Tol], one can introduce
the Toledo invariant of a representation in $\PU$ of an arbitrary group
with discrete and cocompact action on a two-disc (see also [Kre]). Every
cycle of bisectors $(\B_1,\dots,\B_n)$ defines a representation of $H_n$
in $\PU$ given by $r_i\mapsto R_i$. In order to calculate the Toledo
invariant of this representation, we express a K\"ahler primitive in
coordinate-free terms. In fact, such a coordinate-free approach is kind
of a motto in this article. It makes it possible to significantly
simplify the formulae.

\smallskip

We can regard every linear map $\varphi\in\Lin_\Bbb C(\Bbb Cp,V)$ as a
tangent vector $t_{\varphi}\in\T_p\Bbb P$ by defining
$$t_{\varphi}f:=\frac d{d\varepsilon}\Big|_{\varepsilon=0}\hat
f\big(p+\varepsilon\varphi(p)\big)$$
for a local smooth function $f$ on $\Bbb P$ and its lift $\hat f$ to
$V$. In this known way, we identify the tangent space $\T_p\Bbb P$ with
the $\Bbb C$-vector space $\Lin_\Bbb C(\Bbb Cp,V/\Bbb Cp)$ of linear
maps. In particular, for every $p\notin\partial_\infty\Bbb B$,
$${\T}_p\Bbb P\simeq\langle-,p\rangle p^\perp.$$
(Of course, $\langle-,p\rangle p^\perp$ means
$\langle-,p\rangle\otimes_\Bbb Cp^\perp$.) For a nonisotropic $p\in V$
and for $v\in p^\perp$, we denote
$$v_p:=\langle-,p\rangle v\in{\T}_p\Bbb P.$$
Thus, $\T_p\Bbb P$ is equipped with the hermitian form defined by
$\langle v_p,w_p\rangle:=-\langle p,p\rangle\langle v,w\rangle$, where
$v,w\in p^\perp$. We arrive at a positive definite {\it hermitian
metric\/} on the complex hyperbolic plane $\Bbb B$ that coincides,
up~to the scale factor of $4$, with the one introduced in [Gol, p.~77]
(see Corollary 4.1.18).

\medskip

{\bf Lemma 2.1.12 {\rm[Appendix, Propositions 4.5.6 and 4.5.2]}.} {\sl
Let\/ $0\ne c\in V$. Define the\/ $1$-form\/ $P_c$ by the rule
$$P_c(v_p):=-\Im\Big(\frac{\langle p,p\rangle\langle
v,c\rangle}{2\langle p,c\rangle}\Big)$$
for all \/ $p\notin\Bbb Pc^\perp\cup\partial_\infty\Bbb B$. Let\/
$c_1,c_2\in V$, $\langle c_1,c_2\rangle\ne0$. Define the function\/
$f_{c_1,c_2}$ by the rule\footnote{In this lemma, it is better to read
$\Arg$ as a multi-valued function on $\Bbb C\setminus\{0\}$.}
$$f_{c_1,c_2}(p):=\frac12\Arg\Big(\frac{\langle c_1,p\rangle\langle
p,c_2\rangle}{\langle c_1,c_2\rangle}\Big)$$
for all\/ $p\notin\Bbb Pc_1^\perp\cup\Bbb Pc_2^\perp$. Then\/
$dP_c=\omega$ and\/ $P_{c_1}-P_{c_2}=df_{c_1,c_2}$, where, for\/
$p\in\Bbb B$ and\/ $v_p,w_p\in\T_p\Bbb B$, the~K\"ahler form on\/
$\Bbb B$ is given by\/ $\omega(v_p,w_p):=\Im\langle v_p,w_p\rangle$.}

\medskip

{\bf Lemma 2.1.13 {\rm[Appendix, Corollary 4.2.3]}.} {\sl Let\/ $\R$ be
an\/ $\Bbb R$-plane and let\/ $p,c\in\R$. Then\/ $P_c(v_p)=0$ for
every\/ $v_p\in\T_p\R$.}

\medskip

{\bf Lemma 2.1.14 {\rm[Appendix, Lemma 4.1.16]}.} {\sl If\/
$p,g_1,g_2\in\Bbb B$, then\/
$\displaystyle\frac{\langle g_1,p\rangle\langle p,g_2\rangle}{\langle
g_1,g_2\rangle}$
cannot be real nonnegative.}

\medskip

We are now able to calculate the Toledo invariant.

\medskip

{\bf Proposition 2.1.15.} {\sl Let\/ $(\B_1,\dots,\B_n)$ be a cycle of
bisectors and let\/ $\varrho:H_n\to\PU$ be the representation given
by\/ $r_i\mapsto R_i$. Then the Toledo invariant of\/ $\varrho$
satisfies\/ $\tau\equiv n-\displaystyle\frac{\Arg\delta}\pi\mod2$,
where\/ $\delta$ is introduced in Definition\/ {\rm2.1.7}.}

\medskip

{\bf Proof.} We have $R_i=R(m_i)$, where $m_i$ stands for the polar
point to the middle slice $M_i$. Assuming that
$\langle m_i,m_i\rangle=1$, we obtain
$R_i:x\mapsto2\langle x,m_i\rangle m_i-x$ by (2.1.3). Let
$b=b_1\cup\dots\cup b_n$ be a meridional curve of $(\B_1,\dots,\B_n)$
with vertices $s_1,\dots,s_n$. Fixing a representative $s_1\in V$,
define $s_{i+1}=R_is_i\in V$. In particular, $s_{n+1}=\delta s_1$.
Since
$\langle R_ix,x\rangle=2\langle x,m_i\rangle\langle
m_i,x\rangle-\langle x,x\rangle>0$
for every $x\in\Bbb B$, we obtain $\langle s_{i+1},s_i\rangle>0$.

Let $D\subset\Bbb B$ be any disc with $\partial D=b$. In a standard
way, we define a $\varrho$-equivariant continuous map
$\varphi:\Bbb H_\Bbb R^2\to\Bbb B$ such that $\varphi(P)=D$,
$\varphi(v_i)=s_i$, and $\varphi(e_i)=b_i$ (see 2.1.9). The Toledo
invariant of $\varrho$ is defined as
$$\tau:=\frac4{2\pi}\int_{P}\varphi^*\omega$$
(the hermitian metric we introduced differs by the factor of $4$ from
the one in [Tol] (see Corollary 4.1.18)). Taking $c\in\Bbb B$ and
applying Lemma 2.1.12, we obtain
$$\tau=\frac2\pi\int_{D}\omega=\frac2\pi\int_{\partial
D}P_c=\frac2\pi\sum_i\int_{b_i}P_c.$$
By Lemmas 2.1.13 and 2.1.12,
$$\int_{b_i}P_c=\int_{b_i}(P_c-P_{s_i})=\int_{b_i}df_{c,s_i}.$$
This number is the total variation of
$\displaystyle\frac12\Arg\frac{\langle c,p\rangle\langle
p,s_i\rangle}{\langle c,s_i\rangle}$
while $p$ runs over $b_i$ from $s_i$ to $s_{i+1}$. By~Lem\-ma~2.1.14,
$\displaystyle\frac{\langle c,p\rangle\langle p,s_i\rangle}{\langle
c,s_i\rangle}$
cannot be real nonnegative. It follows that
$$\int_{b_i}df_{c,s_i}=\frac12\Arg\frac{\langle c,s_{i+1}\rangle\langle
s_{i+1},s_i\rangle}{\langle c,s_i\rangle}-\frac12\Arg\frac{\langle
c,s_i\rangle\langle s_i,s_i\rangle}{\langle
c,s_i\rangle}=\frac12\Arg\frac{\langle c,s_{i+1}\rangle}{\langle
c,s_i\rangle}-\frac\pi2\eqno{\bold{(2.1.16)}}$$
since $\langle s_{i+1},s_i\rangle>0$ and $\langle s_i,s_i\rangle<0$.
Calculating $\!\!{\mod2}$, we get
$$\tau\equiv\frac1\pi\sum_i\big(\Arg\langle
c,s_{i+1}\rangle-\Arg\langle c,s_i\rangle-\pi\big)\equiv$$
$$\equiv\frac1\pi\big(-n\pi+\Arg\langle c,s_{n+1}\rangle-\Arg\langle
c,s_1\rangle\big)\equiv-n-\frac{\Arg\delta}\pi\mod2.
\eqno{_\blacksquare}$$

Considering the surface groups $G_n$ or $T_n$ of finite index in $H_n$
(see 2.1.9), we obtain the following

\medskip

{\bf Corollary 2.1.17.} {\sl For even\/ $n$, the Toledo invariant of
the induced representation of\/ $G_n$ in\/ $\PU$ satisfies\/
$\tau\equiv-\displaystyle\frac{2\Arg\delta}\pi\mod4$. For odd\/ $n$,
the Toledo invariant of the induced representation of\/ $T_n$ in\/
$\PU$ satisfies\/
$\tau\equiv4n-\displaystyle\frac{4\Arg\delta}\pi\mod8$, where\/
$\delta$ is introduced in Definition\/ {\rm2.1.7}.}

\medskip

{\bf Proof.} When passing to a finite cover, the Toledo invariant gets
multiplied by the degree of the cover. In our case, this means that the
first integral in Proposition 2.1.15 should be taken over two or four
copies of a fundamental domain for $H_n$.\hfill$_\blacksquare$

\medskip

\centerline{\bf2.2.~Discreteness}

\medskip

We intend to prove discreteness of the group $H$ generated by
reflections in middle slices of a cycle of bisectors by showing that a
polyhedron bounded by the cycle is fundamental. To this end, we
reformulate [AGr3, Theorem 3.5] in terms of a simple transversal
configuration of bisectors.

All isometries in this subsection are considered as belonging to $\PU$.

\smallskip

{\bf2.2.1.}~Given oriented bisectors $\B$ and $\B'$ with a common slice
$S$, denote by $n$ and $n'$ the normal vectors at $s\in S\cap\Bbb B$ to
$\B$ and to $\B'$, respectively. Since both $n$ and $n'$ are tangent to
the naturally oriented complex geodesic passing through $s$ and
orthogonal to $S$, it makes sense to measure the oriented angle from
$n$ to $n'$. This angle is said to be the {\it oriented angle\/} from
$\B$ to $\B'$ at $s$.

\medskip

{\bf Definition 2.2.2.} A simple configuration of bisectors
$(\B_1,\dots,\B_n)$ gives rise to the closed solid torus
$\partial_0P:=\B_1\cup\dots\cup\B_n$ {\it fibred by slices\/} that
divides $\overline{\Bbb B}$ into two closed connected parts. The part
$P$ situated on the side of the normal vector to each $\B_i$ is called
the {\it polyhedron\/} of the configuration. The torus
$T:=\partial_\infty\partial_0P$ divides the boundary of $P$ into two
closed parts, $\partial P=\partial_0P\cup\partial_\infty P$.
So,~$T=\partial\partial_\infty P=\partial_0P\cap\partial_\infty P$. The
ideal meridional curves are contained in $T$. In the particular case of
a cycle of bisectors, $\partial_0P$ is also {\it fibred by meridional
curves\/} (hence, $T$ is fibred by the ideal meridional curves).

\medskip

{\bf2.2.3.}~Let $(\B_1,\dots,\B_n)$ be a simple configuration of
bisectors with common slices $S_i$'s, where $\B_i:=\B[S_i,S_{i+1}]$
(the indices are modulo $n$). A {\it face-paring\/} of the polyhedron
$P$ of the configuration is an involution
$\overline{\phantom{m}}:\{1,\dots,n\}\to\{1,\dots,n\}$ and a family of
isometries $I_i\in\PU$ satisfying $I_i\B_i=\B_{\overline i}^-$ and
$I_iI_{\overline i}=1$ for all $i$, where $\B^-$ denotes the segment
$\B$ taken with the opposite orientation.

Define a {\it geometric cycle\/} of common slices as follows. Start
from an ordered triple $(\B_{\overline i_0},S_{l_0},\B_{i_1})$, where
$S_{l_0}$ is a common slice of two different segments
$\B_{\overline i_0},\B_{i_1}$. (In other words, either
$i_1=\overline i_0+1$ with $l_0=i_1$ or $\overline i_0=i_1+1$ with
$l_0=\overline i_0$.) Applying $I_{i_1}$ to $\B_{i_1}$ and to
$S_{l_0}$, we obtain a new triple
$(\B_{\overline i_1},S_{l_1},\B_{i_2})$, where
$S_{l_1}:=I_{i_1}S_{l_0}$ and $\B_{i_2}$ is uniquely determined by the
requirement that $S_{l_1}$ is a common slice of two different segments
$\B_{\overline i_1},\B_{i_2}$. Then we apply $I_{s_2}$ to $\B_{i_2}$
and to $S_{l_1}$, obtain $(\B_{\overline i_2},S_{l_2},\B_{i_3})$, and
so on. When we arrive back at
$(\B_{\overline i_k},S_{l_k},\B_{i_{k+1}})=(\B_{\overline
i_0},S_{l_0},\B_{i_1})$
(not necessarily for the first time), we obtain a {\it cycle\/} of
common slices and its {\it isometry\/} $I:=I_{i_k}\dots I_{i_1}$. (The
least possible $k$ provides a {\it combinatorial\/} cycle.) The cycle
is {\it geometric\/} if its isometry $I$ is the identity and $k$ is
minimal with $I=1$. In this case, $I_{i_k}\dots I_{i_1}=1$ is the {\it
relation\/} of the geometric cycle.

Consider a geometric cycle as above. Pick a point
$s_0\in S_{l_0}\cap\Bbb B$, define
$s_j:=I_{i_j}\dots I_{i_1}s_0\in S_{l_j}$, and~denote by $\alpha_j$ the
oriented angle from $\B_{l_j}$ to $\B_{l_j-1}^-$ at $s_j$. The sum
$\alpha:=\alpha_1+\dots+\alpha_k$ is the {\it total angle\/} of the
geometric cycle at $s_0$.

Strong Simplicity IV [AGr3, Section 3] follows immediately from the
fact that $\B_i$ and $\B_j$ are disjoint in $\overline{\Bbb B}$ for
$j\notin\{i-1,i,i+1\}$ because this implies the metric separability of
$\B_i\cap\Bbb B$ and $\B_j\cap\Bbb B$. Condition~(1) of [AGr3, Theorem
3.5] (requiring that the face-pairing isometries send the interior of
$P$ into the exterior of $P$) follows directly from
$I_i\B_i=\B_{\overline i}^-$. As in [AGr3], we denote by
$N(X,\varepsilon)$ the $\varepsilon$-neighbourhood of $X\subset\Bbb B$.
Condition (3) of [AGr3, Theorem 3.5] for a simple transversal
configuration of bisectors is verified in the following

\medskip

{\bf Lemma 2.2.4 {\rm[Appendix, Lemma 4.3.6]}.} {\sl Let\/ $\B$ and\/
$\B'$ be full bisectors transversal along a common slice\/ $S$. Then,
for every\/ $\vartheta>0$, there exists\/ $\varepsilon>0$ such that\/
$\B'\cap N(\B\cap\Bbb B,\varepsilon)\subset N(S\cap\Bbb B,\vartheta)$.}

\medskip

Now we can reformulate [AGr3, Theorem 3.5] for a simple transversal
configuration of bisectors.

\medskip

{\bf Theorem 2.2.5.} {\sl Let\/ $(\B_1,\dots,\B_n)$ be a simple
transversal configuration of bisectors whose polyhedron\/ $P$ is
equipped with a face-pairing such that every common slice\/ $S_i$
belongs to some geometric cycle. Suppose that, for every geometric
cycle of common slices, there exists a point in a common slice of the
cycle at which the total angle of the cycle equals\/ $2\pi$. Then\/ $P$
is a fundamental polyhedron for the group\/ $G$ generated by the
face-pairing isometries and the defining relations of\/ $G$ are the
relations of the geometric cycles and the relations\/
$I_iI_{\overline i}=1$.}\hfill$_\blacksquare$

\medskip

When dealing with a transversal cycle of bisectors, we define
$\overline i:=i$ and $I_i:=R_i$. In this case, we~obtain a unique
geometric cycle.

\medskip

{\bf Corollary 2.2.6.} {\sl Let\/ $(\B_1,\dots,\B_n)$ be a simple
transversal cycle of bisectors with total angle\/ $2\pi$ at some point
and let\/ $P$ denote the polyhedron of the cycle. Then the subgroup\/
$H\subset\PU$ generated by the reflections in the middle slices of
the\/ $\B_i$'s is discrete and\/ $P$ is a fundamental region in\/
$\Bbb B$ for\/ $H$. The~group\/ $H$ is isomorphic to\/ $H_n$.}
\hfill$_\blacksquare$

\medskip

\centerline{\bf2.3.~Simplicity and transversality}

\medskip

In order to apply Corollary 2.2.6, we need a criterion for the
simplicity of a cycle. We have found no adequate numerical criterion of
this sort since it is rather difficult to decide whether two segments
of bisectors intersect or not (see, for instance, [San]). Fortunately,
by appealing to a simple kind of Poincar\'e duality, the simplicity of
a cycle can be inferred from some extra transversalities (see
Criterion~2.3.9).

All isometries in this subsection are considered as belonging to $\PU$.

\smallskip

{\bf2.3.1.~Notation concerning bisectors.} Every two distinct points
$p_1,p_2\in\overline{\Bbb B}$ determine a closed oriented geodesic
segment $\G[p_1,p_2]$ from $p_1$ to $p_2$ and the corresponding closed
oriented segment of bisector $\B[p_1,p_2]$ whose real spine is
$\G[p_1,p_2]$. (Note that we also use to specify a closed oriented
segment of bisector in terms of two ultraparallel complex geodesics as
in 2.1.1.) For an oriented segment of bisector $\B:=\B[p_1,p_2]$, we
denote by ${\prec}\B{\succ}$ the full bisector of the segment, endowed
with the corresponding orientation. Hence, ${\prec}\B{\succ}=\B[u,v]$,
where $u$ and $v$ are the vertices of ${\prec}\B{\succ}$ taken in a
suitable order. We also denote $\B{\succ}:=\B[p_1,v]$ and
$\B^-:=\B[p_2,p_1]$ (so, $\B^-$ is $\B$ taken with the opposite
orientation). In this way, we deal with three kinds of segments of
bisectors: finite, semifinite, and infinite (full bisectors).

\vskip8pt

\noindent
$\vcenter{\hbox{\epsfbox{Picture.3}}}$

\vskip-119pt

\noindent
$\hskip357pt\vcenter{\hbox{\epsfbox{Picture.4}}}$

\leftskip90pt

\rightskip105pt

\vskip-94pt

Every full oriented bisector $\B$ divides $\overline{\Bbb B}$ into two
{\it half-spaces\/} (closed $4$-balls) $K^+$ and $K^-$, where $K^+$
lies on the side of the normal vector to $\B$. Note that, in general,
the~fundamental polyhedron $P$ in Corollary 2.2.6 is not contained in
the half-space $K_i^+$ {\it bounded\/} by ${\prec}\B_i{\succ}$.

\medskip

{\bf Lemma 2.3.2 {\rm[Appendix, Lemma 4.3.5]}.} {\sl Let\/ $\B$ and\/
$\B'$ be two full bisectors transversal along a common slice\/ $S$.
Then\/ $\B$, $\B'$, and\/ $\partial_\infty\Bbb B$ are transversal and\/
$\B\cap\B'=S$.}

\medskip

Let $\B$ and $\B'$ be two segments of bisectors that begin\break

\leftskip0pt

\vskip-25pt

\noindent
$\hskip352pt\vcenter{\hbox{\epsfbox{Picture.5}}}$

\vskip-57pt

\noindent
with a common slice $S$. Suppose that ${\prec}\B{\succ}$ and
${\prec}\B'{\succ}$ are transversal along~$S$. Then, by Lemma 2.3.2,
${\prec}\B{\succ}$ and ${\prec}\B'{\succ}$ divide $\overline{\Bbb B}$
into four closed $4$-balls. We~define the {\it sector\/} $A$ from
$\B{\succ}$ to $\B'{\succ}$ as follows. If the oriented angle from
$\B{\succ}$ to $\B'{\succ}$ at a point in $S$ does not exceed $\pi$ (it
does not matter at which point in $S$ we measure the angle), then we
put $A:=K^+\cap{K'}^-$, where $K^+$ and ${K'}^+$\break

\vskip-12pt

\rightskip0pt

\noindent
stand respectively for the half-spaces bounded by ${\prec}\B{\succ}$
and ${\prec}\B'{\succ}$. Otherwise, we put $A:=K^+\cup{K'}^-$.

\smallskip

{\bf2.3.3.~Transversal adjacency.} Let $C,M_1,M_2$ be pairwise
ultraparallel complex geodesics. We~connect them by oriented segments
of bisectors $\B_1:=\B[C,M_1]$, $\B:=\B[M_1,M_2]$, and
$\B_2:=\B[M_2,C]$ thus obtaining an {\it oriented triangle\/} of
bisectors $\Delta(C,M_1,M_2)$. The triangle is said to be {\it
transversal\/} if the corresponding full bisectors are transversal
along their common slices. Since ${\prec}\B_1{\succ}$ and
${\prec}\B_2{\succ}$ are transversal along $C$, either the oriented
angle from $\B_1{\succ}$ to $\B_2^-{\succ}$ is less than $\pi$ or the
oriented angle from $\B_2^-{\succ}$ to $\B_1{\succ}$ is less than
$\pi$. We denote by $D$ the sector containing the smaller angle and
call it the {\it interior sector\/} at $C$. So, either
$D=K_1^+\cap K_2^+$ or $D=K_1^-\cap K_2^-$, where $K_i^+$ and $K^+$
denote the half-spaces bounded by ${\prec}\B_i{\succ}$ and
${\prec}\B{\succ}$. Similarly, we define the interior sector $D_i$
between $\B$ and $\B_i$ at $M_i$. By altering the orientation of the
triangle, we can always assume that $D=K_1^+\cap K_2^+$.

Given one more complex geodesic $S$ ultraparallel to $M_1$ and to $M_2$
such that the triangle $\Delta(S,M_2,M_1)$ is also transversal, denote
$\B'_1:=\B[M_1,S]$ and $\B'_2:=\B[S,M_2]$. We say that the triangle
$\Delta(S,M_2,M_1)$ is {\it transversally adjacent\/} to the triangle
$\Delta(C,M_1,M_2)$ if some point in $S$ belongs to $D$ and if the
bisectors ${\prec}\B_i{\succ}$ and ${\prec}\B'_i{\succ}$ are
transversal along $M_i$ for $i=1,2$.

\noindent
$\hskip362pt\vcenter{\hbox{\epsfbox{Picture.10}}}$

\vskip-80pt

\rightskip100pt

\medskip

{\bf Lemma 2.3.4.} {\sl Let\/ $\Delta(C,M_1,M_2)$ be a transversal
triangle of bisectors oriented so that\/ $D=K_1^+\cap K_2^+$, where\/
$D$ stands for the interior sector at\/ $C$. Then the segment\/
$\B:=\B[M_1,M_2]$ is contained in\/ $D$. Moreover, for the interior
sector\/ $D_i$ at\/ $M_i$, we have\/ $D_i=K^+\cap K_i^+$.

Suppose that some transversal triangle\/ $\Delta(S,M_2,M_1)$ is
transversally adjacent to\/ $\Delta(C,M_1,M_2)$. Then\footnote{Note
that $\Delta(S,M_2,M_1)$ can be `inside' or `outside'
$\Delta(C,M_1,M_2)$ (see also Remark 2.3.6).}
$\Delta(S,M_2,M_1)\subset D$.}

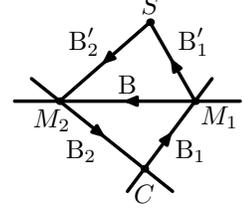
\rightskip0pt

\medskip

{\bf Proof.} Since $\B$ is connected and intersects
${\prec}\B_i{\succ}$ only in $M_i$, the segment $\B$ is contained in
one of the four sectors formed by the ${\prec}\B_i{\succ}$'s. The only
sector that contains both $M_i$'s is $D$. This proves that
$\B\subset D$.

Suppose that $D_1=K^-\cap K_1^-$. Then, by the above statement, we have
$\B_2\subset D_1$, which implies $\B_2\subset\nomathbreak K_1^-$. On
the other hand, $\B_2\subset D\subset K_1^+$. A contradiction. The same
works for $D_2$. So,~$D_i=K^+\cap K_i^+$.

The bisector ${\prec}\B'_1{\succ}$ intersects ${\prec}\B_1{\succ}$ only
in $M_1$. Hence, either $\B'_1{\succ}\subset K_1^+$ or
$\B'_1{\succ}\subset K_1^-$. Since some point of $S$ belongs to
$D\subset K_1^+$ and $S\subset\B'_1{\succ}$, we conclude that
$\B'_1{\succ}\subset K_1^+$. By the same reason,
${\B'_2}^-{\succ}\subset K_2^+$. From $\B\subset D$, we derive that
$\B{\succ}\subset K_1^+$ and $\B^-{\succ}\subset K_2^+$. Consequently,
$D'_1\subset K_1^+$ and $D'_2\subset K_2^+$, where $D'_i$ stands for
the interior sector of $\Delta(S,M_2,M_1)$ at $M_i$. Since
$\B'_2\subset D'_1$ and $\B'_1\subset D'_2$, we obtain
$\B'_2\subset K_1^+$ and $\B'_1\subset K_2^+$. Now, from
$\B'_1{\succ}\subset K_1^+$ and ${\B'_2}^-{\succ}\subset K_2^+$, we
deduce that $\B'_i\subset D$. Therefore,
$\Delta(S,M_2,M_1)\subset D$.\hfill$_\blacksquare$

\medskip

{\bf Definition 2.3.5.} The orientation of the transversal triangle of
bisectors $\Delta(C,M_1,M_2)$ dealt with in Lemma 2.3.4 is said to be
{\it counterclockwise.}

\medskip

{\bf Remark 2.3.6.} {\sl Let\/ $\Delta(C,M_1,M_2)$ and\/
$\Delta(S,M_2,M_1)$ be counterclockwise-oriented transversal triangles.
Suppose that\/ $\Delta(S,M_2,M_1)$ is transversally adjacent to\/
$\Delta(C,M_1,M_2)$. Then\/ $\Delta(S,M_2,M_1)$ and\/
$\Delta(C,M_1,M_2)$ lie in the distinct half-spaces bounded by\/
${\prec}\B[M_1,M_2]{\succ}$.}\hfill$_\blacksquare$

\medskip

{\bf2.3.7.~Centre of a configuration.} Let $(\B_1,\dots,\B_n)$ be a
configuration of bisectors and let $C$ be a complex geodesic
ultraparallel to every middle slice $M_i$. Denote $\B'_i:=\B[C,M_i]$.
If, for every $i$, the full bisectors ${\prec}\B'_{i-1}{\succ}$ and
${\prec}\B'_i{\succ}$ are transversal along their common slice $C$,
then we call $C$ a {\it centre\/} of the configuration. Take $c\in C$.
Denoting by $\beta_i$ the oriented angle from $\B'_{i-1}{\succ}$ to
$\B'_i{\succ}$ at $c$, it is easy to see that
$\beta:=\beta_1+\dots+\beta_n$, the {\it central angle\/} of the
configuration at $C$, does not depend on the choice\break

\vskip-6pt

\noindent
$\vcenter{\hbox{\epsfbox{Picture.8}}}$

\leftskip67pt

\vskip-62pt

\noindent
of $c\in C$ and is an integer multiple of $2\pi$. The condition that a
configuration possesses a centre does not seem too restrictive.

\medskip

{\bf Lemma 2.3.8.} {\sl Let\/ $\B'_i$, $1\le i\le n$, be segments of
bisectors. Suppose that they all have a common slice\/ $C$ {\rm(}they
all begin with\/ $C${\rm)}, that\/ ${\prec}\B'_{i-1}{\succ}$ and\/
${\prec}\B'_i{\succ}$ are transversal along\/ $C$ for all\/ $i$, and
that, for some\/ $c\in C\cap\Bbb B$, the sum of the angles from\/
$\B'_{i-1}{\succ}$ to\/ $\B'_i{\succ}$}\break

\vskip-12pt

\leftskip0pt

\noindent
{\sl at\/ $c$ equals\/ $2\pi$. Then the\/ $\B'_i{\succ}$'s divide\/
$\overline{\Bbb B}$ into\/ $n$ sectors\/ $A'_i$ {\rm(}from\/
$\B'_{i-1}{\succ}$ to\/ $\B'_i{\succ}${\rm)} such that\/
$A'_i\cap A'_j=C$ if\/ $j\notin\{i-1,i,i+1\}$ and\/
$A'_i\cap A'_{i+1}=\B'_i{\succ}$. In other words, there are no extra
intersections besides the obvious ones.}

\medskip

{\bf Proof.} Since $C\cap\Bbb B$ is connected, the mentioned sum of
angles equals $2\pi$ at every point in $C\cap\Bbb B$.
Let~$p\in A'_i\cap A'_j$ be a point that does not belong to the
intersections of $A'_i$ and $A'_j$ listed above. Take the complex
geodesic $S$ passing though $p$ and orthogonal to $C$. By Lemma 2.3.2,
inside $S$, every pairwise intersection of the $\B'_i{\succ}$'s is the
point $S\cap C$. Hence, inside $S$, the sectors $A'_i$'s can only
intersect in the obvious way. A contradiction.\hfill$_\blacksquare$

\medskip

If a configuration is transversal, then $M_{i-1}$ and $M_i$ are
ultraparallel by Lemma 2.3.2. In this case, assuming that the
configuration possesses a centre $C$, we obtain the oriented triangles
$\Delta_i:=\Delta(C,M_{i-1},M_i)$ and
$\Delta'_i:=\Delta(S_i,M_i,M_{i-1})$.

To a certain extent, transversality implies simplicity:

\medskip

{\bf Criterion 2.3.9.} {\sl Let\/ $(\B_1,\dots,\B_n)$ be a transversal
configuration of bisectors possessing a centre\/ $C$. Suppose that the
triangles\/ $\Delta_i$ and\/ $\Delta'_i$ are transversal and that\/
$\Delta'_i$ is transversally adjacent to\/ $\Delta_i$ for every\/ $i$.
If the central angle at some point is\/ $2\pi$ and the angle from\/
$\B'_{i-1}{\succ}$ to\/ $\B'_i{\succ}$ does not\footnote{This condition
is not too restrictive, since only one of these angles can exceed
$\pi$.}
exceed\/ $\pi$ for every\/ $i$, then the cycle is simple.}

\medskip

{\bf Proof.} By Lemma 2.3.8, the segments $\B'_i{\succ}$ divide
$\overline{\Bbb B}$ into $n$ sectors $A'_i$'s whose interiors are
disjoint. The triangle $\Delta'_i$ is transversally adjacent to the
triangle $\Delta_i$, therefore, $\B_{i-1}\cap\B'_{i-1}{\succ}=M_{i-1}$
and $\B_i\cap\B'_i{\succ}=M_i$. By Lemma 2.3.4,
$\B[M_{i-1},S_i]\cup\B[S_i,M_i]\subset A'_i$.\hfill$_\blacksquare$

\medskip

\centerline{\bf2.4.~Fibred polyhedra. Euler number}

\medskip

Every complex hyperbolic disc bundle in the explicit series of examples
we construct in Section 3 is glued from fibred polyhedra. The Euler
number of the resulting bundles will be derived from `fractional Euler
quantities' of these parts.

We work in the PL category.

\smallskip

{\bf2.4.1.~Fibred polyhedra.} Let $P$ be the polyhedron of a simple
configuration of bisectors (see~Definition 2.2.2). In the situation
described in the next lemma, $P$ is said to be {\it fibred.}

\medskip

{\bf Lemma 2.4.2.} {\sl Suppose that the polyhedron\/ $P$ of a simple
configuration of bisectors is a closed\/ $4$-ball, $P\simeq\bold B^4$.
Then\/ $P\simeq\bold B^2\times\bold B^2$ is a bundle with\/
$\partial_0P\simeq\Bbb S^1\times\bold B^2$ being a subbundle\/ if and
only if\/ $\partial_\infty P$ is a solid torus. In this case, the slice
bundle structure of\/ $\partial_0P$ is extendable to\/ $P$.}

\medskip

{\bf Proof.} The arguments are standard. If
$P\simeq\bold B^2\times\bold B^2$, then $P\simeq\bold B^4$ and
$\partial P$ is a sphere $\Bbb S^3$ decomposed into two solid tori
$\partial\bold B^2\times\bold B^2$ and
$\bold B^2\times\partial\bold B^2$ glued along the torus
$\partial\bold B^2\times\partial\bold B^2$. Hence,
$\partial_0P\simeq\partial\bold B^2\times\bold B^2$ implies that
$\partial_\infty P\simeq\bold B^2\times\partial\bold B^2$.

Conversely, if $\partial_\infty P$ is a solid torus, then
$\partial P\simeq\Bbb S^3$ is decomposed into two solid tori glued
along the torus $T$. As is well known,\footnote{For example, one can
use a particular case of Waldhausen's theorem [Sch].}
such a decomposition of $\Bbb S^3$ is topologically unique. Arbitrarily
fibering one of the solid tori by discs (such a fibration is isotopic
to a standard one), we fiber $T$ by circles. Extending the latter
circle bundle, we fiber the other solid torus by circles and obtain the
compatible decompositions $T\simeq\Bbb S^1\times\Bbb S^1$,
$\partial_0P\simeq\Bbb S^1\times\bold B^2$, and
$\partial_\infty P\simeq\bold B^2\times\Bbb S^1$. Since
$P\simeq\bold B^4$ is a cone over $\partial P\simeq\Bbb S^3$, we can
readily extend these decompositions to a compatible decomposition
$P\simeq\bold B^2\times\bold B^2$.\hfill$_\blacksquare$

\medskip

The Dehn Lemma immediately implies the

\medskip

{\bf Remark 2.4.3.} {\sl Let\/ $P$ be the polyhedron of a simple
configuration of bisectors. Then\/ $P$ is fibred\/ {\rm(}hence,
$P\simeq\bold B^4${\rm)} if and only if there exists some simple closed
curve\/ $c\subset T$ contractible in\/ $\partial_\infty P$ such that\/
$c$ intersects each slice of\/ $\partial_0P$ exactly once, i.e., $[c]$
generates\/ $\pi_1(\partial_0P)$.}\hfill$_\blacksquare$

\medskip

The curve $c$ in Remark 2.4.3 is called {\it trivializing.}

\smallskip

\noindent
$\hskip233pt\vcenter{\hbox{\epsfbox{Picture.11}}}$

\rightskip230pt

\vskip-199pt

{\bf2.4.4.~Gluing fibred polyhedra.} Let
$(\B_1,\dots,\allowmathbreak\B_i,\B_{i+1},\dots,\B_n)$, $1\le i<n$, and
$(\B'_m,\dots,\B'_{i+1},\allowmathbreak\B_i^-,\dots,\B_1^-)$,
$1\le i<m$, be simple configurations of bisectors with a common
sequence of bisectors (the segments common to both configurations have
opposite orientations) such that the polyhedra $P_1$ and $P_2$ of the
configurations intersect only in $\B_1\cup\dots\cup\B_i$. Then we can
glue $P_1$ and $P_2$ along $\B_1\cup\dots\cup\B_i$, obtaining the {\it
gluing\/} $P:=P_1\cup P_2$ which is the polyhedron of the simple
configuration
$(\B_{i+1},\dots,\B_n,\allowmathbreak\B'_m,\dots,\B'_{i+1})$. The
polyhedra $P$ and $R_jP$ in Corollary 2.2.6 yield an example of such a
gluing: they are glued along $\B_j$. Suppose that $P_1$ and $P_2$ are
fibred. Then the solid tori $\partial_\infty P_1$ and
$\partial_\infty P_2$ intersect in the annulus
$\partial_\infty\B_1\cup\dots\cup\partial_\infty\B_i$ which is an
annular neighbourhood of a simple curve generating the fundamental
group of each\break

\vskip-12pt

\rightskip0pt

\noindent
solid torus $\partial_\infty P_k$, $k=1,2$. Hence, we can choose a
trivializing curve $c_k$ contractible in $\partial_\infty P_k$ so that
$c_1$ and $c_2$ coincide along the annulus. Thus, we arrive at

\medskip

{\bf Remark 2.4.5.} {\sl Let\/ $P_1$ and\/ $P_2$ be fibred polyhedra
and let\/ $P:=P_1\cup P_2$ be their gluing. Then\/ $P$ is fibred and a
trivializing curve of\/ $P$ can be obtained by\/ {\rm gluing}\/
{\rm(}and removing the common part from\/{\rm)} trivializing curves\/
$c_1$ of\/ $P_1$ and\/ $c_2$ of\/ $P_2$ which coincide along the common
sequence of bisectors.}\hfill$_\blacksquare$

\medskip

{\bf2.4.6.~Euler number.} Let $P$ be a fibred polyhedron of a simple
cycle, equipped with the face-pairing given by reflections in middle
slices, and let $b$ be a meridional curve of the cycle. Clearly, there
exists a simple closed disc $D\subset P$ that intersects each fibre in
$P$ exactly once and such that $b=\partial D$ and
$\overset{\,_\circ}\to D\subset\overset{\,_\circ}\to P$ (where
$\overset{\,_\circ}\to X$ denotes the interior of $X$). We assume that
meridional curves and trivializing curves are {\it oriented\/} with
respect to the orientation of the cycle. We {\it orient\/} the twice
fibred torus $T$ as follows: the first coordinate is the naturally
oriented boundary of a slice and the second is an ideal meridional
curve, already oriented. The orientation of $b$ {\it orients\/} $D$. We
call $D$ an {\it equivariant section\/} of the fibred polyhedron $P$
with respect to reflections in the middle slices. For another
meridional curve $b'$ (clearly, $b\cap b'=\varnothing$), we can find an
equivariant section $D'$ with $\partial D'=b'$ and choose $D$ and $D'$
to be transversal. We call the number $e_P:=\#D\cap D'$ (taking signs
into account) the {\it Euler number\/}\footnote{Standard arguments show
that this number does not depend on the choice of the meridional curves
$b$ and $b'$ and of the equivariant sections $D$ and $D'$.}
of the polyhedron $P$ equipped with its face-pairing (see also 2.4.8).

Clearly, $e_P$ measures the difference between two identifications of
the slice bundle $\partial_0P$ with the product
$\Bbb S^1\times\bold B^2$ : the one given by meridional curves and the
one induced by the trivialization $P\simeq\bold B^2\times\bold B^2$ (in
the torus $T$, fibred by the boundaries of slices, the trivialization
$P\simeq\bold B^2\times\bold B^2$ induces a bundle of trivializing
curves). Obviously, such a difference can be measured in terms of $T$ :

\medskip

{\bf Remark 2.4.7.} {\sl Let\/ $P$ be a fibred polyhedron of a simple
cycle. Then\/ $e_P=\#b\cap c$, where\/ $b$ stands for an ideal
meridional curve and\/ $c$, for a trivializing curve. In other words,
$[b]=e_P\cdot[s]$ in the group\/ $H_1(\partial_\infty P,\Bbb Z)$
generated by\/ $[s]$, where\/ $[s]$ is represented by the boundary of a
naturally oriented slice.}\hfill$_\blacksquare$

\medskip

{\bf2.4.8.}~Let $P$ be a fibred polyhedron of a cycle of bisectors that
satisfies the conditions of Corollary~2.2.6. Then the Euler number of
the complex hyperbolic disc bundle associated to the corresponding
discrete representation of $G_n$ or $T_n$ equals respectively $2e_P$ or
$4e_P$ (see 2.1.9) since the Euler number of a polyhedron is additive
with respect to gluing and the face-pairing respects the meridional
curves. Actually, $e_P$ is the Euler number of the orbifold disc bundle
over $\Bbb S^2(2,\dots,2)$ (see [BSi] for a definition). The Euler
number of an orbifold disc bundle is known to be additive under finite
covers.

\medskip

{\bf Proposition 2.4.9.} {\sl Let\/ $(\B_1,\dots,\B_n)$ be a simple
transversal cycle of bisectors with total angle\/ $2\pi$ at some point
and suppose that the polyhedron\/ $P$ of the cycle is fibred. Then, for
odd\/ $n$, the complex hyperbolic manifold\/ $\Bbb B/T_n$ is
diffeomorphic to a disc bundle over a closed orientable surface of
genus\/ $n-3$ with Euler number\/ $4e_P$. For even\/ $n$, the complex
hyperbolic manifold\/ $\Bbb B/G_n$ is diffeomorphic to a disc bundle
over a closed orientable surface of genus\/ $\frac n2-1$ with Euler
number\/ $2e_P$.}\hfill$_\blacksquare$

\medskip

\centerline{\bf2.5.~Transversal triangles}

\medskip

Transversal triangles of bisectors turn out to be important geometrical
objects. They can serve as building blocks for constructing (fibred)
fundamental polyhedra, and not only in the way used in this paper. They
are naturally equipped with an isometry of their vertices (complex
geodesics), called holonomy. In this subsection, we prove that every
transversal triangle provides a fibred polyhedron (Theorem 2.5.2).

\smallskip

{\bf2.5.1.~Holonomy of a triangle.} Let $\B$ be an oriented full
bisector. The meridional curves induce an identification of the slices
in $\B$ called the {\it slice identification along\/} $\B$. By Lemma
2.1.4, for an oriented segment $\B[S,S']$, the slice identification
between $S$ and $S'$ is induced by the reflection in the middle slice
of $\B[S,S']$. Therefore, this identification is an isometry between
the slices. Clearly, such an identification can also be described via a
suitable one-parameter subgroup in $\PU$.

\smallskip

Let $(\B_1,\B_2,\B_3)$ be a counterclockwise-oriented transversal
triangle of bisectors (see Definition 2.3.5). Denote by $R_i$ the
reflection in the middle slice of $\B_i$ and let $S_1$ stand for the
initial slice of $\B_1$. We~put $\varphi:=R_3R_2R_1$. Consider the
identification of $S_1$ with itself given by the slice identification,
first, along~$\B_1$, then along $\B_2$, and, finally, along $\B_3$.
This identification is obviously induced by $\varphi\in\PU$,
$\varphi S_1=S_1$. We~call $\varphi$ the {\it holonomy\/} of the
triangle $(\B_1,\B_2,\B_3)$. The triangle is said to be {\it elliptic,
parabolic, hyperbolic,} or {\it trivial\/} if $\varphi$, being
restricted to $S_1$, is elliptic, parabolic, hyperbolic, or the
identity.

\vskip10pt

\noindent
$\hskip318pt\vcenter{\hbox{\epsfbox{Picture.12}}}$

\vskip-160pt

\rightskip145pt

Suppose that the triangle is elliptic. For every $s_1\in S_1$,
following the counterclockwise orientation of the triangle, we can draw
a meridional curve $b$ that begins at $s_1$. This curve ends at some
$s'_1\in S_1$. In other terms, $\varphi s_1=s'_1$. We take
$s_1\in\partial_\infty S_1$. Then $s'_1\in\partial_\infty S_1$.
Following the natural orientation of the circle $\partial_\infty S_1$,
we can draw an arc $a\subset\partial_\infty S_1$ from $s'_1$ to $s_1$
obtaining a closed oriented curve $c:=b\cup a\subset T$, where $T$
stands for the torus of the triangle. The curve $c$ is said to be {\it
standard.}

In the case of a hyperbolic triangle, there are two points in
$\partial_\infty S_1$ fixed by $\varphi$. They divide
$\partial_\infty S_1$ into two $\varphi$-invariant parts: the R-{\it
part\/} where $\varphi$ moves the points in the counterclockwise
direction and the L-{\it part\/} where $\varphi$ moves the points in
the clockwise direction. Let $s_1\in\partial_\infty S_1$ be a point in
the interior of the L-part. As above, we can draw a meridional curve
$b$ beginning at $s_1$ and ending at
$s'_1=\varphi s_1\in\partial_\infty S_1$. Obviously, $s'_1$ is also in
the L-part. Again, we draw an arc $a$ from\break

\vskip-12pt

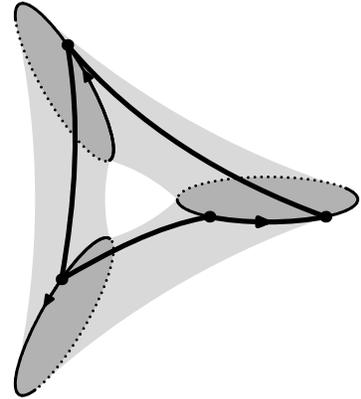
\rightskip0pt

\noindent
$s'_1$ to $s_1$ in the counterclockwise direction. Clearly, $a$ is
contained in the L-part of $\partial_\infty S_1$. We call the closed
oriented curve $c:=b\cup a\subset T$ {\it standard\/} as well. Note
that there are two closed meridional curves in $T$ (they correspond to
the fixed points of $\varphi$), both isotopic to a standard one.

For a parabolic triangle, we distinguish the R-{\it parabolic\/} and
L-{\it parabolic\/} cases. Exactly one point in $\partial_\infty S_1$
is fixed by $\varphi$. The isometry $\varphi$ moves all the other
points in $\partial_\infty S_1$ in the same direction, counterclockwise
for an R-parabolic triangle and clockwise for an L-parabolic triangle.
As above, we~define a {\it standard\/} curve for an L-parabolic
triangle. In $T$, this curve is isotopic to the closed meridional
curve. In the case of an R-parabolic or trivial triangle, we define no
standard curve.

\smallskip

We now extend the definition of the {\it{\rm L}-part\/} of
$\partial_\infty S_1$ to elliptic, parabolic, and trivial triangles.
For~an elliptic triangle, it is the entire $\partial_\infty S_1$. For
an L-parabolic triangle, it is $\partial_\infty S_1$ without the fixed
point of~$\varphi$. For the other two cases, it is empty.

\smallskip

Counterclockwise-oriented transversal triangles of bisectors are simple
configurations by Lemma 2.3.2. Hence, they fall under Definition 2.2.2.
In order to prove that a polyhedron glued from transversal triangles is
fibred and to be able to calculate its Euler number, we need the
following

\medskip

{\bf Theorem 2.5.2.} {\sl Let\/ $(\B_1,\B_2,\B_3)$ be a
counterclockwise-oriented transversal triangle of bisectors. Then the
triangle can be neither trivial nor\/ {\rm R}-parabolic. The
polyhedron\/ $P$ of the triangle is fibred and its standard curve is
trivializing.}

\medskip

The proof of Theorem 2.5.2 is postponed until the end of this
subsection.

\medskip

{\bf Lemma 2.5.3 {\rm[Appendix, Lemma 4.4.18]}.} {\sl Let\/ $C$ be a
complex geodesic and let\/
$\Delta_0:=\Delta(p_1,p_2,p_3)\allowmathbreak\subset\nomathbreak C$ be
a counterclockwise-oriented geodesic triangle,
$p_1,p_2,p_3\in C\cap\Bbb B$. Denote by\/ $S_1$ the complex geodesic
passing through\/ $p_1$ and orthogonal to\/ $C$ and by\/ $R_i$, the
reflection in the complex geodesic passing through the middle point
of\/ $\G[p_i,p_{i+1}]$ and orthogonal to\/ $C$ {\rm(}the indices are
modulo\/ $3${\rm)}. Then the isometry\/ $\varphi:=R_3R_2R_1$, being
restricted to\/ $S_1$, is a rotation about\/ $p_1$ by the angle\/
$-2\Area\Delta_0$.}

\medskip

{\bf Lemma 2.5.4.} {\sl Theorem\/ {\rm2.5.2\/} holds for every triangle
of bisectors with common complex spine.}

\medskip

{\bf Proof.} The first assertion in Theorem 2.5.2, in the case of
common complex spine, follows from Lemma 2.5.3. It remains to prove the
second assertion.

Every triangle of bisectors with common complex spine $C$ is built over
a usual geodesic triangle
$\Delta_0:=\Delta(p_1,p_2,p_3)\subset C\simeq{\Bbb H}_\Bbb C^1$. (The
corresponding triangle of bisectors is simply the complex projective
cone over $\Delta_0$ intersected with $\overline{\Bbb B}$; the vertex
of the cone is the polar point to $C$.) Clearly, $\Delta_0$ is
counterclockwise-oriented. Therefore, $\Area\Delta_0\in(0,\frac\pi4)$
(the metric we use differs by the factor of $4$ from Poincar\'e's one
(see Corollary 4.1.18)). By Lemmas 2.1.4 and 2.5.3, the restriction of
$\varphi$ to the slice $S_1$ is a rotation about $p_1$ by the angle
$-2\Area\Delta_0$. Hence, $\ell a:=2\Area\Delta_0\in(0,\frac\pi2)$ is
the angular measure (with respect to the centre $p_1$) of the arc $a$
ending at $s_1\in\partial_\infty S_1$. We can contract $\Delta_0$ to
$p_1$ inside $\Delta_0$. We put
$\Delta_t:=\Delta\big(p_1,p_2(t),p_3(t)\big)$, $t\in[0,1]$, where
$p_i(t)\in\G[p_i,p_1]$, $p_i(0)=p_i$, and $p_i(1)=p_1$, $i=2,3$. For
every $t\in[0,1]$, the triangle $\Cal C_t$ of bisectors built over
$\Delta_t$ contains the slice $S_1$. We draw the standard curve $c_t$
of $\Cal C_t$ that starts with $s_1$. By Lemma 2.1.5, the meridional
curve depends continuously on its initial point and on the bisector
involved. Since $\partial_\infty P_t\subset\partial_\infty P_0$ and
$\ell a_t\to0$ while $t\to1$, we can see now that $c_0$ contracts to
$s_1$ inside $\partial_\infty P_0$.\hfill$_\blacksquare$

\medskip

{\bf2.5.5.}~Let $\Delta(S_1,S_2,S_3)$ be an oriented triangle of
bisectors. Denote by $g_i$ the polar point to $S_i$, $i=1,2,3$, and
define
$$t_{ij}:=\sqrt{\ta(g_i,g_j)}{\text{ \ \rm for }}i\ne
j,\quad\varkappa:=\frac{\langle g_1,g_2\rangle\langle
g_2,g_3\rangle\langle g_3,g_1\rangle}{\langle g_1,g_1\rangle\langle
g_2,g_2\rangle\langle g_3,g_3\rangle},\quad\varepsilon:=
\frac\varkappa{|\varkappa|},\quad\varepsilon_0:=
\Re\varepsilon,\quad\varepsilon_1:=\Im\varepsilon.$$
The fact that the $S_i$'s are ultraparallel implies $t_{ij}>1$. It
follows from Sylvester's Criterion (see [KoM, p.~113]) that the
determinant of the Gram matrix of $g_1,g_2,g_3$ is nonpositive, that
is, $1+2t_{12}t_{23}t_{31}\varepsilon_0\le t_{12}^2+t_{23}^2+t_{31}^2$,
and that the equality occurs exactly when $g_1,g_2,g_3$ are in the same
projective line. The~same criterion immediately implies that the
numbers $t_{ij}$ and $\varepsilon$ constitute a complete set of
geometric invariants of an oriented triangle and that there exists an
oriented triangle corresponding to any given values of $t_{ij}$'s and
$\varepsilon$ subject to the conditions $t_{ij}>1$, $|\varepsilon|=1$,
and $1+2t_{12}t_{23}t_{31}\varepsilon_0\le t_{12}^2+t_{23}^2+t_{31}^2$.

\medskip

{\bf Criterion 2.5.6 {\rm[Appendix, Criterion 4.4.2]}.} {\sl Suppose
that\/ $1<t_{12}\le t_{23},t_{31}$. Then the triangle\/
$\Delta(S_1,S_2,S_3)$ is transversal if and only if\/
$t_{12}^2\varepsilon_0^2+t_{23}^2+t_{31}^2<
1+2t_{12}t_{23}t_{31}\varepsilon_0$.
A transversal triangle is counterclockwise-oriented if and only if\/
$\varepsilon_1<0$.}

\medskip

{\bf Lemma 2.5.7 {\rm[Appendix, Lemma 4.4.3]}.} {\sl The space of all
counterclockwise-oriented transversal triangles of bisectors, i.e., the
region in\/ $\Bbb R^4$ given by the inequalities
$$|\varepsilon_0|<1,\qquad1<t_{12}\le t_{23},t_{31},\qquad
t_{12}^2\varepsilon_0^2+t_{23}^2+t_{31}^2<
1+2t_{12}t_{23}t_{31}\varepsilon_0\le t_{12}^2+t_{23}^2+t_{31}^2,$$
is path-connected.}

\medskip

{\bf Lemma 2.5.8 {\rm[Appendix, Lemma 4.4.15]}.} {\sl A
counterclockwise-oriented transversal triangle of bisectors can be
neither\/ {\rm R}-parabolic nor trivial.}

\medskip

{\bf Proof of Theorem 2.5.2.} The first assertion in the theorem
corresponds to Lemma 2.5.8. We will prove the second assertion using a
deformation of the parameters $\varepsilon_0,t_{12},t_{23},t_{31}$ of a
given counterclock\-wise-oriented transversal triangle of bisectors. By
Criterion 2.5.6,
$\varepsilon:=\varepsilon_0-i\sqrt{1-\varepsilon_0^2}$. By Lemma 2.5.7,
we can always reach a triangle of the type dealt with in Lemma 2.5.4.
This continuous deformation of a triangle induces a continuous
deformation of the holonomy $\varphi$, of the nonordered pair of the
fixed points of $\varphi$ in the projective line of $S_1$, and, hence,
of the L-part of $\partial_\infty S_1$ (nonempty by Lemma 2.5.8). Now
we can choose $s_1$ in the L-part of $\partial_\infty S_1$ so that
$s_1$ varies continuously during the deformation. The point $s_1$
generates a standard curve $c$. Since, by Lemma 2.1.5, the meridional
curve depends continuously on its initial point and on the bisector
involved, the curve $c$ varies continuously during the deformation.

In other words, we obtain an isotopy of a 3-manifold $F$
($=\partial_\infty P$) in $\Bbb S^3$ ($=\partial_\infty\Bbb B$)
equipped with a simple curve ($=c$) on its boundary
($=\partial\partial_\infty P$) which is a simple torus ($=T$). At the
end of this isotopy, $F$ is a solid torus and the curve on its boundary
is contractible in $F$. Therefore, the initial polyhedron
$\partial_\infty P$ is a solid torus and its standard curve is
contractible.\hfill$_\blacksquare$

\bigskip

\centerline{\bf3.~Series of explicit examples}

\medskip

In this section, we construct an explicit series of discrete groups
with generators $U,W$ and defining relations $U^n=W^n=(U^{-1}W)^2=1$.
Taking a subgroup of finite index, we arrive at a disc bundle over a
closed orientable surface. Then we calculate the Toledo invariant and
find the Euler number of these bundles. Finally, we list and discuss
many particular examples.

\medskip

\centerline{\bf3.1.~A couple of transversal triangles}

\medskip

Following the line of Section 2, we construct two
counterclockwise-oriented transversal triangles $\Delta(C,M_1,M_2)$ and
$\Delta(S_2,M_2,M_1)$ such that $\Delta(S_2,M_2,M_1)$ is transversally
adjacent to $\Delta(C,M_1,M_2)$. Such triangles will be used as
building blocks for fundamental polyhedra. During this subsection,
we~elaborate some explicit and easily verifiable conditions, called the
quadrangle conditions (see 3.1.17), providing the desired properties of
the triangles in question (such as those listed in Proposition 3.1.18).

All isometries in this subsection are considered as belonging to $\SU$.

\smallskip

{\bf3.1.1.}~First, we look for a regular elliptic isometry $W\in\SU$
(see [Gol, p.~203] for the definition) and~for a reflection
$R:=R(m)\in\SU$ in a complex geodesic
$M_1:=\Bbb Pm^\perp\cap\overline{\Bbb B}$ such that $U:=WR$ is a
regular elliptic isometry. In what follows, the orthonormal basis
$q_1,q_2,q_3\in V$ of signature $-++$ is formed by the eigenvectors of
$W$. In this basis,
$$m=\left(\smallmatrix m_1\\m_2\\m_3\endsmallmatrix\right),\qquad
R=\left(\smallmatrix-2m_1^2-1&2m_1m_2&2m_1m_3\\-2m_1m_2&2m_2^2-1&
2m_2m_3\\-2m_1m_3&2m_2m_3&2m_3^2-1\endsmallmatrix\right),\qquad
W=\left(\smallmatrix
w_1^2&0&0\\0&w_2^2&0\\0&0&w_3^2\endsmallmatrix\right),$$
where $m_1,m_2,m_3$ will be nonnegative real numbers such that
$\langle m,m\rangle=1$ and the matrix form of $R$ is easily derivable
from (2.1.3). We will denote by $h_1,h_2,h_3$ the eigenvectors of $U$
and by $u_1^2,u_2^2,u_3^2$ the corresponding eigenvalues.

Take $k,l,n,p\in\Bbb Z$ subject to the conditions
$$0\le k\le l\le n-3,\qquad p=1,2\eqno{\bold{(3.1.2)}}$$
and define
$$u_1:=\exp\frac{(2np-k)\pi i}{3n},\qquad u_2:=\exp\frac{(2np-k-3)\pi
i}{3n},\qquad u_3:=\exp\frac{(2np+2k+3)\pi i}{3n},$$
$$w_1:=\exp\frac{l\pi i}{3n},\quad w_2:=\exp\frac{(l+3)\pi i}{3n},\quad
w_3:=\exp\frac{-(2l+3)\pi i}{3n},\quad
v:=\frac{u_1^2+u_2^2+u_3^2+w_1^2+w_2^2+w_3^2}2.$$
Observe the following straightforward facts: $u_1^2,u_2^2,u_3^2$ are
pairwise distinct; $w_1^2,w_2^2,w_3^2$ are pairwise distinct;
$u_1u_2u_3=w_1w_2w_3=1$; $u_1^{2n}=u_2^{2n}=u_3^{2n}$; and
$w_1^{2n}=w_2^{2n}=w_3^{2n}$.

\smallskip

The inequalities $0\le l\le n-3$ imply the inequalities
$0<\Re\big((w_1^2-w_2^2)w_2\big)$ and
$0<\Re\big((w_1^2-w_3^2)w_3\big)$. Requiring that
$$\Re(vw_2)<\Re
w_2^3,\qquad\Re(vw_3)\le\Re(w_1^2w_3),\eqno{\bold{(3.1.3)}}$$
we obtain $\Re\big((w_1^2-w_2^2)w_2\big)<\Re\big((w_1^2-v)w_2\big)$ and
$0\le\Re\big((w_1^2-v)w_3\big)$. Therefore, we can define
$$m_2:=\sqrt{\frac{\Re\big((w_1^2-v)w_2\big)}
{\Re\big((w_1^2-w_2^2)w_2\big)}},\qquad
m_3:=\sqrt{\frac{\Re\big((w_1^2-v)w_3\big)}
{\Re\big((w_1^2-w_3^2)w_3\big)}},\qquad m_1:=\sqrt{m_2^2+m_3^2-1}$$
such that $m_1>0$ and $m_2>1$.

\medskip

{\bf Lemma 3.1.4.} {\sl Under conditions\/ {\rm(3.1.2)} and\/
{\rm(3.1.3)}, we have\/ $\tr(WR)=u_1^2+u_2^2+u_3^2$.}

\medskip

{\bf Proof.} The equality in the lemma has the form
$$(-2m_2^2-2m_3^2+1)w_1^2+(2m_2^2-1)w_2^2+(2m_3^2-1)w_3^2=
u_1^2+u_2^2+u_3^2$$
which is equivalent to
$$m_2^2(w_2^2-w_1^2)+m_3^2(w_3^2-w_1^2)=v-w_1^2.$$
By definition,
$$m_2^2=\frac{w_2(v-w_1^2)+\overline w_2(\overline v-\overline
w_1^2)}{w_2(w_2^2-w_1^2)+\overline w_2(\overline w_2^2-\overline
w_1^2)}=$$
$$=\frac{\big(w_2^2(v-w_1^2)+(\overline v-\overline
w_1^2)\big)w_1^2w_2^2}{w_1^2w_2^4(w_2^2-w_1^2)+(w_1^2-w_2^2)}=
\frac{w_1^2w_2^4(v-w_1^2)+w_1^2w_2^2(\overline
v-\overline w_1^2)}{(w_1^2w_2^4-1)(w_2^2-w_1^2)}.$$

\smallskip

\noindent
Similarly,
$$m_3^2=\frac{w_1^2w_3^4(v-w_1^2)+w_1^2w_3^2(\overline v-\overline
w_1^2)}{(w_1^2w_3^4-1)(w_3^2-w_1^2)}.$$
It follows from $w_1w_2w_3=1$ that
$$\frac{w_1^2w_2^2}{w_1^2w_2^4-1}+\frac{w_1^2w_3^2}{w_1^2w_3^4-1}=0,
\qquad\frac{w_1^2w_2^4}{w_1^2w_2^4-1}+\frac{w_1^2w_3^4}{w_1^2w_3^4-1}=
1.$$
Now we can see that
$(w_2^2-w_1^2)m_2^2+(w_3^2-w_1^2)m_3^2=v-w_1^2$.\hfill$_\blacksquare$

\medskip

Obviously $R,W\in\SU$. Hence, $U:=WR\in\SU$. By [Gol, p.~204, Theorem
6.2.4] and Lemma~3.1.4, $U$ is a regular elliptic isometry with
eigenvalues $u_1^2,u_2^2,u_3^2$. It follows from
$u_1^{2n}=u_2^{2n}=u_3^{2n}$ and $w_1^{2n}=w_2^{2n}=w_3^{2n}$ that
$U^n$ and $W^n$ belong to the centre of $\SU$.

\smallskip

{\bf3.1.5.}~Let us find the eigenvector $h_i$ of $U$ corresponding to
$u_i^2$. In order to exclude trivial cases, we~require that
$$u_i^2+w_j^2\ne0\quad\text{for all }i,j.\eqno{\bold{(3.1.6)}}$$

{\bf Lemma 3.1.7.} {\sl Under conditions\/ {\rm(3.1.2)}, {\rm(3.1.3)},
and\/ {\rm(3.1.6),}
$\pmatrix\frac{m_1}{u_i^2w_1^{-2}+1}\\\frac{m_2}{u_i^2w_2^{-2}+1}\\
\frac{m_3}{u_i^2w_3^{-2}+1}\endpmatrix$
is the eigenvector of\/ $U$ corresponding to\/ $u_i^2$.}

\medskip

{\bf Proof.} Let
$x=\left(\smallmatrix x_1\\x_2\\x_3\endsmallmatrix\right)\ne0$. We fix
some $i=1,2,3$. The condition that $x$ is the eigenvector of $U=WR$
corresponding to $u_i^2$ is given by the equalities
$2w_j^2m_jf(x)=(u_i^2+w_j^2)x_j$, $j=1,2,3$, where
$f(x):=-m_1x_1+m_2x_2+m_3x_3$. Thus, the inequalities (3.1.6) imply
$f(x)\ne0$ and we can take
$x_j=\displaystyle\frac{m_j}{u_i^2w_j^{-2}+1}$,
$j=1,2,3$.\hfill$_\blacksquare$

\medskip

{\bf3.1.8.}~In order to introduce the desired transversally adjacent
triangles, we put
$$h_1:=\pmatrix\frac{m_1}{u_1^2w_1^{-2}+1}\\\frac{m_2}{u_1^2w_2^{-2}+1}
\\\frac{m_3}{u_1^2w_3^{-2}+1}\endpmatrix,\quad
h_2:=\pmatrix\frac{m_1}{u_2^2w_1^{-2}+1}\\\frac{m_2}{u_2^2w_2^{-2}+1}\\
\frac{m_3}{u_2^2w_3^{-2}+1}\endpmatrix,\quad
M_2:=WM_1,\quad C:=\Bbb Pq_2^\perp\cap\overline{\Bbb B},\quad
S_2:=\Bbb Ph_2^\perp\cap\overline{\Bbb B}$$
and require that
$$\langle h_1,h_1\rangle<0,\qquad1<\ta(m,Wm),\qquad1<\ta(m,h_2).
\eqno{\bold{(3.1.9)}}$$

As is well known, the eigenvectors of a regular elliptic isometry are
pairwise orthogonal. Hence, only one of the $h_i$'s can be negative,
implying that $h_2\notin\overline{\Bbb B}$. So, $S_2$ is a complex
geodesic and $h_1\in S_2$. Obviously, $Wm$ is the polar point to $M_2$.
The last two inequalities in (3.1.9) respectively mean that the complex
geodesics $M_1$ and $M_2$ are ultraparallel and that the complex
geodesics $M_1$ and $S_2$ are ultraparallel [Gol, p.~100]. Since
$\langle q_2,m\rangle=m_2$, we obtain $\ta(q_2,m)=m_2^2>1$ and conclude
that the complex geodesics $C$ and $M_1$ are ultraparallel as well.
Now, it follows from $WC=C$ and from $M_2=WM_1$ that $C$ and $M_2$ are
ultraparallel. Also, from $RM_1=M_1$ and $U=WR$, we deduce that
$M_2=WM_1=WRM_1=UM_1$. The facts that $M_1$ and $S_2$ are
ultraparallel, that $M_2=UM_1$, and~that $S_2=US_2$ imply that $M_2$
and $S_2$ are ultraparallel. Thus, we get the triangles
$\Delta(C,M_1,M_2)$ and $\Delta(S_2,M_2,M_1)$.

The points $q_2,m,h_2,Wm$ are polar to the complex geodesics
$C,M_1,S_2,M_2$. The equalities $M_2=UM_1=WM_1$, $US_2=S_2$, and $WC=C$
imply that $\ta(q_2,m)=\ta(Wm,q_2)$ and $\ta(h_2,Wm)=\ta(m,h_2)$.

We introduce the invariants
$$t_{12}:=t_{31}:=\sqrt{\ta(q_2,m)}=\sqrt{\ta(Wm,q_2)},\quad
t'_{12}:=t'_{31}:=\sqrt{\ta(h_2,Wm)}=\sqrt{\ta(m,h_2)},$$
$$t_{23}:=t'_{23}:=\sqrt{\ta(m,Wm)},$$
$$\varkappa:=\langle q_2,m\rangle\langle m,Wm\rangle\langle
Wm,q_2\rangle,\qquad\varepsilon:=\frac\varkappa{|\varkappa|},\qquad
\varepsilon_0:=\Re\varepsilon,\qquad\varepsilon_1:=\Im\varepsilon,$$
$$\varkappa':=\langle h_2,Wm\rangle\langle Wm,m\rangle\langle
m,h_2\rangle,\qquad\varepsilon':=\frac{\varkappa'}{|\varkappa'|},\qquad
\varepsilon'_0:=\Re\varepsilon',\qquad\varepsilon'_1:=\Im\varepsilon'$$
of the oriented triangles $\Delta(C,M_1,M_2)$ and $\Delta(S_2,M_2,M_1)$
(see 2.5.5) and require that
$$t_{12}^2\varepsilon_0^2+t_{23}^2+t_{31}^2<
1+2t_{12}t_{23}t_{31}\varepsilon_0,\qquad
t_{23}^2\varepsilon_0^2+t_{31}^2+t_{12}^2<1+2t_{23}t_{31}t_{12}
\varepsilon_0,\qquad\varepsilon_1<0,\eqno{\bold{(3.1.10)}}$$
$${t'_{12}}^{\!\!\!\!2}\,{\varepsilon'}^2_{\!\!0}+{t'_{23}}^{\!\!\!\!2}
+{t'_{31}}^{\!\!\!\!2}<1+2t'_{12}t'_{23}t'_{31}\varepsilon'_0,\qquad
{t'_{23}}^{\!\!\!\!2}\,{\varepsilon'}^2_{\!\!0}+{t'_{31}}^{\!\!\!\!2}
+{t'_{12}}^{\!\!\!\!2}<1+2t'_{23}t'_{31}t'_{12}\varepsilon'_0,\qquad
\varepsilon'_1<0.\eqno{\bold{(3.1.11)}}$$
Since $t_{12}=t_{31}$ and $t'_{12}=t'_{31}$, both triangles
$\Delta(C,M_1,M_2)$ and $\Delta(S_2,M_2,M_1)$ are transversal and
counterclockwise-oriented by Criterion 2.5.6.

\medskip

{\bf Criterion 3.1.12 {\rm[Appendix, Criterion 4.3.3]}.} {\sl Let\/
$C,C_1,C_2$ be complex geodesics and let\/ $g,g_1,g_2$ be their polar
points. Suppose that\/ $C$ is ultraparallel to\/ $C_i$ for all\/
$i=1,2$. Then the full bisectors\/ ${\prec}\B[C,C_1]{\succ}$ and\/
${\prec}\B[C,C_2]{\succ}$ are transversal along their common slice\/
$C$ if and only if
$$\bigg|\Re\frac{\langle g_1,g_2\rangle\langle g,g\rangle}{\langle
g_1,g\rangle\langle g,g_2\rangle}-1\bigg|<
\sqrt{1-\frac1{\ta(g,g_1)}}\cdot\sqrt{1-\frac1{\ta(g,g_2)}}.$$}

In the following lemma, we assume that $\Arg$ takes values in
$[0,2\pi)$.

\medskip

{\bf Lemma 3.1.13 {\rm[Appendix, Corollary 4.3.2]}.} {\sl Let\/
$\Cal E\in\SU$ be a regular elliptic isometry and let\/
$e_1,e_2,e_3\in\Bbb P$ be the points corresponding to the eigenvectors
of\/ $\Cal E$ such that\/ $e_1\in\Bbb B$. Denote by\/ $E$ the complex
geodesic with the polar point\/ $e_2$. Let\/ $D$ be a complex geodesic
ultraparallel to\/ $E$. Then the oriented angle from\/ $\B[E,D]$ to\/
$\B[E,\Cal ED]$ at\/ $e_1\in E$ equals\/ $\Arg(\xi_2\xi_1^{-1})$,
where\/ $\xi_i$ stands for the eigenvalue of\/ $\Cal E$ corresponding
to\/ $e_i$, $i=1,2$.}

\medskip

It follows from $M_2=WM_1$, $M_1=U^{-1}M_2$, and Lemma 3.1.13 that both
the angle from $\B[C,M_1]$ to $\B[C,M_2]$ at $q_1\in C$ and the angle
from $\B[S_2,M_2]$ to $\B[S_2,M_1]$ at $h_1\in S_2$ equal
$\displaystyle\frac{2\pi}n$.

We require that
$$\bigg|\Re\frac{\langle h_2,q_2\rangle\langle m,m\rangle}{\langle
h_2,m\rangle\langle m,q_2\rangle}-1\bigg|<\sqrt{1-\frac1{\ta(m,h_2)}}
\cdot\sqrt{1-\frac1{\ta(m,q_2)}}\eqno{\bold{(3.1.14)}}$$
(although it is possible to deduce (3.1.14) from (3.1.2) and (3.1.3)).
By Criterion 3.1.12, the full bisectors ${\prec}\B[M_1,C]{\succ}$ and
${\prec}\B[M_1,S_2]{\succ}$ are transversal along their common slice
$M_1$. It follows from
${\prec}\B[M_1,S_2]{\succ}={\prec}\B[RS_2,M_1]{\succ}$ that the full
bisectors ${\prec}\B[WM_1,WC]{\succ}$ and
${\prec}\B[WRS_2,WM_1]{\succ}$ are transversal along their common slice
$M_2=WM_1$. Since $WC=C$, $WR=U$, and $US_2=S_2$, the full bisectors
${\prec}\B[M_2,C]{\succ}$ and ${\prec}\B[S_2,M_2]{\succ}$ are
transversal along their common slice $M_2$.

\medskip

{\bf Lemma 3.1.15 {\rm[Appendix, Lemma 4.2.15]}.} {\sl Let\/ $g_1,g_2$
be the polar points to the ultraparallel complex geodesics\/ $C_1,C_2$.
A point\/ $x\in\Bbb B$ is on the side of the normal vector to the
oriented bisector\/ ${\prec}\B[C_1,C_2]{\succ}$ if and only if\/
$\Im\displaystyle\frac{\langle g_1,x\rangle\langle
x,g_2\rangle}{\langle g_1,g_2\rangle}\ge0$.}

\medskip

We require that
$$\Im\frac{\langle q_2,h_1\rangle\langle h_1,m\rangle}{\langle
q_2,m\rangle}\ge0,\qquad\Im\frac{\langle Wm,h_1\rangle\langle
h_1,q_2\rangle}{\langle Wm,q_2\rangle}\ge0.\eqno{\bold{(3.1.16)}}$$
By Lemma 3.1.15, conditions (3.1.16) express the fact that $h_1$
belongs to the interior sector at $C$ of the triangle
$\Delta(C,M_1,M_2)$.

\smallskip

{\bf3.1.17.}~For brevity, we call conditions (3.1.2), (3.1.3), (3.1.6),
(3.1.9), (3.1.10), (3.1.11), (3.1.14), (3.1.16) the {\it quadrangle\/}
conditions. Put $s_2:=h_1$ and $q:=q_1$.

\medskip

{\bf Proposition 3.1.18.} {\sl From the quadrangle conditions, we
conclude that both\/ $U$ and\/ $W$ are regular elliptic of order\/ $n$
in\/ $\PU$, that the counterclockwise-oriented transversal triangle\/
$\Delta(S_2,M_2,M_1)$ is transversally adjacent to the
counterclockwise-oriented transversal triangle\/ $\Delta(C,M_1,M_2)$,
and that both the angle from\/ $\B[C,M_1]$ to\/ $\B[C,M_2]$ at\/
$q\in C$ and the angle from\/ $\B[S_2,M_2]$ to\/ $\B[S_2,M_1]$ at\/
$s_2\in S_2$ equal\/ $\displaystyle\frac{2\pi}n$.}\hfill$_\blacksquare$

\medskip

\centerline{\bf3.2.~Examples of disc bundles}

\medskip

From the triangles $\Delta(C,M_1,M_2)$ and $\Delta(S_2,M_2,M_1)$
constructed in the previous subsection, we~obtain a quadrangle $\Cal Q$
whose polyhedron is fibred and constitutes a fundamental region for the
group $K_n$ generated by $U,W$ with the defining relations
$U^n=W^n=(U^{-1}W)^2=1$ in $\PU$. Passing to a subgroup of index $n$,
that is, gluing $n$ copies of the form $W^i\Cal Q$, we arrive at a
transversal simple cycle of bisectors which provides a complex
hyperbolic disc bundle by Proposition 2.4.9. We calculate the Toledo
invariant of the representation provided by the cycle. Using
trivializing curves of the triangles, we find how to calculate the
Euler number of the polyhedron of the cycle.

Actually, the quadrangle $\Cal Q$ gives rise to a complex hyperbolic
disc bundle over the turnover orbifold $\Bbb S^2(n,n,2)$. One can
define the (rational) Euler number of such bundles as in [BSi]. Since
the turnover group $K_n$ possesses a discrete and cocompact action on
the Poincar\'e disc, one can also define the Toledo invariant of
representations $K_n\to\PU$. In our case, both numbers are easily
derivable from those of the cycle glued from copies of $\Cal Q$.

All isometries in this subsection are considered as belonging to $\SU$,
unless otherwise stated.

\smallskip

{\bf3.2.1.}~First, we list the facts following from the quadrangle
condition 3.1.17 that are not included in Proposition 3.1.18 : The
regular elliptic isometry $W$ stabilizes $C$, fixes the point
$q\in C\cap\Bbb B$, and maps $M_1$ onto $M_2$. The regular elliptic
isometry $U$ stabilizes $S_2$, fixes the point $s_2\in S_2\cap\Bbb B$,
and maps $M_1$ onto $M_2$. The isometry $R:=W^{-1}U$ is the reflection
in $M_1$. With a straightforward calculation, we~obtain
$W^nU^{-n}=\delta$, where
$\delta:=\exp\displaystyle\frac{2(k+l+np)\pi i}3$.

Note that the segment of bisector $W\B[M_1,S_2]$ is a continuation of
the segment $\B[S_2,M_2]$. Indeed, the reflection $R^W:=WRW^{-1}$ in
$M_2$ maps $S_2$ onto $WRW^{-1}S_2=WU^{-1}S_2=WS_2$.

\smallskip

{\bf3.2.2.~Discreteness of $K_n$.} Let $\Delta(C,M_1,M_2)$ and
$\Delta(S_2,M_2,M_1)$ be the triangles constructed in Subsection 3.1
and assume that the quadrangle conditions 3.1.17 are valid. We call the
configuration
$$\Cal Q:=\big(\B[C,M_1],\B[M_1,S_2],\B[S_2,M_2],\B[M_2,C]\big)$$
the {\it quadrangle\/} of bisectors. Since the transversal triangles
$\Delta(C,M_1,M_2)$ and $\Delta(S_2,M_2,M_1)$ are transversally
adjacent, the configuration $\Cal Q$ is transversal. By Remark 2.3.6,
$\Cal Q$ is simple and the polyhedron $Q$ of $\Cal Q$ is the gluing of
the polyhedra of the triangles (see 2.4.4). By Theorem 2.5.2 and
Remark~2.4.5, $Q$ is fibred.

Applying Theorem 2.2.5, we will show that $Q\cap\Bbb B$ is a
fundamental polyhedron for $K_n\subset\PU$. Obviously, $W$ maps
$\B[C,M_1]$ onto $\B[C,M_2]$ and $U$ maps $\B[S_2,M_1]$ onto
$\B[S_2,M_2]$. There are three geometric cycles of edges (see 2.2.3 for
the definition). The cycle of $C$ has total angle $2\pi$ at $q\in C$ by
Proposition 3.1.18. The same concerns the cycle of $S_2$ and
$s_2\in S_2$. The geometric cycle of $M_1$ has length $4$ due to the
relation $U^{-1}WU^{-1}W=1$. In order to verify that the total angle at
a point in $M_1\cap\Bbb B$ is $2\pi$, we recall that $W$ sends
$\B[C,M_1]$ onto $\B[C,M_2]$ and that $W\B[M_1,S_2]$ is a continuation
of $\B[S_2,M_2]$ (see 3.2.1). This implies that the sum of two
consecutive angles at the point in question equals $\pi$ giving the
total of $2\pi$.

\smallskip

{\bf3.2.3.~Cycle $\Cal C$, polyhedron $P$, and disc bundle
$M(n,l,k,p)$.} Gluing the copies $W^i\Cal Q$,
$i=\allowmathbreak1,2,\dots,n$, we obtain a simple transversal cycle of
bisectors $\Cal C:=(\B_1,\dots,\B_n)$, where $\B_i:=\B[S_i,S_{i+1}]$
and $S_i:=W^{i-2}S_2$ (the indices are modulo $n$). Indeed, $S_2$ and
$M_2$ are ultraparallel and the reflection in $M_2$ maps $S_2$ onto
$S_3=WS_2$ by 3.2.1. Hence, $S_2$ and $S_3$ are ultraparallel and $M_2$
is the middle slice of $\B_2=\B[S_2,S_3]$. So, $R_i:=R^{W^{i-1}}$ is
the reflection in the middle slice
$M_i:=W^{i-2}M_2$ of $\B_i=W^{i-2}\B_2$~and
$$R_n\dots R_1=W^n(W^{-1}R)^n=W^{n-1}(RW^{-1})^nW=W^{n-1}U^{-n}W=
W^nU^{-n}=\delta\eqno{\bold{(3.2.4)}}$$
by 3.2.1. Since the triangle $\Delta(S_2,M_2,M_1)$ is transversal, the
bisectors ${\prec}\B_1{\succ}$ and ${\prec}\B_2{\succ}$ are transversal
along $S_2$. This implies the transversality of $\Cal C$. The fact that
$\Delta'_2:=\Delta(S_2,M_2,M_1)$ is transversally adjacent to
$\Delta_2:=\Delta(C,M_1,M_2)$ implies that
$\Delta'_i:=W^{i-2}\Delta'_2$ is transversally adjacent to
$\Delta_i:=W^{i-2}\Delta_2$ and that $C=WC$ is a centre of $\Cal C$. By
Proposition 3.1.18, the central angle at $q\in C$ is $2\pi$ because $W$
fixes $q$. By Criterion 2.3.9, $\Cal C$ is simple.

Denote by $P$ the polyhedron of the cycle $\Cal C$. Being glued from
fibred quadrangles, $P$ is fibred.

The total angle of $\Cal C$ at $s_2\in S_2\cap\Bbb B$ equals $2\pi$.
Indeed, the oriented angle from $\B_2{\succ}$ to $\B_1^-{\succ}$ at
$s_2\in S_2$ equals $\displaystyle\frac{2\pi}n$ by Proposition 3.1.18.
Put $s_i:=W^{i-2}s_2$. Then the oriented angle from $\B_i{\succ}$ to
$\B_{i-1}^-{\succ}$ at $s_i\in S_i$ equals $\displaystyle\frac{2\pi}n$.
It remains to observe that
$R_is_i=W^{i-1}RW^{1-i}W^{i-2}s_2=W^{i-1}U^{-1}s_2=W^{i-1}s_2=s_{i+1}$.

\smallskip

By Proposition 2.4.9, we arrive at the complex hyperbolic disc bundle
$M(n,l,k,p)$ diffeomorphic to $\Bbb B/G_n$ if $n$ is even and to
$\Bbb B/T_n$ if $n$ is odd, where $G_n$ or $T_n$ is a torsion-free
subgroup of index $2$ or $4$ in~$H_n$ (see 2.1.10), hence, of index
$2n$ or $4n$ in $K_n$, respectively.

\smallskip

{\bf3.2.5.~Toledo invariant of $\Cal C$.} As in 2.1.11 and Proposition
2.1.15, we can calculate the Toledo invariant of the representation of
$H_n$ (see 2.1.10) defined by the cycle $\Cal C$.

\medskip

{\bf Proposition 3.2.6.} {\sl Suppose that the quadrangle conditions\/
{\rm3.1.17} are valid. Let\/ $t$ be an integer such that\/ $0\le t<3n$
and\/ $t\equiv2np-k-l\mod3n$. Then the Toledo invariant of the
representation of\/ $H_n$ defined by the cycle\/ $\Cal C$ equals\/
$\frac23t-n$.}

\medskip

{\bf Proof.} The point $s_2\in S_2$ generates a meridional curve of
$\Cal C$ with vertices $s_{i+1}:=R_is_i$, where $R_i:=R^{W^{i-1}}$ and
$s_{n+1}=\delta s_1$ (see 3.2.4). Following the proof of Proposition
2.1.15, we take $c:=q$ and, using (2.1.16), obtain
$$\tau=\frac1\pi\sum\limits_{i=1}^n\Big(\Arg\frac{\langle
q,s_{i+1}\rangle}{\langle q,s_i\rangle}-\pi\Big).$$
We have
$$s_i=R_{i-1}\dots R_2s_2=W^{i-2}(RW^{-1})^{i-2}s_2=W^{i-2}U^{2-i}s_2=
u_1^{2(2-i)}W^{i-2}s_2$$
because $s_2=h_1$ is an eigenvector of $U$ with eigenvalue $u_1^2$ (see
3.1.17 and 3.1.1). Since $q$ is an eigenvector of $W$ with eigenvalue
$w_1^2$,
$$\langle q,s_i\rangle=u_1^{2i-4}\langle
q,W^{i-2}s_2\rangle=u_1^{2i-4}\langle
W^{2-i}q,s_2\rangle=u_1^{2i-4}w_1^{4-2i}\langle
q,s_2\rangle=(u_1^2w_1^{-2})^{i-2}\langle q,s_2\rangle.$$
Hence,
$\Arg\displaystyle\frac{\langle q,s_{i+1}\rangle}{\langle
q,s_i\rangle}=\Arg(u_1^2w_1^{-2})=\frac{2t\pi}{3n}$
and $\tau=\frac23t-n$.\hfill$_\blacksquare$

\medskip

{\bf3.2.7.~Holonomy of the quadrangle $\Cal Q$.} In order to calculate
the Euler number of $P$, we need to find explicitly a trivializing
curve of the polyhedron $Q$. Theorem 2.5.2 provides trivializing curves
of $\Delta(C,M_1,M_2)$ and $\Delta(S_2,M_2,M_1)$. Therefore, we only
need to study the contribution coming from the mutual position of the
triangles.

\smallskip

Pick a generic point $z\in\partial_\infty C$ and introduce the
following curves and points:

\vskip6pt

\noindent
$\hskip83pt\vcenter{\hbox{\epsfbox{Picture.15}}}$

$\bullet$ the meridional\footnote{We denote by $x^{-1}$ the (not
necessarily closed) curve $x$ taken with the opposite orientation.}
curve $r_1^{-1}\subset\partial_\infty\B[C,M_1]$ that begins at $z$ and
ends at $p_1\in\partial_\infty M_1$,

$\bullet$ the meridional curve $b_1^1\subset\partial_\infty\B[M_1,S_2]$
that begins at $p_1$ and ends at $y\in\partial_\infty S_2$,

$\bullet$ the naturally oriented simple arc
$a_2\subset\partial_\infty S_2$ that begins at $y$ and ends at
$Uy\in\partial_\infty S_2$,

$\bullet$ the meridional curve $b_2^0\subset\partial_\infty\B[S_2,M_2]$
that begins at $Uy$ and ends at $p_2\in\partial_\infty M_2$,

$\bullet$ the meridional curve $r_2\subset\partial_\infty\B[M_2,C]$
that begins at $p_2$ and ends at $z'\in\partial_\infty C$,

$\bullet$ the naturally oriented simple arc
$c_2\subset\partial_\infty C$ that begins at $z'$ and ends at $z$.

\smallskip

Since $Q$ is fibred, $\partial_\infty Q$ is a solid torus and the group
$H_1(\partial_\infty Q,\Bbb Z)$ is generated by $[c]$, where $c$ stands
for the naturally oriented boundary of $C$. Hence, there exists
$f\in\Bbb Z$ such that $[d_2]=f[c]$, where the closed curve $d_2$ is
given by
$$d_2:=r_1^{-1}\cup b_1^1\cup a_2\cup b_2^0\cup r_2\cup c_2.$$
In order to express $f$ in terms of the holonomies $\varphi$ of
$\Delta(S_2,M_2,M_1)$ and $\varphi'$ of $\Delta(C,M_1,M_2)$,
we~introduce the following notation. Let $S$ be an oriented circle and
let $t_1,t_2,t_3\in S$ be pairwise distinct. We define
$o(t_1,t_2,t_3)=0$ if $t_1,t_2,t_3$ are in the cyclic order of $S$ and
$o(t_1,t_2,t_3)=1$, otherwise. For an isometry $I$ of a complex
geodesic $E$ and for $e\in\partial_\infty E$, we put $\lambda(e,I)=0$
if $e$ belongs to the L-part of $\partial_\infty E$ with respect to $I$
and $\lambda(e,I)=1$, otherwise.

\medskip

{\bf Lemma 3.2.8.}
{\sl$f=\lambda(y,\varphi)+\lambda(z,\varphi')+o(y,Uy,\varphi^{-1}y)-
o(z,Wz,\varphi'z)$.}

\medskip

{\bf Proof.} First, let us assume that $z$ belongs to the L-part of
$\partial_\infty C$ and $y$, to the L-part of $\partial_\infty S_2$. We
introduce a few more curves. The meridional curve
$g\subset\partial_\infty\B[M_1,M_2]$ begins at $p_1$ and ends at some
$p'_2\in\partial_\infty M_2$. By definition, there exists a unique
meridional curve $b'\subset\partial_\infty\B[S_2,M_2]$ such that
$b'\cup g^{-1}\cup b_1^1\cup a'$ is a trivializing curve of
$\Delta(S_2,M_2,M_1)$, where $a'\subset\partial_\infty S_2$ is a simple
arc following the natural orientation of $\partial_\infty S_2$.
Clearly, $b'$ begins at $\varphi^{-1}y$ and ends at $p'_2$. Similarly,
there exists a unique meridional curve
$r'_2\subset\partial_\infty\B[M_2,C]$ such that
$r_1^{-1}\cup g\cup r'_2\cup c'$ is a trivializing curve of
$\Delta(C,M_1,M_2)$, where $c'\subset\partial_\infty C$ is a simple arc
drawn from $\varphi'z$ to $z$ following the natural orientation of
$\partial_\infty C$. So,~$r'_2$~ends at $\varphi'z$.

By Remark 2.4.5, $r_1^{-1}\cup b_1^1\cup a'\cup b'\cup r'_2\cup c'$ is
a trivializing curve of $Q$. In terms of $1$-chains modulo boundaries,
that is, in
$C_1(\partial_\infty Q,\Bbb Z)/\partial C_0(\partial_\infty Q,\Bbb Z)$,
this fact can be written as $-[r_1]+[b_1^1]+[a']+[b']+[r'_2]+[c']=0$.

Denote by $a''\subset\partial_\infty S_2$ the simple arc from $Uy$ to
$\varphi^{-1}y$ that follows the natural orientation of
$\partial_\infty S_2$. It is easy to see that
$[a_2]+[a'']-[a']=o(y,Uy,\varphi^{-1}y)[\partial_\infty S_2]$ in terms
of $1$-chains modulo boundaries. Let $c''\subset\partial_\infty C$
stand for the simple arc from $z'$ to $\varphi'z$ that follows the
natural orientation of $\partial_\infty C$. As~above,
$[c'']+[c']-[c_2]=o(z,z',\varphi'z)[c]$. Since
$-[r_2]-[b_2^0]+[a'']+[b']+[r'_2]-[c'']=0$,
$[\partial_\infty S_2]=[c]$, and
$[d_2]=-[r_1]+[b_1^1]+[a_2]+[b_2^0]+[r_2]+[c_2]$, we have
$$[d_2]=-[r_1]+[b_1^1]+[a_2]+[b_2^0]+[r_2]+[c_2]-[r_2]-[b_2^0]+[a'']+
[b']+[r'_2]-[c'']+$$
$$+[r_1]-[b_1^1]-[a']-[b']-[r'_2]-[c']=[a_2]+[a'']-[a']+[c_2]-[c'']-
[c']=\big(o(y,Uy,\varphi^{-1}y)-o(z,z',\varphi'z)\big)[c].$$

We claim that $z'=Wz$. Indeed, $U$ maps $\B[S_2,M_1]$ onto
$\B[S_2,M_2]$ and $y$ to $Uy$. Hence, it maps $b_1^1$ onto $b_2^0$. In
particular, $p_2=Up_1$. It follows from $Rp_1=p_1$ and $U=WR$ that
$Up_1=Wp_1$. So,~$p_2=Wp_1$. Now, by similar arguments, $W$ maps $r_1$
onto $r_2$, implying $z'=Wz$.

When (say) $y$ belongs to the R-part of $\partial_\infty S_2$, one can
readily understand that the required correction equals
$\lambda(y,\varphi)$.\hfill$_\blacksquare$

\medskip

Clearly, the terms $\lambda(y,\varphi)$ and $\lambda(z,\varphi')$
vanish if the triangles $\Delta(S_2,M_2,M_1)$ and $\Delta(C,M_1,M_2)$
are both elliptic.

\smallskip

{\bf3.2.9.~Euler number.} In order to find the Euler number of $P$, we
need a few auxiliary facts.

Put $U_i:=U^{W^{i-2}}$ and $s_i:=W^{i-2}s_2$. The equalities
$U_2S_2=S_2$ and $U_2s_2=s_2$ imply $U_iS_i=S_i$ and $U_is_i=s_i$ by
3.2.1.

\medskip

{\bf Lemma 3.2.10.} {\sl$R_iU_i^j=U_{i+1}^jR_i$.}

\medskip

{\bf Proof.} It suffices to show that $R_iU_i=U_{i+1}R_i$. The
relations $U=WR$ and $R^2=1$ imply the relations $RW^{-1}U=1$ and
$UR=W$. Since $R_i=R^{W^{i-1}}$,
$$R_iU_i=W^{i-1}RW^{1-i}W^{i-2}UW^{2-i}=W^{i-1}RW^{-1}UW^{2-i}=W,$$
$$U_{i+1}R_i=W^{i-1}UW^{1-i}W^{i-1}RW^{1-i}=W^{i-1}URW^{1-i}=
W.\eqno{_\blacksquare}$$

{\bf Lemma 3.2.11.} {\sl In the complex geodesic\/ $C$, the isometry\/
$W^{-1}$ is the rotation about\/ $q\in C$ by the angle\/
$\beta:=\displaystyle\frac{2(l+1)\pi}n$. In the complex geodesic\/
$S_i$, the isometry\/ $U_i$ is the rotation about\/ $s_i\in S_i$ by the
angle\/ $\alpha:=\displaystyle\frac{2(k+1)\pi}n$.}

\medskip

{\bf Proof.} In the basis $q,q_2,q_3$ of eigenvectors of $W$, the
points in $C\cap\overline{\Bbb B}$ have the form
$\left(\smallmatrix 1\\0\\z\endsmallmatrix\right)$, where $z\in\Bbb C$
and $|z|\le1$. In this way, we identify $C$ with the closed unit disc
in $\Bbb C$ centred at $0$ so that $q$ corresponds to $0$. Hence, in
terms of $z$, $W^{-1}$ acts as the multiplication by $w_1^2w_3^{-2}$,
implying the first assertion. The same arguments work for $U_2=U$ and
$s_2\in S_2$, which implies the second assertion.\hfill$_\blacksquare$

\medskip

{\bf Definition 3.2.12.} Let $a\subset\partial_\infty S_i$ be an arc
(not necessarily simple) that begins at $y$ and ends at $x=U_i^jy$,
$j\in\Bbb Z$. We call such an arc {\it integer.} Assign to an integer
arc $a\subset\partial_\infty S_i$ its {\it angular measure\/} $\ell a$
(with respect to the centre $s_i$) by the following rules:

\smallskip

$\bullet$ For $j=1$ and a simple integer arc
$a\subset\partial_\infty S_i$ drawn in the counterclockwise direction,
we put $\ell a:=\alpha$ (see Lemma 3.2.11).

$\bullet$ If $a\subset\partial_\infty S_i$ is integer, then
$\ell a^{-1}=-\ell a$.

$\bullet$ If $a=a'\cup a''$ with integer
$a',a''\subset\partial_\infty S_i$, then $\ell a=\ell a'+\ell a''$.

$\bullet$ Two integer arcs that are homotopic in $\partial_\infty S_i$
have the same angular measure.

\medskip

{\bf Lemma 3.2.13.} {\sl Let\/ $a\subset\partial_\infty S_i$ be an
integer arc. Then the arc\/ $R_ia\subset\partial_\infty S_{i+1}$ is
integer and\/ $\ell R_ia=\ell a$.}

\medskip

{\bf Proof.} The arc $a$ begins at $y$ and ends at $x=U_i^jy$,
$j\in\Bbb Z$. Hence, the arc $R_ia\subset\partial_\infty S_{i+1}$
begins at $R_iy$ and ends at $R_ix=R_iU_i^jy=U_{i+1}^jR_iy$ by Lemma
3.2.10. It remains to observe that the fact is valid for $j=1$ and any
simple integer arc $a$ drawn in the counterclockwise
direction.\hfill$_\blacksquare$

\medskip

{\bf Proposition 3.2.14.} {\sl Under the quadrangle conditions\/
{\rm3.1.17,} the Euler number of\/ $P$ equals\/ $nf-k-l-2$, where\/ $f$
is calculated in Lemma\/ {\rm3.2.8.}}

\medskip

{\bf Proof.} We introduce the following curves and points (see also
3.2.7) :
$$a_i:=W^{i-2}a_2\subset\partial_\infty S_i,\qquad
b_i^0:=W^{i-2}b_2^0\subset\partial_\infty\B[S_i,M_i],\qquad
b_i^1:=W^{i-1}b_1^1\subset\partial_\infty\B[M_i,S_{i+1}],$$
$$c_i:=W^{i-2}c_2\subset\partial_\infty C,\qquad d_i:=W^{i-2}d_2\subset
W^{i-2}\partial_\infty Q,\qquad
r_i:=W^{i-1}r_1\subset\partial_\infty\B[M_i,C].$$
The facts that $Wp_1=p_2$, $W\B[M_1,C]=\B[M_2,C]$, and $r_2$ is the
meridional curve of $\B[M_2,C]$ generated by $p_2$ imply that
$Wr_1=r_2$. So, the $r_i$'s are well defined. For the same reason, we
can define $p_i=W^{i-1}p_1\in\partial_\infty M_i$. Note that $b_i^0$
ends at $p_i$ and $b_i^1$ begins at $p_i$. Since
$b_i^0\subset\B[S_i,M_i]$ and $b_i^1\subset\B[M_i,S_{i+1}]$ are
meridional curves, they form a meridional curve
$b_i:=b_i^0\cup b_i^1\subset\B[S_i,S_{i+1}]$. We obtain a closed curve
$$\gamma:=a_1\cup b_1\cup\dots\cup a_n\cup b_n\subset T,$$
where $T:=\partial\partial_\infty P$ is the torus of $P$ (see
Definition 2.2.2).

Since $P$ is fibred, $\partial_\infty P$ is a solid torus and the group
$H_1(\partial_\infty P,\Bbb Z)$ is generated by $[s]$, where $s$ stands
for the naturally oriented boundary of an arbitrary slice of $\Cal C$.
Note also that $[s]=[c]$.

The fact that $WP=P$ and $WC=C$ implies $W[c]=[c]$. Therefore,
$[d_1\cup\dots\cup d_n]=[d_1]+\dots+[d_n]=nf[c]$. On the other hand,
$[d_1\cup\dots\cup d_n]=[\gamma]+[c_1\cup\dots\cup c_n]$ and
$[c_1\cup\dots\cup c_n]=(l+1)[c]$ by Lemma 3.2.11 because the angular
measure of every $c_i$ with respect to the centre $q$ equals $\beta$.
Thus, $[\gamma]=(nf-l-1)[c]$.

\vskip3pt

\noindent
$\hskip225pt\vcenter{\hbox{\epsfbox{Picture.13}}}$

\vskip-168pt

\rightskip233pt

Denote by $b$ the meridional curve of $\Cal C$ (see 2.1.8) that
includes $b_1$. We have $b=b'_1\cup\dots\cup b'_n$, where $b'_i$ is a
meridional curve of $\B[S_i,S_{i+1}]$ and $b'_1=b_1$. Let
$a'_1\subset\partial_\infty S_1$ be the initial point of $b_1$ and,
by~induction, let $a'_{i+1}=R_ia'_i\cup a_{i+1}$ for $i=1,\dots,n$.
Using Lemma 2.1.4, it is easy to see that the curves $a'_i\cup b_i$ and
$b'_i\cup R_ia'_i$ are well defined and homotopic in the cylinder
$\partial_\infty\B[S_i,S_{i+1}]$. In terms of $1$-chains modulo
boundaries in the solid torus $\partial_\infty P$, this implies that
$[a'_i]+[b_i]=[b'_i]+[R_ia'_i]$ for all $i$. Lemma 3.2.13 provides the
equality $\ell a'_i=(i-1)\alpha$ (see also Lemma 3.2.11). Indeed, since
$W=UR$ and $a_{i+1}=Wa_i$, it follows that $W=U_{i+1}R_i$ and
$a_{i+1}=U_{i+1}R_ia_i$. The equality $\ell a_i=\alpha$
follows\break

\vskip-12pt

\rightskip0pt

\noindent
from the fact that $R_i$ and $U_{i+1}$ preserve the corresponding
angular measures.

By induction, $\gamma$ is homotopic in $T$ to $b\cup a'_{n+1}$. Hence,
$(nf-l-1)[c]=[\gamma]=[b]+(k+1)[c]=(e_P+k+1)[c]$.\hfill$_\blacksquare$

\medskip

{\bf3.2.15.}~Propositions 2.4.9, 3.2.14, and 3.2.6 have the following

\medskip

{\bf Corollary 3.2.16.} {\sl Under the quadrangle conditions\/
{\rm3.1.17,} the Euler number $e$ of the disc bundle\/ $M(n,l,k,p)$
over an orientable closed surface\/ $\Sigma$ equals\/ $4(nf-k-l-2)$
if\/ $n$ is odd and\/ $2(nf-k-l-2)$ if\/ $n$ is even. Moreover,
$2(\chi+e)\equiv3\tau\mod8n$ for odd\/ $n$ and\/
$2(\chi+e)\equiv3\tau\mod4n$ for even\/ $n$, where\/ $\chi$ is the
Euler characteristic of\/ $\Sigma$ and\/ $\tau$ stands for the Toledo
invariant of the representation of\/ $\pi_1\Sigma$ associated to\/
$M(n,l,k,p)$.}\hfill$_\blacksquare$

\medskip

{\bf Remark 3.2.17.} {\sl The rational Euler number of the disc bundle
over the turnover orbifold\/ $\Bbb S^2(n,n,2)$ provided by the
quadrangle equals\/ $f-\frac{k+l+2}n$. The corresponding Toledo
invariant is\/ $\frac{2t}{3n}-1$ {\rm(}see~Proposition\/ {\rm3.2.6} for
the definition of\/ $t${\rm).}}\hfill$_\blacksquare$

\medskip

\centerline{\bf3.3.~Some interesting examples.}

\medskip

In this subsection, we list some explicit examples of complex
hyperbolic disc bundles of the type $M(n,l,k,p)$ obtained with
straightforward computer calculations. All we need is to find
parameters $n,l,k,p$ satisfying the quadrangle conditions 3.1.17, a
task easily achievable with the help of any computational tool. Our
particular program has a wide margin of error, thus guaranteeing the
validity of our results in the sense that we definitely prefer to lose
an existing example rather than to get a doubtful one. In what follows,
$e$, $g$, $\chi$, and $\tau$ stand respectively for the Euler number of
the disc bundle $M(n,l,k,p)\to\Sigma$, the genus of $\Sigma$, the Euler
characteristic of $\Sigma$, and the Toledo invariant of the
corresponding representation $\pi_1\Sigma\to\PU$.

\smallskip

We have tested all $n\le1001$ and, for every $n$, all possible values
of $l,k,p$. The $\Bbb C$-Fuchsian examples have been discarded by a
direct use of Toledo's rigidity theorem [Tol] (see the beginning of
Introduction). The results are as follows:

\smallskip

$\bullet$ There is no example for $n<9$ or for $n=11,12$. For any other
$n\le1001$, there exists at least one example. The total number of
examples is 308359. There are exactly 89546 examples with integer
Toledo invariant. So, we obtain the first examples of discrete and
faithful representations $\pi_1\Sigma\to\PU$ with fractional Toledo
invariant.

$\bullet$ For every example, both triangles $\Delta(C,M_1,M_2)$ and
$\Delta(S_2,M_2,M_1)$ are elliptic, $o(z,Wz,\varphi'z)=0$,
$f=o(y,Uy,\varphi^{-1}y)$, and $p+f=2$. The last equality can be seen
as some tiny evidence supporting the complex hyperbolic variant of the
GLT-conjecture (see below).

$\bullet$ For every $\Sigma$ with $\chi\Sigma<0$ and for an arbitrary
even integer $\tau$ subject to the Toledo necessary condition
$|\tau|\le|\chi|$, a complex hyperbolic disc bundle was constructed in
[GKL]. Therefore, each of our examples with integer $\tau$ provides a
couple of nonhomeomorphic complex hyperbolic disc bundles over the same
$\Sigma$ and with the same $\tau$. This implies that there exist
discrete and faithful representations $\pi_1\Sigma\to\PU(2,1)$ that lie
in the same connected component of the space of representations [Xia]
but not in the same connected component of the space of discrete and
faithful representations.

$\bullet$ Every example satisfies the inequalities $\tau<0$ and
$\frac12\chi<e<0$. All previously known examples, including those
constructed in [GKL], satisfy the inequalities
$\chi\le e\le\frac12\chi$.

$\bullet$ Every example satisfies the equality $2(\chi+e)=3\tau$. This
equality is a necessary condition for the existence of a holomorphic
section of the bundle (in the sense that there exist a disc bundle
structure on $M$ and a smooth holomorphic surface $\Sigma\subset M$
that intersects every fibre exactly once). It suggests the following
conjecture: Every (or at least one) $M(n,l,k,p)$ possesses a disc
bundle structure admitting a holomorphic section.

$\bullet$ The following table contains all examples with extreme values
of $g$, $e$, and $e/\chi$ :

\newpage

\vskip9pt

\noindent
\hskip59pt\vbox{\offinterlineskip\hrule height1pt\halign{\strut
#&\vrule width1pt#\tabskip=0.3em&#\hfil&\vrule#&\hfil#&\vrule#&\hfil
#&\vrule#&\hfil#&\vrule#&\hfil#&\vrule#&\hfil#\hfil&\vrule
width1pt#\tabskip=0pt\cr
&&\quad Manifold\phantom{\vrule height10pt}&&$g\,$&&$\chi\
\,$&&$e\,\,$&&$\tau\ \ \,$&&Comment&\cr
\noalign{\hrule height1pt}
&&$M(10,6,3,1)$\phantom{\vrule height10pt}&&$4$&&$-6$&&$-2$&&$-5\frac13$&&minimal $g$, maximal $e$&\cr
\noalign{\hrule}
&&$M(9,4,4,1)$\phantom{\vrule height10pt}&&$6$&&$-10$&&$-4$&&$-9\frac13$&&next to minimal $g$, next to maximal $e$&\cr
\noalign{\hrule}
&&$M(9,5,3,1)$\phantom{\vrule height10pt}&&$6$&&$-10$&&$-4$&&$-9\frac13$&&next to minimal $g$, next to maximal $e$&\cr
\noalign{\hrule}
&&$M(9,6,2,1)$\phantom{\vrule height10pt}&&$6$&&$-10$&&$-4$&&$-9\frac13$&&next to minimal $g$, next to maximal $e$&\cr
\noalign{\hrule}
&&$M(14,7,7,1)$\phantom{\vrule height10pt}&&$6$&&$-10$&&$-4$&&$-9\frac13$&&next to minimal $g$, next to maximal $e$&\cr
\noalign{\hrule}
&&$M(14,8,6,1) $\phantom{\vrule height10pt}&&$6$&&$-10$&&$-4$&&$-9\frac13$&&next to minimal $g$, next to maximal $e$&\cr
\noalign{\hrule}
&&$M(14,9,5,1)$\phantom{\vrule height10pt}&&$6$&&$-10$&&$-4$&&$-9\frac13$&&next to minimal $g$, next to maximal $e$&\cr
\noalign{\hrule}
&&$M(14,10,4,1)$\phantom{\vrule height10pt}&&$6$&&$-10$&&$-4$&&$-9\frac13$&&next to minimal $g$, next to maximal $e$&\cr
\noalign{\hrule}
&&$M(14,11,3,1)$\phantom{\vrule height10pt}&&$6$&&$-10$&&$-4$&&$-9\frac13$&&next to minimal $g$, next to maximal $e$&\cr
\noalign{\hrule}
&&$M(14,0,0,2)$\phantom{\vrule height10pt}&&$6$&&$-10$&&$-4$&&$-9\frac13$&&next to minimal $g$, next to maximal $e$&\cr
\noalign{\hrule}
&&$M(16,0,0,2)$\phantom{\vrule height10pt}&&$7$&&$-12$&&$-4$&&$-10\frac23$&&next to  maximal $e$&\cr
\noalign{\hrule}
&&$M(18,0,0,2)$\phantom{\vrule height10pt}&&$8$&&$-14$&&$-4$&&$-12$&&next to  maximal $e$&\cr
\noalign{\hrule}
&&$M(20,0,0,2)$\phantom{\vrule height10pt}&&$9$&&$-16$&&$-4$&&$-13\frac13$&&next to  maximal $e$&\cr
\noalign{\hrule}
&&$M(22,0,0,2)$\phantom{\vrule height10pt}&&$10$&&$-18$&&$-4$&&$-14\frac23$&&next to  maximal $e$&\cr
\noalign{\hrule}
&&$M(24,0,0,2)$\phantom{\vrule height10pt}&&$11$&&$-20$&&$-4$&&$-16$&&next to  maximal $e$&\cr
\noalign{\hrule}
&&$M(26,0,0,2)$\phantom{\vrule height10pt}&&$12$&&$-22$&&$-4$&&$-17\frac13$&&next to  maximal $e$&\cr
\noalign{\hrule}
&&$M(28,0,0,2)$\phantom{\vrule height10pt}&&$13$&&$-24$&&$-4$&&$-18\frac23$&&minimal $e/\chi=\frac16$, next to  maximal $e$&\cr
}\hrule height1pt}

\vskip7pt

In principle, some examples in the above table could be isometric, for
instance, those with $n=14$. However, this is not the case as it is
easy to show that these bundles provide distinct points in the space
$\Cal T_n$ of discrete and faithful representations
$\varrho:\pi_1\Sigma\to\PU$ modulo conjugation by $\PU$. Note that all
examples with $n=14$ have the same known discrete invariants. We
believe that those with $p=1$ are in the same connected component of
$\Cal T_{14}$. Moreover, we expect that there exist new Toledo-like
discrete invariants that distinguish the cases $p=1$ and $p=2$ thus
showing that the corresponding manifolds provide points in different
connected components of $\Cal T_{14}$. Such a behaviour should not be
specific for $n=14$.

The fact that the congruences in Corollary 3.2.16 turn out to be
equalities in our explicit examples is a sort of evidence in support of
the complex hyperbolic GLT-conjecture analogous to the real hyperbolic
one stated in [GLT] : A disc bundle over a closed orientable surface
admits a complex hyperbolic structure if and only if $|e|\le|\chi|$ and
$\chi<0$.

$\bullet$ The two tables below display all 55 examples that satisfy the
inequality $\frac13\chi\le e$ and, therefore, according to [Kui], admit
a real hyperbolic structure:

\newpage

\vskip9pt

\noindent
\hskip43pt\vbox{\offinterlineskip\hrule height1pt\halign{\strut#&\vrule
width1pt#\tabskip=0.3em&#\hfil&\vrule#&\hfil#&\vrule#&\hfil#&\vrule
#&\hfil#&\vrule#&\hfil#&\vrule width1pt#\tabskip=0pt\cr
&&\quad Manifold\phantom{\vrule height10pt}&&$g\,$&&$\chi\
\,$&&$e\,\,$&&$\tau\ \ \,$&\cr
\noalign{\hrule height1pt}
&&$M(10,6,3,1)$\phantom{\vrule height10pt}&&$4$&&$-6$&&$-2$&&$-5\frac13$&\cr
\noalign{\hrule}
&&$M(16,0,0,2)$\phantom{\vrule height10pt}&&$7$&&$-12$&&$-4$&&$-10\frac23$&\cr
\noalign{\hrule}
&&$M(18,0,0,2)$\phantom{\vrule height10pt}&&$8$&&$-14$&&$-4$&&$-12$&\cr
\noalign{\hrule}
&&$M(20,0,0,2)$\phantom{\vrule height10pt}&&$9$&&$-16$&&$-4$&&$-13\frac13$&\cr
\noalign{\hrule}
&&$M(22,0,0,2)$\phantom{\vrule height10pt}&&$10$&&$-18$&&$-4$&&$-14\frac23$&\cr
\noalign{\hrule}
&&$M(22,1,0,2)$\phantom{\vrule height10pt}&&$10$&&$-18$&&$-6$&&$-16$&\cr
\noalign{\hrule}
&&$M(24,0,0,2)$\phantom{\vrule height10pt}&&$11$&&$-20$&&$-4$&&$-16$&\cr
\noalign{\hrule}
&&$M(24,1,0,2)$\phantom{\vrule height10pt}&&$11$&&$-20$&&$-6$&&$-17\frac13$&\cr
\noalign{\hrule}
&&$M(26,0,0,2)$\phantom{\vrule height10pt}&&$12$&&$-22$&&$-4$&&$-17\frac13$&\cr
\noalign{\hrule}
&&$M(26,1,0,2)$\phantom{\vrule height10pt}&&$12$&&$-22$&&$-6$&&$-18\frac23$&\cr
\noalign{\hrule}
&&$M(28,0,0,2)$\phantom{\vrule height10pt}&&$13$&&$-24$&&$-4$&&$-18\frac23$&\cr
\noalign{\hrule}
&&$M(28,1,0,2)$\phantom{\vrule height10pt}&&$13$&&$-24$&&$-6$&&$-20$&\cr
\noalign{\hrule}
&&$M(28,1,1,2)$\phantom{\vrule height10pt}&&$13$&&$-24$&&$-8$&&$-21\frac13$&\cr
\noalign{\hrule}
&&$M(28,2,0,2)$\phantom{\vrule height10pt}&&$13$&&$-24$&&$-8$&&$-21\frac13$&\cr
\noalign{\hrule}
&&$M(17,0,0,2)$\phantom{\vrule height10pt}&&$14$&&$-26$&&$-8$&&$-22\frac23$&\cr
\noalign{\hrule}
&&$M(30,1,1,2)$\phantom{\vrule height10pt}&&$14$&&$-26$&&$-8$&&$-22\frac23$&\cr
\noalign{\hrule}
&&$M(30,2,0,2)$\phantom{\vrule height10pt}&&$14$&&$-26$&&$-8$&&$-22\frac23$&\cr
\noalign{\hrule}
&&$M(32,1,1,2)$\phantom{\vrule height10pt}&&$15$&&$-28$&&$-8$&&$-24$&\cr
\noalign{\hrule}
&&$M(32,2,0,2)$\phantom{\vrule height10pt}&&$15$&&$-28$&&$-8$&&$-24$&\cr
\noalign{\hrule}
&&$M(19,0,0,2)$\phantom{\vrule height10pt}&&$16$&&$-30$&&$-8$&&$-25\frac13$&\cr
\noalign{\hrule}
&&$M(34,1,1,2)$\phantom{\vrule height10pt}&&$16$&&$-30$&&$-8$&&$-25\frac13$&\cr
\noalign{\hrule}
&&$M(34,2,1,2)$\phantom{\vrule height10pt}&&$16$&&$-30$&&$-10$&&$-26\frac23$&\cr
\noalign{\hrule}
&&$M(34,3,0,2)$\phantom{\vrule height10pt}&&$16$&&$-30$&&$-10$&&$-26\frac23$&\cr
\noalign{\hrule}
&&$M(36,1,1,2)$\phantom{\vrule height10pt}&&$17$&&$-32$&&$-8$&&$-26\frac23$&\cr
\noalign{\hrule}
&&$M(36,2,1,2)$\phantom{\vrule height10pt}&&$17$&&$-32$&&$-10$&&$-28$&\cr
\noalign{\hrule}
&&$M(21,0,0,2)$\phantom{\vrule height10pt}&&$18$&&$-34$&&$-8$&&$-28$&\cr
\noalign{\hrule}
&&$M(38,2,1,2)$\phantom{\vrule height10pt}&&$18$&&$-34$&&$-10$&&$-29\frac13$&\cr
\noalign{\hrule}
&&$M(40,2,2,2)$\phantom{\vrule height10pt}&&$19$&&$-36$&&$-12$&&$-32$&\cr
}\hrule height1pt}
\hskip40pt\vbox{\offinterlineskip\hrule height1pt\halign{\strut#&\vrule
width1pt#\tabskip=0.3em&#\hfil&\vrule#&\hfil#&\vrule#&\hfil#&\vrule
#&\hfil#&\vrule#&\hfil#&\vrule width1pt#\tabskip=0pt\cr
&&\quad Manifold\phantom{\vrule height10pt}&&$g\,$&&$\chi\
\,$&&$e\,\,$&&$\tau\ \ \,$&\cr
\noalign{\hrule height1pt}
&&$M(40,3,1,2)$\phantom{\vrule height10pt}&&$19$&&$-36$&&$-12$&&$-32$&\cr
\noalign{\hrule}
&&$M(23,0,0,2)$\phantom{\vrule height10pt}&&$20$&&$-38$&&$-8$&&$-30\frac23$&\cr
\noalign{\hrule}
&&$M(23,1,0,2)$\phantom{\vrule height10pt}&&$20$&&$-38$&&$-12$&&$-33\frac13$&\cr
\noalign{\hrule}
&&$M(42,2,2,2)$\phantom{\vrule height10pt}&&$20$&&$-38$&&$-12$&&$-33\frac13$&\cr
\noalign{\hrule}
&&$M(44,2,2,2)$\phantom{\vrule height10pt}&&$21$&&$-40$&&$-12$&&$-34\frac23$&\cr
\noalign{\hrule}
&&$M(25,0,0,2)$\phantom{\vrule height10pt}&&$22$&&$-42$&&$-8$&&$-33\frac13$&\cr
\noalign{\hrule}
&&$M(25,1,0,2)$\phantom{\vrule height10pt}&&$22$&&$-42$&&$-12$&&$-36$&\cr
\noalign{\hrule}
&&$M(46,3,2,2)$\phantom{\vrule height10pt}&&$22$&&$-42$&&$-14$&&$-37\frac13$&\cr
\noalign{\hrule}
&&$M(27,0,0,2)$\phantom{\vrule height10pt}&&$24$&&$-46$&&$-8$&&$-36$&\cr
\noalign{\hrule}
&&$M(27,1,0,2)$\phantom{\vrule height10pt}&&$24$&&$-46$&&$-12$&&$-38\frac23$&\cr
\noalign{\hrule}
&&$M(52,3,3,2)$\phantom{\vrule height10pt}&&$25$&&$-48$&&$-16$&&$-42\frac23$&\cr
\noalign{\hrule}
&&$M(29,1,0,2)$\phantom{\vrule height10pt}&&$26$&&$-50$&&$-12$&&$-41\frac13$&\cr
\noalign{\hrule}
&&$M(29,1,1,2)$\phantom{\vrule height10pt}&&$26$&&$-50$&&$-16$&&$-44$&\cr
\noalign{\hrule}
&&$M(29,2,0,2)$\phantom{\vrule height10pt}&&$26$&&$-50$&&$-16$&&$-44$&\cr
\noalign{\hrule}
&&$M(31,1,1,2)$\phantom{\vrule height10pt}&&$28$&&$-54$&&$-16$&&$-46\frac23$&\cr
\noalign{\hrule}
&&$M(31,2,0,2)$\phantom{\vrule height10pt}&&$28$&&$-54$&&$-16$&&$-46\frac23$&\cr
\noalign{\hrule}
&&$M(33,1,1,2)$\phantom{\vrule height10pt}&&$30$&&$-58$&&$-16$&&$-49\frac13$&\cr
\noalign{\hrule}
&&$M(35,1,1,2)$\phantom{\vrule height10pt}&&$32$&&$-62$&&$-16$&&$-52$&\cr
\noalign{\hrule}
&&$M(35,2,1,2)$\phantom{\vrule height10pt}&&$32$&&$-62$&&$-20$&&$-54\frac23$&\cr
\noalign{\hrule}
&&$M(37,1,1,2)$\phantom{\vrule height10pt}&&$34$&&$-66$&&$-16$&&$-54\frac23$&\cr
\noalign{\hrule}
&&$M(37,2,1,2)$\phantom{\vrule height10pt}&&$34$&&$-66$&&$-20$&&$-57\frac13$&\cr
\noalign{\hrule}
&&$M(41,2,2,2)$\phantom{\vrule height10pt}&&$38$&&$-74$&&$-24$&&$-65\frac13$&\cr
\noalign{\hrule}
&&$M(41,3,1,2)$\phantom{\vrule height10pt}&&$38$&&$-74$&&$-24$&&$-65\frac13$&\cr
\noalign{\hrule}
&&$M(43,2,2,2)$\phantom{\vrule height10pt}&&$40$&&$-78$&&$-24$&&$-68$&\cr
\noalign{\hrule}
&&$M(45,2,2,2)$\phantom{\vrule height10pt}&&$42$&&$-82$&&$-24$&&$-70\frac23$&\cr
\noalign{\hrule}
&&$M(47,3,2,2)$\phantom{\vrule height10pt}&&$44$&&$-86$&&$-28$&&$-76$&\cr
\noalign{\hrule}
&&$M(53,3,3,2)$\phantom{\vrule height10pt}&&$50$&&$-98$&&$-32$&&$-86\frac23$&\cr
}\hrule height1pt}

\vskip7pt

$\bullet$ The following tables show all examples for $n=101$ ($g=98$,
$\chi=-194$). They illustrate a typical behaviour of $l$, $k$, $p$, and
the Euler number for some fixed $n$ :

\newpage

\vskip9pt

\noindent
\hskip10pt\vbox{\offinterlineskip\hrule height1pt\halign{\strut#&\vrule
width1pt#\tabskip=0.3em&\hfil#\hfil&\vrule width1pt#\tabskip=0pt\cr
&&$e=-96$, $\tau=-193\frac13$\phantom{\vrule height10pt}&\cr
\noalign{\hrule height1pt}
&&$M(101,62,61,1)$\phantom{\vrule height10pt}&\cr
&&$M(101,63,60,1)$\phantom{\vrule height10pt}&\cr
&&$M(101,64,59,1)$\phantom{\vrule height10pt}&\cr
&&$M(101,65,58,1)$\phantom{\vrule height10pt}&\cr
&&$M(101,66,57,1)$\phantom{\vrule height10pt}&\cr
&&$M(101,67,56,1)$\phantom{\vrule height10pt}&\cr
&&$M(101,68,55,1)$\phantom{\vrule height10pt}&\cr
&&$M(101,69,54,1)$\phantom{\vrule height10pt}&\cr
&&$M(101,70,53,1)$\phantom{\vrule height10pt}&\cr
&&$M(101,71,52,1)$\phantom{\vrule height10pt}&\cr
&&$M(101,72,51,1)$\phantom{\vrule height10pt}&\cr
&&$M(101,73,50,1)$\phantom{\vrule height10pt}&\cr
&&$M(101,74,49,1)$\phantom{\vrule height10pt}&\cr
&&$M(101,75,48,1)$\phantom{\vrule height10pt}&\cr
&&$M(101,76,47,1)$\phantom{\vrule height10pt}&\cr
&&$M(101,77,46,1)$\phantom{\vrule height10pt}&\cr
&&$M(101,78,45,1)$\phantom{\vrule height10pt}&\cr
}\hrule height1pt}
\noindent
\hskip8pt\vbox{\offinterlineskip\hrule height1pt\halign{\strut#&\vrule
width1pt#\tabskip=0.3em&\hfil#\hfil&\vrule width1pt#\tabskip=0pt\cr
&&$e=-96$, $\tau=-193\frac13$\phantom{\vrule height10pt}&\cr
\noalign{\hrule height1pt}
&&$M(101,79,44,1)$\phantom{\vrule height10pt}&\cr
&&$M(101,80,43,1)$\phantom{\vrule height10pt}&\cr
&&$M(101,81,42,1)$\phantom{\vrule height10pt}&\cr
&&$M(101,82,41,1)$\phantom{\vrule height10pt}&\cr
&&$M(101,83,40,1)$\phantom{\vrule height10pt}&\cr
&&$M(101,84,39,1)$\phantom{\vrule height10pt}&\cr
&&$M(101,85,38,1)$\phantom{\vrule height10pt}&\cr
&&$M(101,86,37,1)$\phantom{\vrule height10pt}&\cr
&&$M(101,87,36,1)$\phantom{\vrule height10pt}&\cr
&&$M(101,88,35,1)$\phantom{\vrule height10pt}&\cr
&&$M(101,89,34,1)$\phantom{\vrule height10pt}&\cr
&&$M(101,90,33,1)$\phantom{\vrule height10pt}&\cr
&&$M(101,91,32,1)$\phantom{\vrule height10pt}&\cr
&&$M(101,92,31,1)$\phantom{\vrule height10pt}&\cr
&&$M(101,93,30,1)$\phantom{\vrule height10pt}&\cr
&&$M(101,94,29,1)$\phantom{\vrule height10pt}&\cr
&&$M(101,95,28,1)$\phantom{\vrule height10pt}&\cr
}\hrule height1pt}
\noindent
\hskip8pt\vbox{\offinterlineskip\hrule height1pt\halign{\strut#&\vrule
width1pt#\tabskip=0.3em&\hfil#\hfil&\vrule width1pt#\tabskip=0pt\cr
&&$e=-96$, $\tau=-193\frac13$\phantom{\vrule height10pt}&\cr
\noalign{\hrule height1pt}
&&$M(101,96,27,1)$\phantom{\vrule height10pt}&\cr
&&$M(101,97,26,1)$\phantom{\vrule height10pt}&\cr
&&$M(101,98,25,1)$\phantom{\vrule height10pt}&\cr
&&$M(101,11,11,2)$\phantom{\vrule height10pt}&\cr
&&$M(101,12,10,2)$\phantom{\vrule height10pt}&\cr
&&$M(101,13,9,2)$\phantom{\vrule height10pt}&\cr
&&$M(101,14,8,2)$\phantom{\vrule height10pt}&\cr
&&$M(101,15,7,2)$\phantom{\vrule height10pt}&\cr
&&$M(101,16,6,2)$\phantom{\vrule height10pt}&\cr
&&$M(101,17,5,2)$\phantom{\vrule height10pt}&\cr
&&$M(101,18,4,2)$\phantom{\vrule height10pt}&\cr
&&$M(101,19,3,2)$\phantom{\vrule height10pt}&\cr
&&$M(101,20,2,2)$\phantom{\vrule height10pt}&\cr
&&$M(101,21,1,2)$\phantom{\vrule height10pt}&\cr
&&$M(101,22,0,2)$\phantom{\vrule height10pt}&\cr
}\hrule height1pt}
\noindent
\hskip8pt\vbox{\offinterlineskip\hrule height1pt\halign{\strut#&\vrule
width1pt#\tabskip=0.3em&\hfil#\hfil&\vrule width1pt#\tabskip=0pt\cr
&&$e=-92$, $\tau=-190\frac23$\phantom{\vrule height10pt}&\cr
\noalign{\hrule height1pt}
&&$M(101,96,26,1)$\phantom{\vrule height10pt}&\cr
&&$M(101,97,25,1)$\phantom{\vrule height10pt}&\cr
&&$M(101,98,24,1)$\phantom{\vrule height10pt}&\cr
&&$M(101,11,10,2)$\phantom{\vrule height10pt}&\cr
&&$M(101,12,9,2)$\phantom{\vrule height10pt}&\cr
&&$M(101,13,8,2)$\phantom{\vrule height10pt}&\cr
&&$M(101,14,7,2)$\phantom{\vrule height10pt}&\cr
&&$M(101,15,6,2)$\phantom{\vrule height10pt}&\cr
&&$M(101,16,5,2)$\phantom{\vrule height10pt}&\cr
&&$M(101,17,4,2)$\phantom{\vrule height10pt}&\cr
&&$M(101,18,3,2)$\phantom{\vrule height10pt}&\cr
&&$M(101,19,2,2)$\phantom{\vrule height10pt}&\cr
&&$M(101,20,1,2)$\phantom{\vrule height10pt}&\cr
&&$M(101,21,0,2)$\phantom{\vrule height10pt}&\cr
}\hrule height1pt}

\vskip9pt

\noindent
\hskip40pt\vbox{\offinterlineskip\hrule height1pt\halign{\strut#&\vrule
width1pt#\tabskip=0.3em&\hfil#\hfil&\vrule width1pt#\tabskip=0pt\cr
&&$e=-88$, $\tau=-188$\phantom{\vrule height10pt}&\cr
\noalign{\hrule height1pt}
&&$M(101,10,10,2)$\phantom{\vrule height10pt}&\cr
&&$M(101,11,9,2)$\phantom{\vrule height10pt}&\cr
&&$M(101,12,8,2)$\phantom{\vrule height10pt}&\cr
&&$M(101,13,7,2)$\phantom{\vrule height10pt}&\cr
&&$M(101,14,6,2)$\phantom{\vrule height10pt}&\cr
&&$M(101,15,5,2)$\phantom{\vrule height10pt}&\cr
&&$M(101,16,4,2)$\phantom{\vrule height10pt}&\cr
}\hrule height1pt}
\noindent
\hskip37pt\vbox{\offinterlineskip\hrule height1pt\halign{\strut#&\vrule
width1pt#\tabskip=0.3em&\hfil#\hfil&\vrule width1pt#\tabskip=0pt\cr
&&$e=-84$, $\tau=-185\frac13$\phantom{\vrule height10pt}&\cr
\noalign{\hrule height1pt}
&&$M(101,10,9,2)$\phantom{\vrule height10pt}&\cr
&&$M(101,11,8,2)$\phantom{\vrule height10pt}&\cr
&&$M(101,12,7,2)$\phantom{\vrule height10pt}&\cr
}\hrule height1pt}
\noindent
\hskip37pt\vbox{\offinterlineskip\hrule height1pt\halign{\strut#&\vrule
width1pt#\tabskip=0.3em&\hfil#\hfil&\vrule width1pt#\tabskip=0pt\cr
&&$e=-80$, $\tau=-182\frac23$\phantom{\vrule height10pt}&\cr
\noalign{\hrule height1pt}
&&$M(101,9,9,2)$\phantom{\vrule height10pt}&\cr
}\hrule height1pt}

\bigskip

\centerline{\bf4.~Appendix: technical tools, elementary algebraic and
geometric background}

\medskip

This section contains proofs of all facts that were just stated in the
previous sections. The core of the method we use here is a
coordinate-free approach in complex hyperbolic geometry. It originates
mainly from the works [San] and [HSa]. With this approach, we
significantly reduce the amount of calculations involved in the proofs.
It happened to be convenient to work in a slightly extended generality,
using extensively the positive part of $\Bbb P$. This either does not
alter or indeed simplifies the proofs. It is worthwhile mentioning that
the pseudo-Riemannian geometry of the positive part is important {\it
per se,} as it is nothing but the geometry of complex geodesics; it is
studied more systematically in [AGr1].

Behind several results of this section there is an important and
curious geometrical interpretation that we plan to discuss in other
articles. It is related to some kind of parallel displacement. (The
explicit formulae can be found in [AGr1] and [AGr2].) This concerns the
slice identification, the angle between cotranchal bisectors, and the
K\"ahler primitive, for instance.

Sometimes, we prove well-known results. In these cases, the direct
references do not simplify the exposition since, usually, we need
something developed in the proof.

For the convenience of the reader unfamiliar with complex hyperbolic
geometry, we have collected in Subsection 4.1 a series of simple
well-known facts. Most of them can be found in [Gol] or [San]. Some of
the notation and definitions introduced in this subsection differ from
the standard ones.

\medskip

\centerline{\bf4.1.~Preliminaries: hermitian metric, complex geodesics,
geodesics, $\Bbb R$-planes, and bisectors}

\medskip

{\bf4.1.1.}~Fix, once and for all, a three-dimensional $\Bbb C$-vector
space $V$ equipped with a hermitian form $\langle-,-\rangle$ of
signature $++-$. (In order to deal with the Poincar\'e disc, one can
take $\dim_\Bbb CV=2$ and signature $+-$.) Depending on the context, we
will use the elements of $V$ to denote the points in the complex
projectivization $\Bbb P$ of $V$. In general, $\Bbb PW$ denotes the
complex projectivization of $W\subset V$.

As in 2.1.1, we define
$$\Bbb B:=\big\{p\in\Bbb P\mid\langle p,p\rangle<0\big\},\qquad
\partial_\infty\Bbb B:=\big\{p\in\Bbb P\mid\langle p,p\rangle=0\big\},
\qquad\overline{\Bbb B}:=\Bbb B\cup\partial_\infty\Bbb B.$$

We recall that every linear map $\varphi\in\Lin_\Bbb C(\Bbb Cp,V)$ can
be regarded as a tangent vector $t_{\varphi}\in\T_p\Bbb P$ by defining
$$t_{\varphi}f:=\frac d{d\varepsilon}\Big|_{\varepsilon=0}\hat
f\big(p+\varepsilon\varphi(p)\big)$$
for a local smooth function $f$ on $\Bbb P$ and its lift $\hat f$ to
$V$. Indeed, for small $\varepsilon\in\Bbb R$, we have a smooth curve
$c(\varepsilon)\in\Bbb P$ whose lift to $V$ is given by
$c_0(\varepsilon):=p+\varepsilon\varphi(p)$. The above $t_{\varphi}f$
is nothing but $\dot c(0)f$. In this way,
$${\T}_p\Bbb P\simeq{\Lin}_\Bbb C(\Bbb Cp,V/\Bbb Cp).$$

Let $p$ be nonisotropic, i.e., $p\notin\partial_\infty\Bbb B$. We have
the orthogonal decomposition
$$V=\Bbb Cp\oplus p^\perp,\qquad v=\pi'[p]v+\pi[p]v,$$
where the orthogonal projections
$$\pi'[p]v:=\frac{\langle v,p\rangle}{\langle p,p\rangle}p\in\Bbb
Cp,\qquad\pi[p]v:=v-\frac{\langle v,p\rangle}{\langle p,p\rangle}p\in
p^\perp$$
do not depend on the choice of a representative $p\in V$. The natural
identification
$${\T}_p\Bbb P\simeq\langle-,p\rangle p^\perp$$
provides a nondegenerate hermitian form on the tangent space
$\T_p\Bbb P$ :
$$\langle v_p,w_p\rangle:=-\langle p,p\rangle\langle
v,w\rangle,\eqno{\bold{(4.1.2)}}$$
where $v,w\in p^\perp$ and, for a chosen representative $p\in V$,
$$v_p:=\langle-,p\rangle v\in{\T}_p\Bbb P.\eqno{\bold{(4.1.3)}}$$
Obviously, this hermitian form depends smoothly on a nonisotropic $p$.

\medskip

{\bf Lemma 4.1.4.} {\sl Let $c:[a,b]\to\Bbb P$ be a smooth curve and
let $c_0:[a,b]\to V$ be a smooth lift of $c$.
If~$c(t_0)\notin\partial_\infty\Bbb B$, then
$\dot c(t_0)=\big\langle-,c_0(t_0)\big\rangle
\displaystyle\frac{\pi\big[c(t_0)\big]\dot
c_0(t_0)}{\big\langle c_0(t_0),c_0(t_0)\big\rangle}$
is the tangent vector to\/ $c$ at\/ $c(t_0)$.}

\medskip

{\bf Proof.} Denoting
$k:=\displaystyle\frac{\big\langle\dot
c_0(t_0),c_0(t_0)\big\rangle}{\big\langle
c_0(t_0),c_0(t_0)\big\rangle}$,
we have
$$\dot c(t_0)f=\frac d{dt}\Big|_{t=t_0}\hat f\big(c_0(t)\big)=\frac
d{d\varepsilon}\Big|_{\varepsilon=0}\hat f\big(c_0(t_0)+\varepsilon\dot
c_0(t_0)\big)=$$
$$=\frac d{d\varepsilon}\Big|_{\varepsilon=0}\hat
f\Big(c_0(t_0)+\varepsilon\pi'\big[c(t_0)\big]\dot
c_0(t_0)+\varepsilon\pi\big[c(t_0)\big]\dot c_0(t_0)\Big)=\frac
d{d\varepsilon}\Big|_{\varepsilon=0}\hat f\Big((1+\varepsilon
k)c_0(t_0)+\varepsilon\pi\big[c(t_0)\big]\dot c_0(t_0)\Big)=$$
$$=\frac d{d\varepsilon}\Big|_{\varepsilon=0}\hat
f\Big(c_0(t_0)+\frac\varepsilon{1+\varepsilon k}\pi\big[c(t_0)\big]\dot
c_0(t_0)\Big)=\frac d{d\varepsilon}\Big|_{\varepsilon=0}\hat
f\Big(c_0(t_0)+\varepsilon\pi\big[c(t_0)\big]\dot
c_0(t_0)\Big).\eqno{_\blacksquare}$$

{\bf4.1.5.~Tance.} For $p_1,p_2\in\Bbb P$, we define
$$\ta(p_1,p_2):=\frac{\langle p_1,p_2\rangle\langle
p_2,p_1\rangle}{\langle p_1,p_1\rangle\langle p_2,p_2\rangle}.$$
When one of $p_1,p_2$ is isotropic, our convention is that
$\ta(p_1,p_2):=+\infty$ if $\langle p_1,p_2\rangle\neq0$ and
$\ta(p_1,p_2):=1$ if $\langle p_1,p_2\rangle=0$. The concept of tance
can be used for expressing distances (see, for example,
Corollary~4.1.18 and Lemma 4.1.8).

Let $p_1,p_2\in\Bbb B$. Then $\ta(p_1,p_2)\ge1$ and $\ta(p_1,p_2)=1$ if
and only if $p_1=p_2$. Given subsets $X,Y\subset\Bbb B$, define the
tance between $X$ and $Y$ as
$$\ta(X,Y):=\inf\big\{\ta(x,y)\mid x\in X,\ y\in Y\big\}.$$

\smallskip

{\bf4.1.6.~Polar points to projective lines.} For every projective line
$\L\subset\Bbb P$, there exists a unique point $p\in\Bbb P$ such that
$\L=\Bbb Pp^\perp$. We call $p$ the {\it polar\/} point to $\L$. By
definition, the {\it signature\/} of $\L$ is simply that of the
hermitian form on the $\Bbb C$-vector subspace in $V$ corresponding to
$\L$. The signature of $\Bbb Pp^\perp$ is respectively $++$, $+-$, $+0$
in the cases $p\in\Bbb B$, $p\notin\overline{\Bbb B}$,
$p\in\partial_\infty\Bbb B$. Denote by
$$\L{\wr}p_1,p_2{\wr}:=\Bbb P(\Bbb Cp_1+\Bbb Cp_2)$$
the projective line spanned by two distinct points $p_1,p_2\in\Bbb P$.

\medskip

{\bf Lemma 4.1.7 {\rm[Gol, p.~100]}.} {\sl Let\/ $p_1,p_2\in\Bbb P$ be
distinct. Then\/ $\Bbb Pp_1^\perp$ and\/ $\Bbb Pp_2^\perp$ intersect in
the polar point\/ $p$ to\/ $\L{\wr}p_1,p_2{\wr}$. Moreover,

\smallskip

$\bullet$ $p\in\partial_\infty\Bbb B$ if and only if\/
$\ta(p_1,p_2)=1$;

$\bullet$ $p\in\Bbb B$ if and only if\/ $\ta(p_1,p_2)<1$ and\/
$p_1,p_2\notin\overline{\Bbb B}$.}

\medskip

{\bf Proof.} The first assertion is obvious. Sylvester's Criterion
[KoM, p.~113] and the definition of $\ta(p_1,p_2)$ imply that the
signature of $\L{\wr}p_1,p_2{\wr}$ is $+0$ if and only if
$\ta(p_1,p_2)=1$. This proves the second assertion. Suppose that
$p\in\Bbb B$. Then $p_1,p_2\notin\overline{\Bbb B}$ because
$p_1,p_2\in p^\perp$. Applying Sylvester's Criterion to the Gram matrix
of $p_1,p_2,p$, we infer that $\ta(p_1,p_2)<1$. Conversely,
$\ta(p_1,p_2)<1$ and $p_1,p_2\notin\overline{\Bbb B}$ imply
$p\in\Bbb B$ by the same criterion.\hfill$_\blacksquare$

\medskip

The {\it complex geodesics\/} we deal with in Sections 1--3 are by
definition of the form $\Bbb Pp^\perp\cap\overline{\Bbb B}$ with
positive~$p$. Two distinct complex geodesics are said to be {\it
ultraparallel\/} if they do not intersect. By~Lem\-ma~4.1.7, this is
equivalent to the inequality $\ta(p_1,p_2)>1$ for their polar points
$p_1$ and $p_2$.

The tance between a point and a complex geodesic is given in the
following

\medskip

{\bf Lemma 4.1.8 {\rm[Gol, p.~197, Theorem 6.1.1]}.} {\sl Let\/
$p\notin\overline{\Bbb B}$ and let\/ $q\in\Bbb B$. Then\/
$\ta(q,\Bbb Pp^\perp\cap\Bbb B)=1-\ta(p,q)$.}

\medskip

{\bf Proof.} It follows from
$\big\langle\pi[p]q,\pi[p]q\big\rangle=\big\langle\pi[p]q,q\big\rangle=
\langle q,q\rangle\big(1-\ta(p,q)\big)<0$
that $\pi[p]q\in\Bbb Pp^\perp\cap\Bbb B$ and that
$\ta\big(\pi[p]q,q\big)=1-\ta(p,q)$. We choose an orthonormal basis
$e_1,e_2,e_3$ in $V$ such that $e_1$ and $e_2$ represent $p$ and
$\pi[p]q$, respectively. We can assume that $q=ce_1+e_2$, $|c|<1$.
Every point in $\Bbb Pp^\perp\cap\Bbb B$ can be written in the form
$p(z):=e_2+ze_3$, $|z|<1$. It remains to observe that
$$\ta\big(p(z),q\big)=\frac1{\big(-1+|z|^2\big)\big(|c|^2-1\big)}=
\frac1{\big(1-|z|^2\big)\big(1-|c|^2\big)}.\eqno{_\blacksquare}$$

{\bf4.1.9.~Geodesics and $\Bbb R$-planes.} Let $S\subset V$ be an
$\Bbb R$-vector subspace such that the hermitian form is real and
nondegenerate on $S$. It is immediate that
$\Bbb CS\simeq\Bbb C\otimes_\Bbb RS$. Therefore,
$\Bbb PS=\Bbb P_\Bbb RS$. If~$\dim_\Bbb RS=2$, we call
$$\G S:=\Bbb PS=\Bbb P_\Bbb R^1S\simeq\Bbb S^1$$
an {\it extended geodesic\/} (or simply a {\it geodesic\/}) [San]. If
$\dim_\Bbb RS=3$, we call
$$\R S:=\Bbb PS=\Bbb P_\Bbb R^2S$$
an {\it extended\/} $\Bbb R$-{\it plane\/} (or simply an $\Bbb R$-{\it
plane\/}) [Gol, p.~80]. When writing $\G S$ or $\R S$ (or even
$\Bbb PS$) for a geodesic or real plane, we assume that $S$ is an
$\Bbb R$-vector subspace such that the hermitian form is real and
nondegenerate on $S$. Every geodesic $\G S$ spans its projective line
$\L\G S:=\Bbb P(\Bbb CS)$. By definition, the {\it signature\/} of
$\G S$ is that of $\L\G S$. It coincides with the signature of the
(hermitian) form on $S$ and can be $++$ or $+-$.

\medskip

{\bf Remark 4.1.10.} {\sl Let\/ $\L$ be a projective line and let\/
$p\in\L$ be nonisotropic. Then there exists a unique\/ $q\in\L$
orthogonal to\/ $p$, i.e., $\langle p,q\rangle=0$. If\/ $\G$ is a
geodesic such that\/ $p\in\G\subset\L$, then\/
$q\in G$.}\hfill$_\blacksquare$

\medskip

Sylvester's Criterion implies the following

\medskip

{\bf Remark 4.1.11.} {\sl Suppose that\/ $\ta(g_1,g_2)\ne0,1$ for
some\/ $g_1,g_2\in\Bbb P$. Then there exists a unique extended
geodesic\/ $\G{\wr}g_1,g_2{\wr}$ containing\/ $g_1$ and\/
$g_2$.}\hfill$_\blacksquare$

\medskip

When using the notation $\G{\wr}g_1,g_2{\wr}$, we assume by default
that $\ta(g_1,g_2)\ne0,1$.

\smallskip

{\bf4.1.12.}~In order to write down the equation of the geodesic
$\G{\wr}g_1,g_2{\wr}$, we put
$$b(x,g_1,g_2):=\frac{\langle g_1,x\rangle\langle x,g_2\rangle}{\langle
g_1,g_2\rangle}-\frac{\langle g_2,x\rangle\langle x,g_1\rangle}{\langle
g_2,g_1\rangle}.$$
Note that $b(x,g_1,g_2)$ does not depend on the choice of
representatives $g_1,g_2\in V$.

\medskip

{\bf Lemma 4.1.13.} {\sl The geodesic\/ $\G{\wr}g_1,g_2{\wr}$ is given
in\/ $\L{\wr}g_1,g_2{\wr}$ by the equation\/ $b(x,g_1,g_2)=0$.}

\medskip

{\bf Proof.} $\G{\wr}g_1,g_2{\wr}=\G S$ for a suitable real subspace
$S\subset V$. The equality $b(x,g_1,g_2)=0$ holds for every $x\in S$
since the hermitian form is real on $S$ and we can take $g_1,g_2\in S$.

Conversely, suppose that
$\displaystyle\frac{\langle g_1,x\rangle\langle x,g_2\rangle}{\langle
g_1,g_2\rangle}\in\Bbb R$
for some $x\in\L{\wr}g_1,g_2{\wr}$. We can assume that $g_1,g_2\in S$
and that $x=g_1+cg_2$ for some $c\in\Bbb C$. Then
$0\ne\langle g_1,g_2\rangle\in\Bbb R$ and
$$0\equiv\langle g_1,x\rangle\langle x,g_2\rangle\equiv c\langle
g_1,g_1\rangle\langle g_2,g_2\rangle+\overline c\langle
g_1,g_2\rangle\langle g_1,g_2\rangle\mod\Bbb R.$$
If $c\notin\Bbb R$, then
$\langle g_1,g_1\rangle\langle g_2,g_2\rangle=\langle
g_1,g_2\rangle\langle g_1,g_2\rangle$.
In other words, $\ta(g_1,g_2)=1$.\hfill$_\blacksquare$

\medskip

{\bf4.1.14.}~Every geodesic $\G$ of signature $+-$ possesses exactly
two isotropic points $v_1,v_2\in\G$ called the {\it vertices\/} of
$\G$. There is a canonical way to describe such a geodesic in terms of
its vertices:

\medskip

{\bf Lemma 4.1.15 {\rm[Gol, p.~155]}.} {\sl Let\/ $\G$ be a geodesic of
signature\/ $+-$, let\/ $v_1,v_2$ stand for the vertices of\/ $\G$, and
let\/ $g\in\G$ be a point of signature\/ $\pm$. Then we can choose
representatives\/ $v_1,v_2,g\in V$ such that\/ $g=v_1+v_2$ and\/
$\left(\smallmatrix0&\pm\frac12\\\pm\frac12&0\endsmallmatrix\right)$ is
the Gram matrix of\/ $v_1,v_2$. The formula\/ $g(t):=t^{-1}v_1+tv_2$,
$t>0$, describes all points in\/ $\G$ of signature\/ $\pm$. For all\/
$t$, we have\/ $\big\langle g(t),g(t)\big\rangle=\pm1$.}

\medskip

{\bf Proof.} We can always choose representatives $v_1,v_2,g\in V$
providing the required Gram matrix, $g\in\Bbb Rv_1+\Bbb Rv_2$, and
$\langle g,g\rangle=\pm1$. Replacing $g$ by $-g$ if necessary, we
obtain $g=t^{-1}v_1+tv_2$ for a suitable $t>0$. It remains to take
$t^{-1}v_1,tv_2$ for $v_1,v_2$.\hfill$_\blacksquare$

\medskip

{\bf Lemma 4.1.16.} {\sl If\/ $p,g_1,g_2\in\Bbb B$, then\/
$\displaystyle\frac{\langle g_1,p\rangle\langle p,g_2\rangle}{\langle
g_1,g_2\rangle}$
cannot be real nonnegative.}

\medskip

{\bf Proof.} It suffices to show that
$\Re\displaystyle\frac{\langle g_1,p\rangle\langle
p,g_2\rangle}{\langle g_1,g_2\rangle}<0$
assuming that $g_1\ne g_2\ne p$. Taking the negative point $\pi[q]p$ in
place of $p$, where $q\notin\overline{\Bbb B}$ denotes the polar point
to $\L{\wr}g_1,g_2{\wr}$, we can assume that $p\in\L{\wr}g_1,g_2{\wr}$.
Choose representatives such that $\langle g_i,g_i\rangle=-1$,
$\langle g_1,g_2\rangle=a>1$, and $p=g_1+cg_2$ for some $c\in\Bbb C$.
Since $\langle p,p\rangle<0$, we have $2a\Re c<1+|c|^2$. Now,
$$\Re\big(\langle g_1,p\rangle\langle p,g_2\rangle\big)=(a^2+1)\Re
c-a\big(1+|c|^2\big)<
\frac{(1-a^2)\big(1+|c|^2\big)}{2a}<0.\eqno{_\blacksquare}$$

{\bf4.1.17.~Length of geodesics.} Let $\G[g_1,g_2]\subset\Bbb B$ denote
the oriented segment of geodesic joining the points $g_1,g_2\in\Bbb B$.
By Lemma 4.1.15, we can parameterize a lift of $\G[g_1,g_2]$ to $V$ as
$c_0(t)=e^{-t}v_1+e^tv_2$, $t\in[0,a]$, for some $a\ge0$. It is easy to
see that $\dot c_0(t)=-e^{-t}v_1+e^tv_2$ is orthogonal to $c_0(t)$.
By~Lemma~4.1.4, $\big\langle\dot c(t),\dot c(t)\big\rangle=1$.
Therefore,
$\ell\G[g_1,g_2]=\displaystyle\int\limits_0^a\sqrt{\big\langle\dot
c(t),\dot c(t)\big\rangle}\,dt=a$.
A straightforward calculation shows that
$\ta(g_1,g_2)=\displaystyle\Big(\frac{e^{-a}+e^a}2\Big)^2$. We obtain
the following result.

\medskip

{\bf Corollary 4.1.18.} {\sl For\/ $g_1,g_2\in\Bbb B$, we have\/
$\cosh^2\big(\dist(g_1,g_2)\big)=\ta(g_1,g_2)$.}\hfill$_\blacksquare$

\medskip

Now, we can see that the metric on $\Bbb B$ coincides, up to the scale
factor of $4$, with the one introduced in [Gol, p.~77].

\smallskip

\vskip5pt

\noindent
$\hskip265pt\vcenter{\hbox{\epsfbox{Pictures.2}}}$

\rightskip200pt

\vskip-222pt

{\bf4.1.19.~Bisectors.} Let $\G S$ be a geodesic and let $f$ stand for
the polar point to its projective line $\L\G S$. The~projective cone
$\B$ over $\G S$ with vertex $f$ is said to be an {\it extended
bisector\/} (or simply a {\it bisector\/}). The nonisotropic point $f$
is the {\it focus\/} of $\B$. Clearly, $f$ is the only singular point
in $\B$. The geodesic $\G S$ is the {\it real spine\/} of $\B$, the
projective line $\L\G S$ is the {\it complex spine\/} of $\B$, and the
projective line $\L{\wr}g,f{\wr}$ is the {\it slice\/} of $\B$
containing $g\in\G S$. In this way, we obtain the {\it slice
decomposition\/} of~$\B$ : every point in $\B$ different from the focus
belongs to a unique slice of $\B$. By Remark 4.1.10, the polar point to
every slice of $\B$ belongs to the real spine of $\B$. The~vertices of
the real spine of $\B$ (if exist) are the {\it vertices\/} of $\B$.

We can also describe the bisector $\B$ as $\B=\Bbb PW$, where
$W:=S+\Bbb Cf$ is a $4$-dimensional real subspace in $V$. This
description immediately provides the {\it meridional decomposition\/}
of $\B$ : Let us fix some representative $f\in V$. For~every
$0\ne c\in\Bbb C$, the hermitian form on $S':=S+\Bbb Rcf$ is real and
nondegenerate. The $\Bbb R$-plane $\R S'\subset\B$ is said\break

\vskip-12pt

\rightskip0pt

\noindent
to be a {\it meridian\/} of $\B$. Two different meridians intersect
exactly in the focus and in the real spine. Obviously, every point of a
bisector belongs to some meridian. Such a meridian is unique if the
point is different from the focus and does not lie on the real spine.

It is easy to see that the intersection of a meridian of $\B$ with the
slice $\L{\wr}g,f{\wr}$ of $\B$ is a geodesic (unless $g$ is a vertex
of the real spine of $\B$) that contains $g$.

It is also possible to define a bisector as the hypersurface
equidistant (equitant) from two different points. In this article, we
do not need such a definition.

\medskip

{\bf4.1.20.~Remark.} {\sl Our definition of a bisector differs slightly
from the common one. A standardly defined bisector is, in our terms,
$\B\cap\overline{\Bbb B}$, where\/ $\B$ is a bisector with positive
focus. In terms of\/ {\rm[Gol, p.~248],} the bisectors in our sense are
extors, however, not every extor is a bisector in our sense.}

\medskip

We denote by $\B{\wr}g_1,g_2{\wr}$ the bisector with the real spine
$\G{\wr}g_1,g_2{\wr}$. When using the notation $\B{\wr}g_1,g_2{\wr}$,
we assume by default that $\ta(g_1,g_2)\ne0,1$.

Lemma 4.1.13 implies the following

\medskip

{\bf Proposition 4.1.21.} {\sl The bisector\/ $\B{\wr}g_1,g_2{\wr}$ is
given in\/ $\Bbb P$ by the equation\/
$b(x,g_1,g_2)=0$.}\hfill$_\blacksquare$

\medskip

\vskip4pt

\noindent
$\hskip219pt\vcenter{\hbox{\epsfbox{Pictures.4}}}$

\rightskip240pt

\vskip-106pt

{\bf Lemma 4.1.22 {\rm[Gol, p.~165, Theorem 5.2.4]}.} {\sl Let\/
$p_1,p_2\in\Bbb P$ be such that\/ $\ta(p_1,p_2)\ne0,1$ and let\/
$S_i:=\Bbb Pp_i^\perp$, $i=1,2$. Then there exists a unique bisector\/
$\B$ such that\/ $S_1$ and\/ $S_2$ are slices of\/ $\B$. The real spine
of such\/ $\B$ is $\G{\wr}p_1,p_2{\wr}$.}

\medskip

{\bf Proof.} The intersection $f$ of the $S_i$'s has to be the focus of
the bisector $\B$ in question. By~Lemma 4.1.7, $\L{\wr}p_1,p_2{\wr}$
is necessarily the complex spi\-ne of $\B$. Since the polar point to
every\break

\vskip-12pt

\rightskip0pt

\noindent
slice of $\B$ belongs to the real spine of $\B$, this real spine
contains the points $p_1$ and $p_2$. By Remark 4.1.11, the real spine
of $\B$ has to equal $\G{\wr}p_1,p_2{\wr}$. By Remark 4.1.10, the point
$g_i\in\L{\wr}p_1,p_2{\wr}$ orthogonal to $p_i$ belongs to
$\G{\wr}p_1,p_2{\wr}$, $i=1,2$. In other words, $S_i=\L{\wr}g_i,f{\wr}$
is a slice of $\B{\wr}p_1,p_2{\wr}$.\hfill$_\blacksquare$

\medskip

In particular, two ultraparallel complex geodesics $S_1$ and $S_2$ are
slices of a unique bisector $\B$. Obviously, the focus $f$ of such $\B$
is positive. The projective cone over $\G[g_1,g_2]$ with vertex $f$ is
the segment $\B[S_1,S_2]$ of $\B$, where $g_i$ stands for the
intersection of $S_i$ with the real spine of $\B$. Using the
parameterization given in Lemma 4.1.15 of (positive) polar points to
the slices of $\B$, it is easy to show that there exists a complex
geodesic $M$ called the {\it middle slice\/} of the segment
$\B[S_1,S_2]$ such that the reflection in $M$ exchanges $S_1$ and
$S_2$.

\medskip

{\bf Lemma 4.1.23.} {\sl Suppose that the Gram matrix of\/
$p_1,p_2\in V$ has the form\/
$\left(\smallmatrix1&t\\t&1\endsmallmatrix\right)$ with\/ $t>1$. Then\/
$m:=\displaystyle\frac{p_1+p_2}{\sqrt{2t+2}}$ is the polar point to the
middle slice\/ $M$ of the segment\/ $\B[S_1,S_2]$ and\/
$\langle m,m\rangle=1$, where\/ $S_i:=\Bbb Pp_i^\perp$, $i=1,2$.}

\medskip

{\bf Proof.} It is easy to see that $\langle m,m\rangle=1$. We have
(see (2.1.3) for the definition of $R(m)$)
$$R(m)p_1=2\langle p_1,m\rangle m-p_1=\frac{2\langle
p_1,p_1+p_2\rangle}{2t+2}(p_1+p_2)-p_1=p_2.\eqno{_\blacksquare}$$

For different points $v_1,v_2\in\partial_\infty\Bbb B$ and
$p\notin\partial_\infty\Bbb B$, the $\eta$-invariant [Gol, p.~231] is
defined as
$$\eta(v_1,v_2,p):=\frac{\langle v_1,p\rangle\langle
p,v_2\rangle}{\langle v_1,v_2\rangle\langle
p,p\rangle}.\eqno{\bold{(4.1.24)}}$$
The tance between a point and a bisector can be expressed in terms of
the $\eta$-invariant:

\medskip

{\bf Lemma 4.1.25 {\rm[San, p.~97, Example 2]}.} {\sl Let\/
$\B{\wr}v_1,v_2{\wr}$ be a bisector with vertices\/
$v_1,v_2\in\partial_\infty\Bbb B$ and let\/
$p\in\Bbb B\setminus\B{\wr}v_1,v_2{\wr}$. Then\/
$\ta\big(p,\B{\wr}v_1,v_2{\wr}\cap\Bbb B\big)=
1-\Re\eta(v_1,v_2,p)+\big|\eta(v_1,v_2,p)\big|$.}

\medskip

{\bf Proof.} We can assume that $\langle p,p\rangle=-1$,
$\langle v_1,v_2\rangle=\frac12$, and $r:=|z_1|=|z_2|$, where
$z_i:=\langle p,v_i\rangle$. By~Lemma 4.1.15, the polar points to the
slices of signature $+-$ are parameterized by $g(t):=t^{-1}v_1+tv_2$,
$t>0$, with $\big\langle g(t),g(t)\big\rangle=1$. By Lemma 4.1.8,
$$\ta\big(p,\Bbb Pg(t)^\perp\cap\Bbb
B\big)=1-\ta\big(p,g(t)\big)=1+2\Re(\overline
z_1z_2)+t^{-2}r^2+t^2r^2.$$
This function has a unique minimum at $t=1$. It remains to observe that
$\eta(v_1,v_2,p)=-2\overline z_1z_2$.\hfill$_\blacksquare$

\medskip

\centerline{\bf4.2.~Bisectors}

\medskip

In this subsection, we explicitly describe the tangent space (at
generic points) to subspaces of the form $\Bbb PW$, where $W\subset V$
is an $\Bbb R$-vector subspace. This tool provides a handy description
of the normal vector to a bisector. Then we prove some known facts
concerning transversalities. Finally, we observe that a meridional
curve (see 2.1.1) depends continuously on its segment of bisector and
initial point.

\smallskip

{\bf4.2.1.}~Let $W\subset V$ be an $\Bbb R$-vector subspace. We define
$$d:=\min_{0\ne w\in W}\dim_\Bbb R\Bbb Cw\cap W,\qquad D:=\{w\in
W\mid\dim_\Bbb R\Bbb Cw\cap W=d\}.$$
It is easy to show that $D$ is open in $W$, that $\Bbb PD$ is smooth
and open in $\Bbb PW$, and that $\dim_\Bbb R\Bbb PD=\dim_\Bbb RW-d$.

\medskip

{\bf Lemma 4.2.2.} {\sl Let\/ $p\in D$ and let\/
$\varphi\in\Lin_\Bbb C(\Bbb Cp,V)$. Then\/ $t_\varphi\in\T_p\Bbb PW$ if
and only if\/ $\varphi(p)\in\Bbb Cp+W$.}

\medskip

{\bf Proof.} Suppose that $\varphi(p)\in W$. For small
$\varepsilon\in\Bbb R$, we have a smooth curve
$c(\varepsilon)\in\Bbb PD$ whose lift to $V$ is given by
$c_0(\varepsilon):=p+\varepsilon\varphi(p)\in W$. Therefore, the
tangent vector $t_\varphi$ to $c$ belongs to $\T_p\Bbb PD=\T_p\Bbb PW$.
Such tangent vectors form the whole $\T_p\Bbb PD$ since
$$\dim_\Bbb R\Bbb PD=\dim_\Bbb RW-d=\dim_\Bbb RW/\Bbb Cp\cap
W.\eqno{_\blacksquare}$$

{\bf Corollary 4.2.3.} {\sl Let\/ $\L$ be a projective line, let\/
$p\in\L$ be nonisotropic, and let\/ $q\in\L$ be orthogonal to\/ $p$.
Then\/ $\T_p\L=\langle-,p\rangle\Bbb Cq$. Let\/ $\Bbb PS$ be an
$\Bbb R$-plane or geodesic and let\/ $p\in S$ be nonisotropic. Then\/
$\T_p\Bbb PS=\langle-,p\rangle(p^\perp\cap S)$.}\hfill$_\blacksquare$

\medskip

As it is easy to see, Lemma 2.1.13 follows from the second part of
Corollary 4.2.3.

\smallskip

{\bf4.2.4.~Orthogonal projective lines.} Let $\L_1$ be a projective
line and let $p\in\L_1$ be nonisotropic. Using Corollary 4.2.3, it is
easy to show that there exists a unique projective line $\L_2$ passing
through $p$ and orthogonal to $\L_1$ in the sense that $\T_p\L_1$ and
$\T_p\L_2$ are orthogonal. Obviously, $\L_2=\Bbb Pp_1^\perp$ for some
$p_1\in\L_1$. This means that two distinct projective lines
intersecting in a nonisotropic point $p$ are orthogonal at $p$ if and
only if their polar points are orthogonal. In this case,
$\T_p\Bbb P=\T_p\L_1\oplus\T_p\L_2$. In~particular, the complex spine
of a bisector is orthogonal to every slice of signature different from
$+0$. This implies that the intersection of such a slice with a
meridian is a geodesic orthogonal to the real spine.

\smallskip

{\bf4.2.5.~Some simple transversalities.} In order to show a few simple
facts concerning transversalities, we need the following

\medskip

{\bf Lemma 4.2.6.} {\sl Let\/ $p\in\partial_\infty\Bbb B$ and let\/
$\varphi\in\Lin_\Bbb C(\Bbb Cp,V)$. Then\/
$t_\varphi\in\T_p\partial_\infty\Bbb B$ if and only if\/
$\Re\big\langle\varphi(p),p\big\rangle=0$.}

\medskip

{\bf Proof.} We take $q\in V$ such that $\langle p,q\rangle\ne0$. In
some neighbourhood of $p$, the $3$-sphere $\partial_\infty\Bbb B$ is
given by the equation
$h(x):=\displaystyle\frac{\langle x,x\rangle}{\langle q,x\rangle\langle
x,q\rangle}=0$.
The tangent vector $t_\varphi$ is simply $\dot c(0)$, where
$c(\varepsilon)$ is a smooth curve in $\Bbb P$ whose lift to $V$ is
given by $c_0(\varepsilon)=p+\varepsilon\varphi(p)$. It follows from
$$\dot c(0)h=\frac d{d\varepsilon}\Big|_{\varepsilon=0}h
\big(p+\varepsilon\varphi(p)\big)=\frac
d{d\varepsilon}\Big|_{\varepsilon=0}\frac{\big\langle
p+\varepsilon\varphi(p),p+\varepsilon\varphi(p)\big\rangle}{\big\langle
q,p+\varepsilon\varphi(p)\big\rangle\big\langle
p+\varepsilon\varphi(p),q\big\rangle}=
\frac{\big\langle\varphi(p),p\big\rangle+\big\langle
p,\varphi(p)\big\rangle}{\langle p,q\rangle\langle q,p\big\rangle}$$
that $\dot c(0)h=0$ if and only if
$\Re\big\langle\varphi(p),p\big\rangle=0$. It remains to observe that
such $t_\varphi$'s form a $3$-dimensional $\Bbb R$-vector subspace in
$\T_p\Bbb P$.\hfill$_\blacksquare$

\medskip

{\bf Remark 4.2.7 {\rm(folklore)}.} {\sl Every projective line of
signature $+-$ is transversal to $\partial_\infty\Bbb B$.}

\medskip

{\bf Proof.} Let $p\notin\overline{\Bbb B}$ and let
$q\in\Bbb Pp^\perp\cap\partial_\infty\Bbb B$. Fix some point
$q'\in\Bbb Pp^\perp$ different from $q$. So,
$\Bbb Pp^\perp=\L{\wr}q,q'{\wr}$ and we can assume that
$0\ne\langle q,q'\rangle\in\Bbb R$. Suppose that
$t_\varphi\in\T_q\Bbb Pp^\perp\cap\T_q\partial_\infty\Bbb B$, where
$\varphi\in\Lin_\Bbb C(\Bbb Cq,V)$. By Lemma 4.2.2,
$\varphi(q)\in p^\perp=\Bbb Cq+\Bbb Cq'$. We can change $\varphi$ by
adding to $\varphi(q)$ an arbitrary element in $\Bbb Cq$ since this
does not alter $t_\varphi$. Hence, we can assume that
$\varphi(q)\in\Bbb Cq'$. By Lemma 4.2.6,
$\varphi(q)\in\Bbb Riq'$.\hfill$_\blacksquare$

\medskip

{\bf Corollary 4.2.8.} {\sl Every bisector\/ $\B$ is transversal to\/
$\partial_\infty\Bbb B$.}

\medskip

{\bf Proof.} Denote by $\G S$ the real spine of $\B$ and by $f$, the
focus of $\B$. Let $p\in\B\cap\partial_\infty\Bbb B$. Obviously,
$p\ne f$. We can assume $p=s+cf$ for some $0\ne s\in S$ and
$c\in\Bbb C$. It follows from $p\in\partial_\infty\Bbb B$ that
$\langle s,s\rangle+|c|^2\langle f,f\rangle=0$. Take some
$s'\in S\setminus\Bbb Rs$.

Suppose that $t_\varphi\in\T_p\B\cap\T_p\partial_\infty\Bbb B$, where
$\varphi\in\Lin_\Bbb C(\Bbb Cp,V)$. By Lemmas 4.2.2 and 4.2.6,
$\varphi(p)\in\Bbb Cs+\Bbb Rs'+\Bbb Cf$ and
$\Re\big\langle\varphi(p),p\big\rangle=0$. Adding to $\varphi(p)$ an
arbitrary element from $\Bbb Cp$ does not change $t_\varphi$.
Therefore, we can assume that $\varphi(p)\in\Bbb Rs'+\Bbb Cf$.

If $c=0$, then $p=s$ is isotropic and $\langle s',s\rangle\ne0$. In
this case, the equality $\Re\big\langle\varphi(p),p\big\rangle=0$
implies $\varphi(p)\in\Bbb Cf$ and we conclude that
$\dim_\Bbb R(\T_p\B\cap\T_p\partial_\infty\Bbb B)\le2$. If $c\ne0$,
then $s$ is nonisotropic and we can assume that $c=1$ and take
$s'\in S\setminus\Bbb Rs$ orthogonal to $s$. In this case, the equality
$\Re\big\langle\varphi(p),p\big\rangle=0$ implies
$\varphi(p)\in\Bbb Rs'+\Bbb Rif$.\hfill$_\blacksquare$

\medskip

{\bf4.2.9.}~We now introduce the equation of the tangent space to the
bisector $\B{\wr}g_1,g_2{\wr}$. Define
$$t(v,p,g_1,g_2):=\frac{\langle g_1,v\rangle\langle
p,g_2\rangle+\langle g_1,p\rangle\langle v,g_2\rangle}{\langle
g_1,g_2\rangle}-\frac{\langle g_2,v\rangle\langle p,g_1\rangle+\langle
g_2,p\rangle\langle v,g_1\rangle}{\langle g_2,g_1\rangle}.$$
Note that $t(v,p,g_1,g_2)$ does not depend on the choice of
representatives $v,p,g_1,g_2\in V$ that provides the same linear map
$\langle-,p\rangle v$.

\medskip

{\bf Lemma 4.2.10.} {\sl Let\/ $p\in\B{\wr}g_1,g_2{\wr}$ be different
from the focus\/ $f$ of\/ $\B{\wr}g_1,g_2{\wr}$ and let\/
$\varphi\in\Lin_\Bbb C(\Bbb Cp,V)$. Then\/
$t_\varphi\in\T_p\B{\wr}g_1,g_2{\wr}$ if and only if\/
$t\big(\varphi(p),p,g_1,g_2\big)=0$.}

\medskip

{\bf Proof.} The equality $t\big(\varphi(p),p,g_1,g_2\big)=0$ remains
valid if we choose other representatives
$g_1,g_2,p\allowmathbreak\in V$. Hence, we can assume that
$g_1,g_2\in S$ and $p\in W$, where $\G{\wr}g_1,g_2{\wr}=\G S$,
$\B{\wr}g_1,g_2{\wr}=\Bbb PW$, and~$W:=S+\Bbb Cf$. By Lemma 4.2.2,
$t_\varphi\in\T_p\B{\wr}g_1,g_2{\wr}$ if and only if
$\varphi(p)\in\Bbb Cp+W$.

Suppose that $\varphi(p)\in\Bbb Cp+W$. Let us show that
$t\big(\varphi(p),p,g_1,g_2\big)=0$. Clearly, $t(s,p,g_1,g_2)=0$ for
all $s\in S$. Since $\langle g_i,f\rangle=0$ for $i=1,2$, we obtain
$t(x,p,g_1,g_2)=0$ for all $x\in\Bbb Cf$. For every $c\in\Bbb C$,
$$t(cp,p,g_1,g_2)=\frac{\overline c\langle g_1,p\rangle\langle
p,g_2\rangle+c\langle g_1,p\rangle\langle p,g_2\rangle}{\langle
g_1,g_2\rangle}-\frac{\overline c\langle g_2,p\rangle\langle
p,g_1\rangle+c\langle g_2,p\rangle\langle p,g_1\rangle}{\langle
g_2,g_1\rangle}=(c+\overline c)b(p,g_1,g_2)=0$$
by Proposition 4.1.21. Hence, $t\big(\varphi(p),p,g_1,g_2\big)=0$.

Conversely, suppose that $t\big(\varphi(p),p,g_1,g_2\big)=0$. The point
$p$ is different from $f$. Interchanging $g_1$ and $g_2$ if necessary,
we can assume that $p=g_1+rg_2+c'f$ for some $r\in\Bbb R$ and
$c'\in\Bbb C$. As was shown, $t(\Bbb Cp+\Bbb Cf,p,g_1,g_2)=0$. It is
easy to see that $V=\Bbb Cg_2+\Bbb Cp+\Bbb Cf$. Hence, we can assume
that $\varphi(p)=cg_2$ with $c\in\Bbb C$. Since
$t(cg_2,g_2,g_1,g_2)=(c+\overline c)b(g_2,g_1,g_2)=0$ and
$\langle g_1,g_2\rangle\in\Bbb R$, we obtain
$$t(cg_2,p,g_1,g_2)=t(cg_2,g_1,g_1,g_2)=(\overline c-c)\frac{\langle
g_1,g_2\rangle\langle g_2,g_1\rangle-\langle g_1,g_1\rangle\langle
g_2,g_2\rangle}{\langle g_1,g_2\rangle}.$$
By Sylvester's Criterion, $t(cg_2,p,g_1,g_2)=0$ implies
$c\in\Bbb R$.\hfill$_\blacksquare$

\medskip

{\bf Proposition 4.2.11.} {\sl Let\/ $p\in\B{\wr}g_1,g_2{\wr}$ be
nonisotropic and different from the focus $f$ of\/
$\B{\wr}g_1,g_2{\wr}$ and let\/ $g\in\G{\wr}g_1,g_2{\wr}$ stand for the
polar point to the slice of\/ $\B{\wr}g_1,g_2{\wr}$ containing\/ $p$.
Then
$$0\ne n(p,g_1,g_2):=\langle-,p\rangle i\Big(\frac{\langle
p,g_2\rangle}{\langle g_1,g_2\rangle}g_1-\frac{\langle
p,g_1\rangle}{\langle g_2,g_1\rangle}g_2\Big)\in\langle-,p\rangle\Bbb
Cg$$
is a normal vector to\/ $\B{\wr}g_1,g_2{\wr}$ at\/ $p$.}

\medskip

{\bf Proof.} If $n(p,g_1,g_2)=0$, then $\langle p,g_i\rangle=0$ for
$i=1,2$, implying $p=f$. A contradiction.

Obviously, $n(p,g_1,g_2)\in\langle-,p\rangle f^\perp$. By Proposition
4.1.21, $b(p,g_1,g_2)=0$, implying
$n(p,g_1,g_2)\in\langle-,p\rangle p^\perp$. Clearly, $\L(p,f)$ is the
slice of $\B{\wr}g_1,g_2{\wr}$ containing $p$. Therefore,
$p^\perp\cap f^\perp=\Bbb Cg$.

Take an arbitrary vector $v_p\in\T_p\Bbb P$, $v\in p^\perp$ (see
(4.1.3) for the definition). By Lemma 4.2.10,
$v_p\in\T_p\B{\wr}g_1,g_2{\wr}$ if and only if $t(v,p,g_1,g_2)=0$.
Since
$$\big\langle n(p,g_1,g_2),v_p\big\rangle=-\langle p,p\rangle
i\Big(\frac{\langle g_1,v\rangle\langle p,g_2\rangle}{\langle
g_1,g_2\rangle}-\frac{\langle g_2,v\rangle\langle p,g_1\rangle}{\langle
g_2,g_1\rangle}\Big)$$
and $p$ is nonisotropic, the equality $t(v,p,g_1,g_2)=0$ is equivalent
to $\Re\big\langle n(p,g_1,g_2),v_p\big\rangle=0$.\hfill$_\blacksquare$

\medskip

{\bf4.2.12.}~We introduce the orientation of $\G{\wr}g_1,g_2{\wr}$ as
follows. Let $i=1,2$. If $g_i\in\Bbb B$, then define
$q_i\in\G{\wr}g_1,g_2{\wr}$ as the point orthogonal to $g_i$.
Otherwise, put $q_i:=g_i$. Denote by $a\subset\G{\wr}g_1,g_2{\wr}$ the
open arc from $q_2$ to $q_1$ that contains no point orthogonal to
$q_i$, $i=1,2$, and such that $a\cap\overline{\Bbb B}=\varnothing$.
The~{\it orientation\/} of $\G{\wr}g_1,g_2{\wr}$ is that of $a$.

\medskip

{\bf Remark 4.2.13.} {\sl Let\/ $\G{\wr}g_1,g_2{\wr}$ be a geodesic,
let\/ $p_i\in\G{\wr}g_1,g_2{\wr}$ be orthogonal to\/ $g_i$, and let\/
$p$ denote the polar point to\/ $\L{\wr}g_1,g_2{\wr}$. Suppose that\/
$p_1,p_2\notin\Bbb B$. We have the following cases\/{\rm:}

\smallskip

$\bullet$ $p\in\Bbb B$. Then\/ $g_1,g_2,p_1,p_2$ are in a cyclic order
in the circle\/ $\G{\wr}g_1,g_2{\wr}$.

$\bullet$ $p\notin\overline{\Bbb B}$. Then\/ $g_1,g_2,v_2,p_2,p_1,v_1$
are in a cyclic order in the circle\/ $\G{\wr}g_1,g_2{\wr}$, where\/
$v_1,v_2$ are the vertices of\/ $\G{\wr}g_1,g_2{\wr}$ taken in an
appropriate order. {\rm(}Note that we admit\/ $g_i=v_i$ or\/
$p_i=v_i$.{\rm)}}\hfill$_\blacksquare$

\medskip

The bisector $\B{\wr}g_1,g_2{\wr}$ is {\it oriented\/} with respect to
the orientation of its real spine $\G{\wr}g_1,g_2{\wr}$ and the natural
(complex) orientation of its slices. Using Remark 4.2.13, it is easy to
see that the introduced orientation of a bisector is compatible with
the one introduced in 2.1.1.

\smallskip

{\bf4.2.14.}~We will show that the normal vector in Proposition 4.2.11
is in fact normal to the oriented bisector $\B{\wr}g_1,g_2{\wr}$. It
follows from Proposition 4.1.21 that every bisector
$\B{\wr}g_1,g_2{\wr}$ divides $\Bbb P$ into two {\it extended\/}
half-spaces given by the inequalities
$$\Im\frac{\langle g_1,x\rangle\langle x,g_2\rangle}{\langle
g_1,g_2\rangle}\ge0,\qquad\Im\frac{\langle g_1,x\rangle\langle
x,g_2\rangle}{\langle g_1,g_2\rangle}\le0.$$

{\bf Lemma 4.2.15.} {\sl Let\/ $\B{\wr}g_1,g_2{\wr}$ be an oriented
bisector, let\/ $p_i\in\G{\wr}g_1,g_2{\wr}$ be the point orthogonal
to\/~$g_i$, $i=1,2$, and let\/ $p\in\B{\wr}g_1,g_2{\wr}$ be a
nonisotropic point different from the focus $f$ of\/
$\B{\wr}g_1,g_2{\wr}$. Denote by $\sigma_f$ the signature of $f$. Then

\smallskip

$\bullet$ the inequalities\/
$\displaystyle\Im\frac{\langle g_1,x\rangle\langle
x,g_2\rangle}{\langle g_1,g_2\rangle}\ge0$
and\/
$\displaystyle\Im\frac{\langle p_1,x\rangle\langle
x,p_2\rangle}{\langle p_1,p_2\rangle}\ge0$
are equivalent,

$\bullet$ the vector\/ $n(p,g_1,g_2)$ is normal to the oriented
bisector\/ $\B{\wr}g_1,g_2{\wr}$,

$\bullet$ the extended half-space given by\/
$\sigma_f\displaystyle\Im\frac{\langle g_1,x\rangle\langle
x,g_2\rangle}{\langle g_1,g_2\rangle}\ge0$
lies on the side of\/ $n(p,g_1,g_2)$.}

\medskip

{\bf Proof.} By Remark 4.2.13, $\G{\wr}g_1,g_2{\wr}$ and
$\G{\wr}p_1,p_2{\wr}$ coincide as oriented geodesics. By
Proposition~4.1.21, we can exclude the case
$x\in\B{\wr}g_1,g_2{\wr}=\B{\wr}p_1,p_2{\wr}$ and reduce the first
assertion to the last~two.

Fix $x\not\in\B{\wr}g_1,g_2{\wr}$. We will vary $p,g_1,g_2$ keeping the
following conditions: $p$ belongs to the bisector and satisfies
$f\ne p\notin\partial_\infty\Bbb B$; $g_1,g_2$ belong to the real spine
and satisfy $\ta(g_1,g_2)\ne0,1$. During the deformation, the
orientation of $\B{\wr}g_1,g_2{\wr}$ remains the same, the normal
vector $n(p,g_1,g_2)$ never vanishes, and
$\displaystyle\Im\frac{\langle g_1,x\rangle\langle
x,g_2\rangle}{\langle g_1,g_2\rangle}\ne0$.
Therefore, we can assume that $p\in\G{\wr}g_1,g_2{\wr}$, that $g_1,g_2$
are vertices if $\sigma_f=+1$, and that $p=g_2$ if $\sigma_f=-1$.

\smallskip

The case $\sigma_f=+1$. By Lemma 4.1.15, we can choose representatives
such that $\langle g_1,g_2\rangle=\pm\frac12$ and $p=g_1+g_2$. By
Proposition 4.2.11, $n(p,g_1,g_2)=\langle-,p\rangle i(g_1-g_2)$. By
Lemma 4.1.15, a lift of the oriented geodesic $\G{\wr}g_1,g_2{\wr}$ is
given by $c_0(t):=e^{\pm t}g_1+e^{\mp t}g_2$, $t\in\Bbb R$. Obviously,
$c_0(0)=p$. By Lemma~4.1.4, $\dot c(0)=\langle-,p\rangle(g_1-g_2)$.
Since $i\dot c(0)=n(p,g_1,g_2)$, the vector $n(p,g_1,g_2)$ is normal to
the oriented bisector $\B{\wr}g_1,g_2{\wr}$.

The curve
$n_0(t):=p+ti\langle p,p\rangle(g_1-g_2)=(1\pm ti)g_1+(1\mp ti)g_2$ is
a lift of a smooth curve $n(t)\in\Bbb P$ with $\dot n(0)=n(p,g_1,g_2)$.
It remains to observe that
$\Im\displaystyle\frac{\big\langle g_1,n(t)\big\rangle\big\langle
n(t),g_2\big\rangle}{\langle g_1,g_2\rangle}=t$.

\smallskip

The case $\sigma_f=-1$. We can choose representatives such that
$\langle g_i,g_i\rangle=1$ and $\langle g_1,g_2\rangle=r>0$.
By~Sylvester's Criterion, $r<1$. By Proposition 4.2.11,
$n(p,g_1,g_2)=\langle-,p\rangle i(\frac1rg_1-g_2)$. A lift of the
oriented geodesic $\G{\wr}g_1,g_2{\wr}$ is given by
$c_0(t):=(1-t)g_2+tg_1$, $t\in[0,1]$. (Note that
$\big\langle c_0(t),g_2\big\rangle=0$ would imply $t=\frac1{1-r}>1$.)
By Lemma 4.1.4, $\dot c(0)=\langle-,p\rangle(g_1-rg_2)$. The vector
$n(p,g_1,g_2)$ is normal to the oriented bisector $\B{\wr}g_1,g_2{\wr}$
because $i\dot c(0)=rn(p,g_1,g_2)$.

The curve $n_0(t):=g_2+ti(g_1-rg_2)=tig_1+(1-rti)g_2$ is a lift of a
smooth curve $n(t)\in\Bbb P$ with $\dot n(0)=rn(p,g_1,g_2)$. It remains
to observe that
$\Im\displaystyle\frac{\big\langle g_1,n(t)\big\rangle\big\langle
n(t),g_2\big\rangle}{\langle
g_1,g_2\rangle}=-t\Big(\frac{1-r^2}r\Big)$.\hfill$_\blacksquare$

\medskip

{\bf Lemma 4.2.16.} {\sl An\/ {\rm(}ideal\/{\rm)} meridional curve
depends continuously on its segment of bisector and initial point.}

\medskip

{\bf Proof.} Let $\B[g_1,g_2]$ be a segment of bisector with focus
$f\notin\overline{\Bbb B}$ and $g_1,g_2\in\Bbb B$, let $S_1$ denote the
initial slice of $\B[g_1,g_2]$, i.e., the one containing $g_1$, let
$p_1\in\G{\wr}g_1,g_2{\wr}$ be orthogonal to $g_1$, and let
$m_0\in\G{\wr}g_1,g_2{\wr}$ denote the polar point to the middle slice
of $\B[g_1,g_2]$. The meridional curve $b$ of $\B[g_1,g_2]$ generated
by $p\in S_1\cap\overline{\Bbb B}$ can be parameterized by
$\G[g_1,g_2]$ as well as by $\G[p_1,m_0]$, where $m\in\G[p_1,m_0]$ is
the polar point to the middle slice of $\B[g_1,g]$ (in other words, the
parameters $g,m$ are related by the equality $g=R(m)g_1$). So, by Lemma
2.1.4,
$$b(g):=-R(m)p.$$
Since $p\ne f$, we can take $p$ in the form $p=cf+g_1$ with
$c\in\Bbb C$. We assume that $\langle g_1,g_1\rangle=-1$. Then
$\langle g,g\rangle=-1$ because $g=R(m)g_1$. It follows from
$\langle m,m\rangle>0$ and (2.1.3) that
$$\langle g,g_1\rangle=\big\langle
R(m)g_1,g_1\big\rangle=2\frac{\langle g_1,m\rangle\langle
m,g_1\rangle}{\langle m,m\rangle}-\langle g_1,g_1\rangle=2\frac{\langle
g_1,m\rangle\langle m,g_1\rangle}{\langle m,m\rangle}+1>0,$$
implying
$\sqrt{\ta(g,g_1)}\displaystyle\frac{\langle p,g_1\rangle}{\langle
g,g_1\rangle}=-1$.
We have $R(m)f=-f$ since $\langle m,f\rangle=0$. It follows that
$b(g)=-R(m)p=cf-g=\pi[g_1]p-g$. Therefore,
$$b(g)=\pi[g_1]p+\sqrt{\ta(g,g_1)}\frac{\langle p,g_1\rangle}{\langle
g,g_1\rangle}g.$$
Note that this formula does not depend on the choice of representatives
$g_1,g\in V$ and is linear in $p$.\hfill$_\blacksquare$

\newpage

\centerline{\bf4.3.~Bisectors with common slice}

\medskip

Using the formula for the normal vector to a bisector, we calculate in
this subsection the angle between two cotranchal (= having a common
slice) bisectors at a point in their common slice (for a geometric
interpretation of this angle, see [AGr1]). Then we present a numerical
criterion of transversality of such a pair of bisectors. Finally, we
prove some properties of transversal cotranchal bisectors including
Condition (3) of [AGr3, Theorem 3.5].

\medskip

{\bf Lemma 4.3.1.} {\sl Let\/
$g\in\Bbb P\setminus\partial_\infty\Bbb B$ and\/ $g_1,g_2\in\Bbb P$ be
such that\/ $\ta(g,g_i)\ne0,1$ and let\/ $p\in\Bbb Pg^\perp$ be
nonisotropic and different from the foci of\/ $\B{\wr}g,g_i{\wr}$,
$i=1,2$. Then
$$\big\langle n(p,g,g_2),n(p,g,g_1)\big\rangle=-\langle
p,p\rangle\langle g,g\rangle\frac{\langle g_1,p\rangle\langle
p,g_2\rangle}{\langle g_1,g\rangle\langle g,g_2\rangle}.$$}

{\bf Proof.} By Proposition 4.2.11,
$n(p,g,g_1)=\langle-,p\rangle i\displaystyle\frac{\langle
p,g_1\rangle}{\langle g,g_1\rangle}g$
and
$n(p,g,g_2)=\langle-,p\rangle i\displaystyle\frac{\langle
p,g_2\rangle}{\langle g,g_2\rangle}g$.
The~result follows from (4.1.2).\hfill$_\blacksquare$

\medskip

The normal vectors $n(p,g,g_1)$ and $n(p,g,g_2)$ are both tangent to
the naturally oriented projective line that passes through $p$ and is
orthogonal to the common slice $\Bbb Pg^\perp$ of the oriented
bisectors $\B{\wr}g,g_i{\wr}$, $i=1,2$. Hence, the oriented angle from
$n(p,g,g_1)$ to $n(p,g,g_2)$ is simply
$\Arg\big\langle n(p,g,g_2),n(p,g,g_1)\big\rangle$. (Recall that the
function $\Arg$ takes values in $[0,2\pi)$.)

\medskip

{\bf Corollary 4.3.2.} {\sl Let\/ $\Cal E\in\SU$ be a regular elliptic
isometry and let\/ $e_1,e_2,e_3\in\Bbb P$ be the points corresponding
to the eigenvectors of\/ $\Cal E$ such that\/ $e_1\in\Bbb B$. Denote
by\/ $E$ the complex geodesic with the polar point\/ $e_2$. Let\/ $D$
be a complex geodesic ultraparallel to\/ $E$. Then the oriented angle
from\/ $\B[E,D]$ to\/ $\B[E,\Cal ED]$ at\/ $e_1\in E$ equals\/
$\Arg(\xi_2\xi_1^{-1})$, where\/ $\xi_i$ stands for the eigenvalue of\/
$\Cal E$ corresponding to\/ $e_i$, $i=1,2$.}

\medskip

{\bf Proof.} Denote by $d$ the polar point to $D$. The points $d$ and
$e_1\in E$ cannot be orthogonal, nor can $d$ and $e_2$, since $D$ and
$E$ are ultraparallel. Therefore, we can choose representatives
$e_1,e_2\in V$ such that $\langle d,e_1\rangle>0$ and
$\langle d,e_2\rangle>0$. By Lemma 4.3.1, the angle in question equals
$$\Arg\frac{\langle d,e_1\rangle\langle e_1,\Cal Ed\rangle}{\langle
d,e_2\rangle\langle e_2,\Cal Ed\rangle}=\Arg\frac{\langle\Cal
E^{-1}e_1,d\rangle}{\langle\Cal
E^{-1}e_2,d\rangle}=\Arg\frac{\xi_1^{-1}\langle
e_1,d\rangle}{\xi_2^{-1}\langle
e_2,d\rangle}=\Arg(\xi_2\xi_1^{-1}).\eqno{_\blacksquare}$$

{\bf Criterion 4.3.3.} {\sl Let\/ $g\notin\overline{\Bbb B}$ and\/
$g_1,g_2\in\Bbb P$ be such that $\ta(g,g_i)>1$ for\/ $i=1,2$. The
bisectors\/ $\B{\wr}g,g_1{\wr}$ and\/ $\B{\wr}g,g_2{\wr}$ are
transversal along\/ $\Bbb Pg^\perp\cap\overline{\Bbb B}$ if and only if
$$\bigg|\Re\frac{\langle g_1,g_2\rangle\langle g,g\rangle}{\langle
g_1,g\rangle\langle g,g_2\rangle}-1\bigg|<\sqrt{1-\frac1{\ta(g,g_1)}}
\cdot\sqrt{1-\frac1{\ta(g,g_2)}}.$$}

{\bf Proof.} We can choose representatives such that
$\langle g,g\rangle=\langle g_1,g\rangle=\langle g,g_2\rangle=1$.
Hence, $g_i=g+b_i$ with $b_i\in g^\perp$. (Since
$b_i\in\G{\wr}g,g_i{\wr}\cap\Bbb Pg^\perp$, the slice of
$\B{\wr}g,g_i{\wr}$ containing $b_i$ is $\Bbb Pg^\perp$.) The
inequality $\ta(g,g_i)>1$ implies $\langle b_i,b_i\rangle<0$.
Normalizing the $b_i$'s, we obtain
$d_i:=\displaystyle\frac{b_i}{\sqrt{-\langle b_i,b_i\rangle}}$ with
$\langle d_i,d_i\rangle=-1$. Let $d\in g^\perp$ be orthogonal to $d_2$,
normalized so that $\langle d,d\rangle=1$. Since $d,d_2$ form an
orthonormal basis in~$g^\perp$, we have $d_1=rd+cd_2$ for suitable
$r,c\in\Bbb C$. Choosing an appropriate representative $d\in V$, we~can
assume that $r\ge0$. From $\langle d_1,d_1\rangle=-1$, we obtain
$r=\sqrt{|c|^2-1}$,
$$d_1=\sqrt{|c|^2-1}d+cd_2,\qquad|c|\ge1.$$
Every point $p\in\Bbb Pg^\perp\cap\overline{\Bbb B}$ has the form
$p=zd+d_2$ with $|z|\le1$. It follows from
$\langle p,g_2\rangle=\langle zd+d_2,g+b_2\rangle=\langle
d_2,b_2\rangle$
that $0\ne\langle p,g_2\rangle\in\Bbb R$. Note that
$$\frac{\langle p,g_1\rangle}{\sqrt{-\langle
b_1,b_1\rangle}}=\frac{\langle p,b_1\rangle}{\sqrt{-\langle
b_1,b_1\rangle}}=\langle p,d_1\rangle=\langle
zd+d_2,\sqrt{|c|^2-1}d+cd_2\rangle=z\sqrt{|c|^2-1}-\overline c.$$

Let $\varphi\in\Lin_\Bbb C(\Bbb Cp,V)$. By Lemma 4.2.10,
$t_\varphi\in\T_p\B{\wr}g,g_i{\wr}$ if and only if
$t\big(\varphi(p),p,g,g_i\big)=0$. In view of $p\in g^\perp$, this is
equivalent to
$$\big\langle g,\varphi(p)\big\rangle\langle p,g_i\rangle\in\Bbb R.$$
By Lemma 4.2.2, $t_\varphi\in\T_p\Bbb Pg^\perp$ if and only if
$\varphi(p)\in g^\perp$. Therefore, considering $t_\varphi$ modulo
$\T_p\Bbb Pg^\perp$, we~can assume that $\varphi(p)\in\Bbb Cg$.

Suppose that $\varphi(p)=ag\ne0$, $a\in\Bbb C$. Since
$0\ne\langle p,g_2\rangle\in\Bbb R$, the condition
$\big\langle g,\varphi(p)\big\rangle\langle p,g_2\rangle\in\Bbb R$ is
equivalent to $a\in\Bbb R$. Consequently, the condition
$t_\varphi\in\T_p\B{\wr}g,g_1{\wr}\cap\T_p\B{\wr}g,g_2{\wr}$ is
equivalent to the requirements that $a\in\Bbb R$ and
$\langle p,g_1\rangle\in\Bbb R$. In terms of $z$, the latter takes the
form $z\sqrt{|c|^2-1}-\overline c\in\Bbb R$. The~existence of $z$ with
$|z|\le1$ meeting the last condition means that
$\sqrt{|c|^2-1}\ge|\Im c|$. This is equivalent to $|\Re c|\ge1$. We can
see that the bisectors $\B{\wr}g,g_1{\wr}$ and $\B{\wr}g,g_2{\wr}$ are
transversal along $\Bbb Pg^\perp\cap\overline{\Bbb B}$ if and only if
$|\Re c|<1$.

It remains to observe that
$\displaystyle\sqrt{1-\frac1{\ta(g,g_i)}}=\sqrt{1-\langle
g_i,g_i\rangle}=\sqrt{-\langle b_i,b_i\rangle}$
and that
$$-c\sqrt{-\langle b_1,b_1\rangle}\sqrt{-\langle
b_2,b_2\rangle}=\langle d_1,d_2\rangle\sqrt{-\langle
b_1,b_1\rangle}\sqrt{-\langle b_2,b_2\rangle}=\langle
b_1,b_2\rangle=\langle g_1,g_2\rangle-1=\frac{\langle
g_1,g_2\rangle\langle g,g\rangle}{\langle g_1,g\rangle\langle
g,g_2\rangle}-1.\eqno{_\blacksquare}$$

In order to prove some fundamental properties of transversal cotranchal
bisectors, we need the following technical

\medskip

{\bf Remark 4.3.4.} {\sl Let\/
$v_1,v'_1,v_2,v'_2\in\partial_\infty\Bbb B$ be such that\/
$v_i\ne v'_i$ for all $i=1,2$ and let\/
$g\in\big(\G{\wr}v_1,v'_1{\wr}\cap\G{\wr}v_2,v'_2{\wr}\big)\setminus
\overline{\Bbb B}$.
Then there exist representatives\/ $v_1,v'_1,v_2,v'_2,g\in V$ such
that\/ $\langle v_i,v'_i\rangle=\frac12$ and\/ $g=v_i+v'_i$ for all\/
$i$. Such representatives satisfy the condition\/
$\Re\langle v_1,v_2\rangle\le\frac14$ or\/
$\Re\langle v'_1,v_2\rangle\le\frac14$. Moreover, the~inequality\/
$\Re\langle v_1,v_2\rangle\le2\big|\langle v_1,v_2\rangle\big|^2$
holds, the equality being valid exactly when\/ $v_1,v'_1,v_2,v'_2,g$
belong to a same projective line.}

\medskip

{\bf Proof.} By Lemma 4.1.15, we take representatives $g_1,g_2\in V$ of
$g$ and representatives $v_i,v'_i\in V$ such that
$\langle v_i,v'_i\rangle=\frac12$ and $g_i=v_i+v'_i$ for $i=1,2$. Since
$|g_1|=|g_2|$, we have $g_1=cg_2$ with $|c|=1$. It remains to take
$cv_2,cv'_2$ in place of $v_2,v'_2$.

It follows from $\langle g,v_2\rangle=\frac12$ that
$\langle v_1,v_2\rangle+\langle v'_1,v_2\rangle=\frac12$, implying that
$\Re\langle v_1,v_2\rangle\le\frac14$ or
$\Re\langle v'_1,v_2\rangle\le\frac14$.

The rest follows from Sylvester's Criterion: the determinant of the
Gram matrix of $g,v_1,v_2$ equals
$\frac12\Re\langle v_1,v_2\rangle-\big|\langle
v_1,v_2\rangle\big|^2$.\hfill$_\blacksquare$

\medskip

{\bf Lemma 4.3.5 {\rm(compare with [Hsi, p.~97, Theorem 3.3])}.} {\sl
Let\/ $\B_1$ and\/ $\B_2$ be two bisectors with positive foci and let\/
$S$ be a common slice of\/ $\B_1$ and\/ $\B_2$ of signature\/ $+-$.
Suppose that\/ $\B_1$ and\/ $\B_2$ are transversal along\/
$S\cap\overline{\Bbb B}$. Then\/ $\B_1\cap\overline{\Bbb B}$,
$\B_2\cap\overline{\Bbb B}$, and\/ $\partial_\infty\Bbb B$ are
transversal and\/
$\B_1\cap\B_2\cap\overline{\Bbb B}=S\cap\overline{\Bbb B}$.}

\medskip

{\bf Proof.} Let $g\notin\overline{\Bbb B}$ denote the polar point to
$S$ and let $v_i,v'_i\in\partial_\infty\Bbb B$ denote the vertices of
$\B_i$, $i=1,2$. By Remark 4.3.4, we can assume that $\Re u\le2|u|^2$,
$\Re u\le\frac14$, $\langle v_i,v'_i\rangle=\frac12$, and $g=v_i+v'_i$
for $i=1,2$, where $u:=\langle v_1,v_2\rangle$. Since
$\B_i=\B{\wr}g,v_i{\wr}$, we have $\Re u>0$ by Criterion 4.3.3.

Every point in $\G{\wr}g,v_i{\wr}$ but $v_i$ has the form
$$g_i(t_i):=g+(t_i-1)v_i$$
with $t_i\in\Bbb R$. It is easy to see that
$\big\langle g_i(t_i),g_i(t_i)\big\rangle=t_i$. Therefore, $g_i(t_i)$,
$t_i>0$, runs over all polar points to the slices of $\B_i$ of
signature $+-$.

Let us show that the slices of signature $+-$ of the bisectors $\B_1$
and $\B_2$ that are different from their common slice $S$ do not
intersect in $\overline{\Bbb B}$. By Lemma 4.1.7, it suffices to prove
that the inequality $\ta\big(g_1(t_1),g_2(t_2)\big)>1$ is valid for
$t_1,t_2>0$ unless $t_1=t_2=1$. Since
$$2\big\langle g_1(t_1),g_2(t_2)\big\rangle=t_1+t_2+2(t_1-1)(t_2-1)u$$
and $\big\langle g_i(t_i),g_i(t_i)\big\rangle=t_i$, we obtain
$$4t_1t_2\ta\big(g_1(t_1),g_2(t_2)\big)=
(t_1+t_2)^2+4(t_1+t_2)(t_1-1)(t_2-1)\Re u+4(t_1-1)^2(t_2-1)^2|u|^2.$$
Hence, the inequality $\ta\big(g_1(t_1),g_2(t_2)\big)>1$ is equivalent
to the inequality
$$(t_1-t_2)^2+4(t_1+t_2)(t_1-1)(t_2-1)\Re
u+4(t_1-1)^2(t_2-1)^2|u|^2>0$$
which can be rewritten in the form
$$(t_1-t_2)^2+2(t_1^2-1)(t_2^2-1)\Re
u+2(t_1-1)^2(t_2-1)^2\big(2|u|^2-\Re u\big)>0.$$
Since $\Re u>0$ and $2|u|^2\ge\Re u$, this inequality is valid if
$t_1,t_2>1$ or if $0<t_1,t_2<1$. Therefore, we can assume that
$0<t_1<1<t_2$. In this case, taking into account the inequalities
$\Re u\le\frac14$ and $2|u|^2\ge\Re u$, it suffices to observe that
$(t_1-t_2)^2>\frac12(1-t_1^2)(t_2^2-1)$. We have shown that the slices
of signature $+-$ of $\B_1$ and $\B_2$ do not intersect.

The bisectors $\B_1$ and $\B_2$ cannot have a common vertex $v$.
Otherwise, they would have the same real spine $\G{\wr}g,v{\wr}$ and
would coincide.

Note that $u$ cannot be real. This would contradict the inequalities
$0<\Re u\le\frac14$ and $\Re u\le2|u|^2$. It~is easy to see that no
vertex of one bisector belongs to the other bisector. Indeed,
$v_2\in\B_1$ means that $\big\langle g_1(t_1),v_2\big\rangle=0$ for a
suitable $t_1>0$. This can be written as $\frac12+(t_1-1)u=0$, implying
that $u$ is real. Similarly, $v'_2\in\B_1$ implies that
$\frac{t_1}2+(1-t_1)u=0$ for a suitable $t_1>0$.

Thus, $\B_1\cap\B_2\cap\overline{\Bbb B}=S\cap\overline{\Bbb B}$.

The transversality follows from Corollary 4.2.8 and Remark
4.2.7.\hfill$_\blacksquare$

\medskip

{\bf Lemma 4.3.6.} {\sl Let\/ $\B_1$ and\/ $\B_2$ be bisectors with
positive foci. Suppose that they possess a common slice\/ $S$ of
signature\/ $+-$ and that they are transversal along\/
$S\cap\overline{\Bbb B}$. Then, for any\/ $\vartheta>0$, there exists
some\/ $\varepsilon>0$ such that, for every\/ $p\in\B_2\cap\Bbb B$, the
inequality\/ $\ta(p,\B_1\cap\Bbb B)<1+\varepsilon^2$ implies the
inequality\/ $\ta(p,S\cap\Bbb B)<1+\vartheta^2$.}

\medskip

{\bf Proof.} Let $g\notin\overline{\Bbb B}$ denote the polar point to
$S$ and let $v_i,v'_i\in\partial_\infty\Bbb B$ denote the vertices of
$\B_i$, $i=1,2$. By Remark 4.3.4, we can assume that $\Re u\le2|u|^2$,
$\Re u\le\frac14$, $\langle v_i,v'_i\rangle=\frac12$, and $g=v_i+v'_i$
for $i=1,2$, where $u:=\langle v_1,v_2\rangle$. In particular, we have
$\langle v_i,g\rangle=\langle v'_i,g\rangle=\frac12$. (The equality
$\Re u=2|u|^2$ means that $\B_1$ and $\B_2$ have the same complex
spine. In this case, the fact is easier and will be proven later.)
Since $\B_i=\B{\wr}g,v_i{\wr}$, we have $\Re u>0$ by Criterion 4.3.3.

We assume that $\Re u<2|u|^2$. Therefore, $g,v_1,v_2$ are
$\Bbb C$-linearly independent. The point
$$f:=v_1+(4u-1)v_2-2ug$$
is the focus of $\B_2$ because
$\langle f,v_2\rangle=\langle f,g\rangle=0$. Define
$$k:=\sqrt{4|u|^2-2\Re u}.$$
It is easy to see that
$$\langle f,f\rangle=\langle v_1,f\rangle=(4\overline u-1)u-\overline
u=k^2,\qquad\langle f,v'_1\rangle=\langle f,g-v_1\rangle=-k^2.$$

Let us parameterize the points in $\B_2\cap\Bbb B$. By Lemma 4.1.15,
the negative points in the real spine $\G{\wr}v_2,v'_2{\wr}$ of $\B_2$
are parameterized by $g(t):=t^{-1}v_2-tv'_2$, $t>0$. Since
$\big\langle\frac fk,\frac fk\big\rangle=1$, every point in
$\B_2\cap\Bbb B$ takes the form $p(t,z):=\frac zkf+g(t)$, where $t>0$
and $|z|<1$. (The point $p(t,z)$ belongs to the slice of $\B_2$ that
contains $g(t)$.) Calculating straightforwardly, we obtain
$$\big\langle p(t,z),g\big\rangle=\big\langle
g(t),g\big\rangle=\frac{t^{-1}-t}{2},\qquad\big\langle
v_1,g(t)\big\rangle=\big\langle
v_1,t^{-1}v_2-t(g-v_2)\big\rangle=(t^{-1}+t)u-\frac t2,$$
$$\big\langle g(t),v'_1\big\rangle=\big\langle
g(t),g-v_1\big\rangle=\frac{t^{-1}}2-(t^{-1}+t)\overline
u,\qquad\big\langle g(t),g(t)\big\rangle=-1,\qquad\big\langle
p(t,z),p(t,z)\big\rangle=|z|^2-1,$$
$$\big\langle v_1,p(t,z)\big\rangle=k\overline z+(t^{-1}+t)u-\frac
t2,\qquad\big\langle
p(t,z),v'_1\big\rangle=-kz+\frac{t^{-1}}2-(t^{-1}+t)\overline u.$$

By Lemma 4.1.8,
$$\ta\big(p(t,z),S\cap\Bbb B\big)=1-\ta\big(p(t,z),g\big)=
1+\frac{(t^{-1}-t)^2}{4\big(1-|z|^2\big)}=1+h^2,\qquad
h:=\frac{|t^{-1}-t|}{2\sqrt{1-|z|^2}}.$$

In order to express explicitly the tance
$\ta\big(p(t,z),\B_1\cap\Bbb B\big)$, we need to calculate the
$\eta$-invariant (4.1.24)~:
$$\eta\big(v_1,v'_1,p(t,z)\big)=\frac{\big(2(t^{-1}+t)u+2k\overline
z-t\big)\big(2(t^{-1}+t)\overline
u+2kz-t^{-1}\big)}{2\big(1-|z|^2\big)}=a+ib,\eqno{\bold{(4.3.7)}}$$
where
$$a:=\frac{1+(t^{-1}+t)^2k^2+4k^2|z|^2+2(t^{-1}+t)k
\big((4u_0-1)z_0-4u_1z_1\big)}{2\big(1-|z|^2\big)},$$
$$b:=(t^{-1}-t)\cdot\frac{(t^{-1}+t)u_1-kz_1}{1-|z|^2},\qquad
u=u_0-iu_1,\qquad z=z_0-iz_1,\qquad u_0,u_1,z_0,z_1\in\Bbb R.$$

We have to estimate the tance $\ta\big(p(t,z),\B_1\cap\Bbb B\big)$ in
terms of $h$. To this end, we introduce
$$c:=|4u_0-1|+4|u_1|,\qquad d:=(4u_0-1)z_0-4u_1z_1,\qquad
e:=\sqrt{1-|z|^2},$$
rewrite $a$ and $b$ (using the equality $t^{-1}+t=2\sqrt{e^2h^2+1}$) as
$$a=\frac{1+8k^2+4k^2e^2h^2-4k^2e^2+4kd\sqrt{e^2h^2+1}}{2e^2},\qquad
b=\pm h\frac{4u_1\sqrt{e^2h^2+1}-2kz_1}e,\eqno{\bold{(4.3.8)}}$$
and observe that the inequalities
$$0<c,\qquad|d|\le c,\qquad0<e\le1,\qquad0\le h,\qquad0<k<2|u_1|$$
are valid. Only the last inequality requires a proof: It follows from
$0<u_0\le\frac14<\frac12$ that $(\frac12-2u_0)^2<\frac14$. Hence,
$k=\sqrt{4u_1^2+4u_0^2-2u_0}=\sqrt{4u_1^2+(\frac12-2u_0)^2-
\frac14}<2|u_1|$.

\smallskip

Suppose that $\ta\big(p(t,z),\B_1\cap\Bbb B\big)<1+\varepsilon^2$ for
some $\varepsilon>0$. By Lemma 4.1.25,
$\ta\big(p(t,z),\B_1\cap\Bbb B\big)=1-a+|a+ib|$ in view of (4.3.7).
Therefore, $b^2<\varepsilon^4+2a\varepsilon^2$, which can be written in
the form
$$\big(4hu_1\sqrt{e^2h^2+1}-2hkz_1\big)^2-\big(2\varepsilon
k\sqrt{e^2h^2+1}+\varepsilon
d\big)^2<\varepsilon^2\big(1+4k^2+\varepsilon^2e^2-4k^2e^2-d^2\big)$$
taking the expressions (4.3.8) into account. Using an inequality of the
type $\big(|A|-|B|\big)^2-\big(|C|+|D|\big)^2\le(A-B)^2-(C+D)^2$ and
the fact that $|d|\le c$ and $e^2\le1$, we conclude that
$$\big(4h|u_1|\sqrt{e^2h^2+1}-2hk|z_1|\big)^2-\big(2\varepsilon
k\sqrt{e^2h^2+1}+\varepsilon
c\big)^2<\varepsilon^2\big(1+4k^2+\varepsilon^2\big).$$
In view of $|z_1|<1\le\sqrt{e^2h^2+1}$ and $k<2|u_1|$, this implies
that
$$\big(4h|u_1|\sqrt{e^2h^2+1}-2hk\big)^2-\big(2\varepsilon
k\sqrt{e^2h^2+1}+\varepsilon
c\big)^2<\varepsilon^2\big(1+4k^2+\varepsilon^2\big).$$
This inequality can be converted into
$$\Big(\big(4h|u_1|-2\varepsilon k\big)\sqrt{e^2h^2+1}-(2hk+\varepsilon
c)\Big)\Big(\big(4h|u_1|+2\varepsilon
k\big)\sqrt{e^2h^2+1}-(2hk-\varepsilon
c)\Big)<\varepsilon^2\big(1+4k^2+\varepsilon^2\big).
\eqno{\bold{(4.3.9)}}$$

\smallskip

Given $\vartheta>0$, we have to find some $\varepsilon>0$ such that the
inequality $\ta(p,\B_1\cap\Bbb B)<1+\varepsilon^2$ implies
$h<\vartheta$. First, we require that $\varepsilon\le\vartheta$. We
assume now that $4h|u_1|-2\varepsilon k\ge0$ because, otherwise, the
inequality $h<\varepsilon$ follows from $0<k<2|u_1|$ and implies that
$h<\vartheta$. Next, we require that
$\varepsilon\le\displaystyle\frac{4|u_1|-2k}{c+2k}\vartheta$. Now, we
assume that $4h|u_1|-2\varepsilon k\ge2hk+\varepsilon c$ (otherwise,
the inequality $h<\vartheta$ follows) and conclude that
$$\big(4h|u_1|-2\varepsilon k\big)\sqrt{e^2h^2+1}-(2hk+\varepsilon
c)\ge0$$
in view of $\sqrt{e^2h^2+1}\ge1$. Requiring that $\varepsilon\le1$, we
can deduce from (4.3.9) that
$$\Big(\big(4h|u_1|-2\varepsilon k\big)\sqrt{e^2h^2+1}-(2hk+\varepsilon
c)\Big)^2<\varepsilon^2(2+4k^2),$$
which, in its turn, implies
$4h|u_1|-2\varepsilon k-2hk-\varepsilon c<\varepsilon\sqrt{2+4k^2}$ due
to $1\le\sqrt{e^2h^2+1}$ and
$4h|u_1|-2\varepsilon k\ge2hk+\varepsilon c$. We obtain
$h\big(4|u_1|-2k\big)<\varepsilon(\sqrt{2+4k^2}+2k+c)$ and, therefore,
$h<\displaystyle\frac{\sqrt{2+4k^2}+2k+c}{4|u_1|-2k}\varepsilon$.
Finally, we require that
$\varepsilon\le
\displaystyle\frac{4|u_1|-2k}{\sqrt{2+4k^2}+2k+c}\vartheta$.

The complete list of requirements concerning $\varepsilon$ is the
following:
$$\varepsilon\le\vartheta,\qquad\varepsilon\le\frac{4|u_1|-2k}{c+2k}
\vartheta,\qquad\varepsilon\le1,\qquad\varepsilon\le\frac{4|u_1|-2k}
{\sqrt{2+4k^2}+2k+c}\vartheta.$$
(Note that none of these inequalities involves $t$ or $z$.)

\smallskip

Now, we consider the case when $\B_1$ and $\B_2$ have the same complex
spine. Let $f$ be the common focus of $\B_1$ and $\B_2$. We assume
$\langle f,f\rangle=1$. Every point in $\B_2\cap\Bbb B$ has the form
$p(t,z):=zf+t^{-1}v_2-tv'_2$, $t>0$, $|z|<1$. We have
$$\big\langle p(t,z),g\big\rangle=\frac{t^{-1}-t}2,\qquad\big\langle
v_1,p(t,z)\big\rangle=\langle
v_1,t^{-1}v_2-t(g-v_2)\big\rangle=(t^{-1}+t)u-\frac t2,$$
$$\big\langle p(t,z),v'_1\big\rangle=\big\langle
p(t,z),g-v_1\big\rangle=\frac{t^{-1}}2-(t^{-1}+t)\overline
u,\qquad\big\langle p(t,z),p(t,z)\big\rangle=|z|^2-1.$$

By Lemma 4.1.8,
$$\ta\big(p(t,z),S\cap\Bbb B\big)=1-\ta\big(p(t,z),g\big)=
1+\frac{(t^{-1}-t)^2}{4(1-|z|^2)}=1+h^2,\qquad
h:=\frac{|t^{-1}-t|}{2\sqrt{1-|z|^2}}.$$

Taking into account that $\Re u=2|u|^2$, we calculate the
$\eta$-invariant (4.1.24) :
$$\eta\big(v_1,v'_1,p(t,z)\big)=\frac{\big(2(t^{-1}+t)u-t\big)
\big(2(t^{-1}+t)\overline u-t^{-1}\big)}{2\big(1-|z|^2\big)}=\frac
{1+2i(t^2-t^{-2})u_1}{2\big(1-|z|^2\big)}=a+ib,\eqno{\bold{(4.3.10)}}$$
where
$$a:=\frac1{2\big(1-|z|^2\big)},\qquad
b:=\frac{(t^2-t^{-2})u_1}{1-|z|^2},\qquad u_1:=\Im u.$$

Suppose that $\ta\big(p(t,z),\B_1\cap\Bbb B\big)<1+\varepsilon^2$ for
some $\varepsilon>0$. By Lemma 4.1.25,
$\ta\big(p(t,z),\B_1\cap\Bbb B\big)=1-a+|a+ib|$ in view of (4.3.10).
Therefore, $b^2<\varepsilon^4+2a\varepsilon^2$, which can be written in
the form
$$\frac{(t^{-2}-t^2)^2u_1^2}{1-|z|^2}<
\varepsilon^2\Big(\varepsilon^2\big(1-|z|^2\big)+1\Big).$$
Since $|z|<1$ and $t+t^{-1}\ge2$, we obtain
$16h^2u_1^2<\varepsilon^2(\varepsilon^2+1)$. Note that $u$ cannot be
real: this would contradict $\Re u=2|u|^2$ and $0<\Re u\le\frac14$.
Requiring that $\varepsilon^2\le15$ and $\varepsilon\le|u_1|\vartheta$,
we arrive at $h<\vartheta$.\hfill$_\blacksquare$

\medskip

\centerline{\bf4.4.~Triangles of bisectors}

\medskip

This subsection begins with a numerical criterion of transversality of
an oriented triangle of bisectors. Then we establish the
path-connectedness of the space of oriented transversal triangles.
Finally, we~calculate the trace of the holonomy of an oriented triangle
(see 2.5.1) and prove that counterclockwise-oriented transversal
triangles can be neither R-parabolic nor trivial.

\smallskip

{\bf4.4.1.}~Throughout this subsection, $g_1,g_2,g_3$ denote positive
points such that $\ta(g_i,g_j)>1$ for $i\ne j$. Put
$$t_{ij}:=\sqrt{\ta(g_i,g_j)}{\text{ \ \rm for }}i\ne
j,\qquad\varkappa:=\frac{\langle g_1,g_2\rangle\langle
g_2,g_3\rangle\langle g_3,g_1\rangle}{\langle g_1,g_1\rangle\langle
g_2,g_2\rangle\langle
g_3,g_3\rangle},\qquad\varepsilon:=\frac\varkappa{|\varkappa|},\qquad
G:=\left(\smallmatrix1&t_{12}&t_{31}\overline\varepsilon\\
t_{12}&1&t_{23}\\t_{31}\varepsilon&t_{23}&1\endsmallmatrix\right),$$
$$\qquad\varepsilon_0:=\Re\varepsilon,\qquad\varepsilon_1:=
\Im\varepsilon,\qquad d:=\det
G=1+2t_{12}t_{23}t_{31}\varepsilon_0-t_{12}^2-t_{23}^2-t_{31}^2,\qquad
S_i:=\Bbb Pg_i^\perp.$$
We fix suitable representatives such that $G$ is the Gram matrix of
$g_1,g_2,g_3\in V$. The numbers $t_{ij}$'s and $\varepsilon$ are
invariant under the action of $\PU$ on the ordered triples
$(g_1,g_2,g_3)$ and constitute a complete set of invariants. Indeed, if
the triples $(g_1,g_2,g_3)$ and $(g'_1,g'_2,g'_3)$ have the same Gram
matrix, then there exists some $X\in\U$ such that $Xg_i=g'_i$. By
Sylvester's Criterion [KoM, p.~113], all the values of the $t_{ij}$'s
and $\varepsilon$ subject to the conditions $t_{ij}>1$,
$|\varepsilon|=1$, and $d\le0$ are possible. The equality $d=0$ means
that $g_1,g_2,g_3$ lie on the same projective line.

Recall that the triangle of bisectors $\Delta(S_1,S_2,S_3)$ is {\it
transversal\/} if the bisectors $\B{\wr}g_i,g_j{\wr}$ and
$\B{\wr}g_i,g_k{\wr}$ are transversal along their common slice
$S_i\cap\overline{\Bbb B}$ for all $i,j,k$. The triangle
$\Delta(S_1,S_2,S_3)$ is said to be {\it counterclockwise-oriented\/}
if the oriented angle from $\B[S_2,S_3]$ to $\B[S_2,S_1]$ at some point
$p\in S_2\cap\Bbb B$ does not exceed $\pi$ (see Definition 2.3.5).

\medskip

{\bf Criterion 4.4.2.} {\sl A transversal triangle of bisectors is
counterclockwise-oriented if and only if\/ $\varepsilon_1<0$. Suppose
that\/ $1<t_{12}\le t_{23},t_{31}$. Then\/ $\Delta(S_1,S_2,S_3)$ is
transversal if and only if\/
$t_{12}^2\varepsilon_0^2+t_{23}^2+t_{31}^2<
1+2t_{12}t_{23}t_{31}\varepsilon_0$.}

\medskip

{\bf Proof.} For the first assertion, we measure the oriented angle
from $\B[S_2,S_3]$ to $\B[S_2,S_1]$ at the point
$p:=\pi[g_2]g_1=g_1-t_{12}g_2\in S_2$. (It follows from $t_{12}>1$ that
$p\in S_2\cap\Bbb B$.) By Lemma 4.3.1, this angle equals
$$\Arg\frac{\langle g_3,p\rangle\langle p,g_1\rangle}{\langle
g_3,g_2\rangle\langle
g_2,g_1\rangle}=\Arg\big((t_{31}\varepsilon-t_{12}t_{23})
(1-t_{12}^2)\big)=\Arg(t_{12}t_{23}-t_{31}\varepsilon),$$
where $\Arg$ takes values in $[0,2\pi)$. The triangle is
counterclockwise-oriented exactly when
$\Im(t_{12}t_{23}-t_{31}\varepsilon)>0$, that is, when
$\varepsilon_1<0$.

By Criterion 4.3.3, the transversality of $\B{\wr}g_1,g_2{\wr}$ and
$\B{\wr}g_2,g_3{\wr}$ along $S_2\cap\overline{\Bbb B}$ is equivalent to
the inequality
$$\bigg|\Re\frac{t_{31}\overline\varepsilon}{t_{12}t_{23}}-1\bigg|<
\sqrt{1-\frac1{t_{12}^2}}\cdot\sqrt{1-\frac1{t_{23}^2}}.$$
This inequality is equivalent to
$(t_{31}\varepsilon_0-t_{12}t_{23})^2<\sqrt{(t_{12}^2-1)(t_{23}^2-1)}$,
that is, to
$$t_{12}^2+t_{23}^2+t_{31}^2\varepsilon_0^2<
1+2t_{12}t_{23}t_{31}\varepsilon_0.$$
By symmetry, the other two transversalities are equivalent to the
inequalities
$$t_{12}^2+t_{23}^2\varepsilon_0^2+t_{31}^2<1+2t_{12}t_{23}t_{31}
\varepsilon_0,\qquad t_{12}^2\varepsilon_0^2+t_{23}^2+t_{31}^2<
1+2t_{12}t_{23}t_{31}\varepsilon_0.$$
The last inequality implies the other two because
$t_{12}\le t_{23},t_{31}$.\hfill$_\blacksquare$

\medskip

The path-connectedness of the space of oriented transversal triangles
of bisectors is the subject of the following

\medskip

{\bf Lemma 4.4.3.} {\sl The region given in\/ $\Bbb R^4(e,t_1,t_2,t_3)$
by the inequalities\/ $1<t_1\le t_2,t_3$ and\/
$t_1^2e^2+t_2^2+t_3^2<1+2t_1t_2t_3e\le t_1^2+t_2^2+t_3^2$ is
path-connected.}

\medskip

{\bf Proof.} The inequalities in the lemma imply $e^2<1$ and, hence,
$t_1e<t_2t_3$. The inequality $t_1^2e^2+t_2^2+t_3^2<1+2t_1t_2t_3e$ can
be rewritten as $(t_2t_3-t_1e)^2<(t_2^2-1)(t_3^2-1)$. Therefore, it is
equivalent to the inequalities
$t_2t_3-\sqrt{(t_2^2-1)(t_3^2-1)}<t_1e\le t_2t_3$. Now, from the
inequality $1\le t_2t_3-\sqrt{(t_2^2-1)(t_3^2-1)}$ implied by
$1<t_2t_3$, we conclude that $1<t_1e$. In particular, $0<e<1$.

The inequalities
$$1<t_1\le t_2,t_3,\qquad t_2t_3-\sqrt{(t_2^2-1)(t_3^2-1)}<t_1e\le
t_2t_3,\qquad1+2t_1t_2t_3e\le t_1^2+t_2^2+t_3^2\eqno{\bold{(4.4.4)}}$$
are equivalent to those in the lemma. We have seen that they imply
$0<e<1<t_1e$.

\smallskip

Suppose that $e,t_1,t_2,t_3$ satisfy the inequalities in the lemma.
Keeping the inequalities (4.4.4), we~will vary the $t_i$'s until we
reach $t_1=t_2=t_3$ or $1+2t_1t_2t_3e=t_1^2+t_2^2+t_3^2$. Without loss
of generality, we~can assume that $t_2\le t_3$. It follows from the
inequality $t_2<t_2t_3e$ (implied by $1<t_1e$ and $t_1\le t_3$) that
the function $f(x):=x^2+t_2^2+t_3^2-2xt_2t_3e-1$ is decreasing in
$x\in[t_1,t_2]$. Increasing $t_1$ and keeping the inequalities (4.4.4)
and $t_2\le t_3$, we can reach a point such that either $t_1=t_2$ or
$f(t_1)=0$. The latter means that $1+2t_1t_2t_3e=t_1^2+t_2^2+t_3^2$.

In the case of $t_1=t_2$, we can rewrite the condition $t_1=t_2$ and
the inequalities in the lemma as
$$1<t_1=t_2\le t_3,\qquad2t_1^2(t_3e-1)\le
t_3^2-1<t_1^2(2t_3e-e^2-1).\eqno{\bold{(4.4.5)}}$$
As we saw, these conditions imply that $0<e<1<t_1e$. Hence, $0<t_3e-1$
and $0<2t_3e-e^2-1$. Now, we increase $t_1=t_2$ keeping conditions
(4.4.5) and reach a point where $t_1=t_2=t_3$ or
$2t_1^2(t_3e-1)=t_3^2-1$. The latter means that
$1+2t_1t_2t_3e=t_1^2+t_2^2+t_3^2$.

In the case of $t_1=t_2=t_3$, the inequalities (4.4.4) take the form
$$1<t_1=t_2=t_3,\qquad1<t_1e\le t_1^2,\qquad
t_1e\le\frac{3t_1^2-1}{2t_1^2}.\eqno{\bold{(4.4.6)}}$$
It follows from $\displaystyle\frac{3t_1^2-1}{2t_1^2}<t_1^2$ (implied
by $1<t_1$) that we can increase $e$ and keep conditions (4.4.6) until
$t_1e=\displaystyle\frac{3t_1^2-1}{2t_1^2}$. Again, we reach a point
satisfying the equality $1+2t_1t_2t_3e=t_1^2+t_2^2+t_3^2$.

It remains to show that the region given by the inequalities in the
lemma and by the equality $1+2t_1t_2t_3e=t_1^2+t_2^2+t_3^2$ is
path-connected. Excluding $e$ with the help of the equality, we can
describe this region as the one given in $\Bbb R^3(t_1,t_2,t_3)$ by the
inequalities
$$1<t_1\le t_2,t_3,\qquad
t_1^2+t_2^2+t_3^2<1+2t_1t_2t_3.\eqno{\bold{(4.4.7)}}$$
Assuming that $t_2\le t_3$, we obtain $t_2<t_1t_3$. Therefore, keeping
the inequalities (4.4.7), we can increase $t_2$ until $t_2=t_3$. Now,
the inequalities (4.4.7) are equivalent to the conditions
$1<t_1\le t_2=t_3$.\hfill$_\blacksquare$

\medskip

{\bf4.4.8.~Trace of holonomy.} Let us recall the definition of the
holonomy $\varphi$ of $\Delta(S_1,S_2,S_3)$. Put
$$m_1:=\frac{g_1+g_2}{\sqrt{2t_{12}+2}},\qquad
m_2:=\frac{g_2+g_3}{\sqrt{2t_{23}+2}},\qquad m_3:=\frac{\varepsilon
g_1+g_3}{\sqrt{2t_{31}+2}}.\eqno{\bold{(4.4.9)}}$$
By Lemma 4.1.23, $R_i:=R(m_i)$ is the reflection (see (2.1.3) for the
definition) in the middle slice of the segment $\B[S_i,S_{i+1}]$ and
$\langle m_i,m_i\rangle=1$ (the indices are modulo $3$). Define
$\varphi:=R_3R_2R_1\in\SU$. It is easy to see that $R_1g_1=g_2$,
$R_2g_2=g_3$, and $R_3g_3=\varepsilon g_1$. Hence,
$$\varphi g_1=\varepsilon g_1.\eqno{\bold{(4.4.10)}}$$
In particular, $g_1^\perp$ is stable under $\varphi$. Therefore, the
isometry $\varphi$ induces an isometry $\psi\in\SU(g_1^\perp)$ of the
complex geodesic $\Bbb Pg_1^\perp$.

\medskip

{\bf Lemma 4.4.11.}
{\sl$|\tr\psi|=\sqrt{2(1+\varepsilon_0)}\Big(1-\displaystyle\frac
d{(t_{12}+1)(t_{23}+1)(t_{31}+1)}\Big)$.}

\medskip

{\bf Proof}. Since $\varphi\in\SU$ and $\varphi g_1=\varepsilon g_1$,
it suffices to show that
$$\tr\varphi=\varepsilon-(1+\overline\varepsilon)\Big(1-\frac
d{(t_{12}+1)(t_{23}+1)(t_{31}+1)}\Big).\eqno{\bold{(4.4.12)}}$$
Indeed, (4.4.12) implies that
$\tr\psi=\pm(\sqrt{\varepsilon}+\overline{\sqrt{\varepsilon}})
\Big(1-{\displaystyle\frac d{(t_{12}+1)(t_{23}+1)(t_{31}+1)}}\Big)$
because $\det\varphi=1$ and $\varphi g_1=\varepsilon g_1$.

\smallskip

Define the linear maps $\varphi_i\in\Lin_\Bbb C(V,V)$ by the rule
$\varphi_ix:=2\langle x,m_i\rangle m_i$ and denote
$g_{ij}:=\langle m_i,m_j\rangle$, $i,j=1,2,3$. Obviously,
$R_i=\varphi_i-1$ and
$$\varphi_j\varphi_ix=4\langle x,m_i\rangle
g_{ij}m_j,\qquad\varphi_k\varphi_j\varphi_ix=8\langle x,m_i\rangle
g_{ij}g_{jk}m_k.$$
Considering the orthogonal decompositions
$V=\Bbb Cm_i\oplus m_i^\perp$, we obtain
$$\tr\varphi_i=2,\qquad\tr(\varphi_j\varphi_i)=4g_{ij}g_{ji},
\qquad\tr(\varphi_k\varphi_j\varphi_i)=8g_{ij}g_{jk}g_{ki}.$$
It follows from
$\varphi=R_3R_2R_1=\varphi_3\varphi_2\varphi_1-\varphi_3\varphi_2-
\varphi_3\varphi_1-\varphi_2\varphi_1+\varphi_3+\varphi_2+\varphi_1-1$
that
$$\tr\varphi=8g_{12}g_{23}g_{31}-4g_{23}g_{32}-4g_{13}g_{31}-4g_{12}
g_{21}+3\eqno{\bold{(4.4.13)}}$$
(compare with [Pra]). From (4.4.9), we obtain
$$g_{12}=\frac{1+t_{12}+t_{23}+t_{31}\overline\varepsilon}
{2\sqrt{(t_{12}+1)(t_{23}+1)}},\qquad
g_{23}=\frac{1+t_{23}+t_{31}+t_{12}\overline\varepsilon}
{2\sqrt{(t_{23}+1)(t_{31}+1)}},\qquad
g_{31}=\frac{\varepsilon+t_{31}\varepsilon+t_{12}\varepsilon+t_{23}}
{2\sqrt{(t_{31}+1)(t_{12}+1)}},$$
$$4g_{23}g_{32}=4\ta(m_2,m_3)=\frac{(1+t_{23}+t_{31})^2+t_{12}^2+
2t_{12}(1+t_{23}+t_{31})\varepsilon_0}{(t_{23}+1)(t_{31}+1)}=$$
$$=2+\frac{2t_{12}(1+t_{23}+t_{31}+t_{23}t_{31})\varepsilon_0-d}
{(t_{23}+1)(t_{31}+1)}=2+2t_{12}\varepsilon_0-
\frac{d(t_{12}+1)}{(t_{12}+1)(t_{23}+1)(t_{31}+1)}.$$
By symmetry,
$$4g_{13}g_{31}=2+2t_{23}\varepsilon_0-\frac{d(t_{23}+1)}{(t_{12}+1)
(t_{23}+1)(t_{31}+1)},\qquad4g_{12}g_{21}=2+2t_{31}\varepsilon_0-
\frac{d(t_{31}+1)}{(t_{12}+1)(t_{23}+1)(t_{31}+1)}.$$
It is possible to show with a straightforward calculation that
$$8g_{12}g_{23}g_{31}(t_{12}+1)(t_{23}+1)(t_{31}+1)=$$
$$=(t_{12}+1)(t_{23}+1)(t_{31}+1)\big(2(t_{12}+t_{23}+t_{31})
\varepsilon_0+2+\varepsilon-\overline\varepsilon\big)-
d(t_{12}+t_{23}+t_{31}+2-\overline\varepsilon).$$
In other words,
$$8g_{12}g_{23}g_{31}=2(t_{12}+t_{23}+t_{31})\varepsilon_0+2+
\varepsilon-\overline\varepsilon-\frac{d(t_{12}+t_{23}+t_{31}+2-
\overline\varepsilon)}{(t_{12}+1)(t_{23}+1)(t_{31}+1)}.$$
Now, from (4.4.13), we infer (4.4.12).\hfill$_\blacksquare$

\medskip

{\bf4.4.14.~Parabolic holonomy.} Recall that a triangle
$\Delta(S_1,S_2,S_3)$ is {\it parabolic\/} or {\it trivial\/} if the
isometry $\psi$ (see 4.4.8) is parabolic or the identity. An oriented
parabolic triangle is L-{\it parabolic\/} if $\psi$ moves the points in
$S_1\cap\partial_\infty\Bbb B$ (different from the fixed one) in the
clockwise direction.

\medskip

{\bf Lemma 4.4.15.} {\sl Let\/ $\Delta(S_1,S_2,S_3)$ be a
counterclockwise-oriented transversal triangle. Suppose that\/
$\Delta(S_1,S_2,S_3)$ is parabolic or trivial. Then\/
$\Delta(S_1,S_2,S_3)$ is\/ {\rm L}-parabolic.}

\medskip

{\bf Proof.} Denote $p:=\pi[g_1]g_2\in S_1$. It follows from $t_{12}>1$
that $p\in S_1\cap\Bbb B$. Let $q\in S_1\cap\partial_\infty\Bbb B$
denote a fixed point of $\varphi$. We assume that
$\langle p,q\rangle=1$ and define
$p':=\langle p,p\rangle q-p\in p^\perp$. Since
$\langle p',p'\rangle=-\langle p,p\rangle$, every point in
$S_1\cap\overline{\Bbb B}$ has the form
$$p(z):=zp'+p,\qquad z\in\Bbb C,\qquad|z|\le1,$$
in the orthogonal basis $p,p'\in g_1^\perp$. Clearly, $p(1)\sim q$,
where $\sim$ stands for $\Bbb C^*$-proportionality.

Since $\varphi$ is either parabolic or trivial on $S_1=\L{\wr}p,q{\wr}$
with fixed point $q$, we have
$$\varphi p=up+uivq,\qquad\varphi q=uq\eqno{\bold{(4.4.16)}}$$
for suitable $u,v\in\Bbb C$, $u\ne0$. This implies that
$$\tr\varphi=2u+\varepsilon$$
because $\varphi g_1=\varepsilon g_1$ by (4.4.10). From
$\varphi\in\SU$, $\langle p,q\rangle=1$, and $\langle q,q\rangle=0$, we
conclude that
$$1=\langle p,q\rangle=\langle\varphi p,\varphi
q\rangle=|u|^2,\qquad\langle p,p\rangle=\langle\varphi p,\varphi
p\rangle=|u|^2\big(\langle p,p\rangle-i\overline v+iv\big)=\langle
p,p\rangle-2\Im v.$$
Hence, $v\in\Bbb R$.

Let us understand how does $\varphi$ act on $S_1$ in terms of $z$ :
$$\varphi\big(p(z)\big)=z\varphi\big(\langle p,p\rangle
q-p\big)+\varphi p\sim z\langle p,p\rangle q+(1-z)(p+ivq)=$$
$$=\frac{z\langle p,p\rangle+(1-z)iv}{\langle p,p\rangle}\big(\langle
p,p\rangle q-p\big)+\frac{\langle p,p\rangle+(1-z)iv}{\langle
p,p\rangle}p\sim p(z'),$$
where
$$z':=\frac{z\langle p,p\rangle+(1-z)iv}{\langle p,p\rangle+(1-z)iv}.$$
In particular, $\varphi\big(p(-1)\big)\sim p(z_0)$ with
$z_0:=\displaystyle\frac{2iv-\langle p,p\rangle}{2iv+\langle
p,p\rangle}$
and
$\Im z_0=\displaystyle\frac{4v\langle p,p\rangle}{\big|2iv+\langle
p,p\rangle\big|^2}$.
We see that $\varphi$ moves $p(-1)$ in the clockwise direction exactly
when $\Im z_0>0$, i.e., if and only if $v<0$. In other words, the
triangle is L-parabolic if and only if $v<0$. It is easy to see from
(4.4.16) that the triangle is trivial if and only if $v=0$.

\smallskip

Thus, we have to show that $v<0$.

\smallskip

It follows from $\langle p,q\rangle=1$ and (4.4.16) that
$\langle\varphi p,p\rangle=u\big(\langle p,p\rangle+iv\big)$, i.e.,
$v=\Im\big(u^{-1}\langle\varphi p,p\rangle\big)$. By~(4.4.12),
$$\tr\varphi=2u+\varepsilon=\varepsilon-(1+\overline\varepsilon)
\Big(1-\frac d{(t_{12}+1)(t_{23}+1)(t_{31}+1)}\Big)$$
with $d\le0$. Consequently, $v$ and
$\Im\big(-(1+\varepsilon)\langle\varphi p,p\rangle\big)$ have the same
signs.

Let us calculate $\langle\varphi p,p\rangle$ in terms of the Gram
matrix $G$. We only need to find $\varphi g_2$ since $p=g_2-t_{12}g_1$
and $\varphi g_1=\varepsilon g_1$ by (4.4.10). Using (4.4.9) and
(2.1.3), we see that
$$R_1x=\langle x,g_1+g_2\rangle\frac{g_1+g_2}{t_{12}+1}-x,\qquad
R_2x=\langle x,g_2+g_3\rangle\frac{g_2+g_3}{t_{23}+1}-x,\qquad
R_3x=\langle x,\varepsilon g_1+g_3\rangle\frac{\varepsilon
g_1+g_3}{t_{31}+1}-x$$
for all $x\in V$. So, $R_1g_2=g_1$,
$R_2g_1=\displaystyle\frac{t_{12}+t_{31}\overline\varepsilon}{t_{23}+1}
(g_2+g_3)-g_1$,
and
$$\varphi g_2=\Big\langle\frac{t_{12}+t_{31}\overline\varepsilon}
{t_{23}+1}(g_2+g_3)-g_1,\varepsilon g_1+g_3\Big\rangle\frac{\varepsilon
g_1+g_3}{t_{31}+1}-\frac{t_{12}+t_{31}\overline\varepsilon}{t_{23}+1}
(g_2+g_3)+g_1=$$
$$=\Big(\frac{t_{12}+t_{31}\overline\varepsilon}{t_{23}+1}
(t_{12}\overline\varepsilon+t_{31}+t_{23}+1)-\overline\varepsilon-
t_{31}\overline\varepsilon\Big)\frac{\varepsilon g_1+g_3}{t_{31}+1}-
\frac{t_{12}+t_{31}\overline\varepsilon}{t_{23}+1}(g_2+g_3)+g_1=$$
$$=\frac{(t_{12}+t_{31}\overline\varepsilon)
(t_{12}\overline\varepsilon+t_{31})+(t_{23}+1)(t_{12}-
\overline\varepsilon)}{(t_{23}+1)(t_{31}+1)}(\varepsilon
g_1+g_3)-\frac{t_{12}+t_{31}\overline\varepsilon}{t_{23}+1}(g_2+g_3)+
g_1.$$
Therefore, taking into account that
$\langle g_1,g_2-t_{12}g_1\rangle=0$, we obtain
$$(t_{23}+1)(t_{31}+1)\langle\varphi p,p\rangle=(t_{23}+1)(t_{31}+1)
\langle\varphi g_2,g_2-t_{12}g_1\rangle=$$
$$=(t_{23}+1)(t_{31}+1)\Big\langle\frac{(t_{12}+t_{31}\overline
\varepsilon)(t_{12}\overline\varepsilon+t_{31})+(t_{23}+1)(t_{12}-
\overline\varepsilon)}{(t_{23}+1)(t_{31}+1)}g_3-\frac{t_{12}+t_{31}
\overline\varepsilon}{t_{23}+1}(g_2+g_3),g_2-t_{12}g_1\Big\rangle=$$
$$=(t_{12}+t_{31}\overline\varepsilon)(t_{12}\overline\varepsilon+
t_{31})(t_{23}-t_{12}t_{31}\varepsilon)+(t_{23}+1)
(t_{12}-\overline\varepsilon)(t_{23}-t_{12}t_{31}\varepsilon)-$$
$$-(t_{12}+t_{31}\overline\varepsilon)(t_{31}+1)
(1-t_{12}^2+t_{23}-t_{12}t_{31}\varepsilon)=$$
$$=t_{12}t_{23}t_{31}\overline\varepsilon^2+(t_{12}^2t_{23}+
t_{12}^2t_{31}-t_{23}^2-t_{23}t_{31}-t_{31}^2-t_{23}-
t_{31})\overline\varepsilon+t_{12}^3+t_{12}t_{23}^2+t_{12}t_{23}t_{31}+
t_{12}t_{31}^2-t_{12}-t_{12}^2t_{23}t_{31}\varepsilon.$$
Hence,
$$\Im\big(-(t_{23}+1)(t_{31}+1)(1+\varepsilon)\langle\varphi
p,p\rangle\big)=t_{12}^2t_{23}t_{31}\Im\varepsilon^2-
t_{12}t_{23}t_{31}\Im\overline\varepsilon^2-$$
$$-(t_{12}^2t_{23}+t_{12}^2t_{31}+t_{12}t_{23}t_{31}-t_{23}^2-
t_{23}t_{31}-t_{31}^2-t_{23}-t_{31})\Im\overline\varepsilon-$$
$$-(t_{12}^3+t_{12}t_{23}^2+t_{12}t_{23}t_{31}+t_{12}t_{31}^2-t_{12}-
t_{12}^2t_{23}t_{31})\Im\varepsilon=$$
$$=(t_{12}+1)\varepsilon_1\big(1+2t_{12}t_{23}t_{31}\varepsilon_0-
t_{12}^2-t_{23}^2-t_{31}^2+(t_{12}-1)(t_{23}+1)(t_{31}+1)\big)=$$
$$=(t_{12}+1)\varepsilon_1\big(d+(t_{12}-1)(t_{23}+1)(t_{31}+1)\big).$$

By Lemma 4.4.11,
$$\sqrt{2(1+\varepsilon_0)}\Big(1-\frac
d{(t_{12}+1)(t_{23}+1)(t_{31}+1)}\Big)=2$$
because the triangle is parabolic or trivial. So,
$$d=\Big(1-\frac2{\sqrt{2(1+\varepsilon_0)}}\Big)(t_{12}+1)(t_{23}+1)
(t_{31}+1).$$
We see that the sign of $v$ is that of
$\varepsilon_1\Big(t_{12}-\displaystyle\frac{t_{12}+1}
{\sqrt{2(1+\varepsilon_0)}}\Big)$.
By Criterion 4.4.2, $\varepsilon_1<0$ since the triangle is
counterclockwise-oriented. It remains to prove that
$0<t_{12}-\displaystyle\frac{t_{12}+1}{\sqrt{2(1+\varepsilon_0)}}$,
i.e., that $2<(t_{12}-1)^2+2t_{12}^2\varepsilon_0$. It suffices to show
that $1<t_{12}\varepsilon_0$.

In order see that $1<t_{12}\varepsilon_0$, note that the inequality
$1\le t_{23}t_{31}-\sqrt{(t_{23}^2-1)(t_{31}^2-1)}$ follows from
$1<t_{23}t_{31}$. By Criterion 4.4.2, we have
$t_{12}^2\varepsilon_0^2+t_{23}^2+t_{31}^2<
1+2t_{12}t_{23}t_{31}\varepsilon_0$
or, equivalently,
$(t_{23}t_{31}-t_{12}\varepsilon_0)^2<(t_{23}^2-1)(t_{31}^2-1)$. The
inequality $t_{12}\varepsilon_0\le t_{23}t_{31}$ implies
$t_{23}t_{31}-\sqrt{(t_{21}^2-1)(t_{31}^2-1)}<
t_{12}\varepsilon_0$.~\hfill$_\blacksquare$

\medskip

{\bf4.4.17.~$\Bbb C$-plane triangles.} Assume that $d=0$. This means
that the triangle of bisectors $\Delta(S_1,S_2,S_3)$ is built over a
usual geodesic triangle $\Delta(p_1,p_2,p_3)$ situated in a complex
geodesic $C$ of signature $+-$ which is actually the common complex
spine of the bisectors. In other words, $\B[S_i,S_j]=\B[p_i,p_j]$,
where $p_i$ is orthogonal to $g_i$ in $C$ (see Remark 4.1.10).

The triangle $\Delta(S_1,S_2,S_3)$ is always transversal. This follows
from [Mos, p.~186, Lemma 2.3.4] (Mostow shows that the angle between
cospinal bisectors is constant) or from Criterion 4.4.2 :
The~inequality
$t_{12}^2\varepsilon_0^2+t_{23}^2+t_{31}^2<
1+2t_{12}t_{23}t_{31}\varepsilon_0$
follows from
$1+2t_{12}t_{23}t_{31}\varepsilon_0=t_{12}^2+t_{23}^2+t_{31}^2$ and
$\varepsilon_0^2<1$. (Note that $\varepsilon_0^2=1$ would imply that
the points $g_1,g_2,g_3$ belong to the same geodesic by Lemma 4.1.13,
implying that $p_1,p_2,p_3$ belong to the same geodesic by Remark
4.1.10.) Obviously, $\Delta(S_1,S_2,S_3)$ is counterclockwise-oriented
exactly when $\Delta(p_1,p_2,p_3)$ is counterclockwise-oriented in $C$.

\medskip

{\bf Lemma 4.4.18.} {\sl Suppose that\/ $\Delta(S_1,S_2,S_3)$ is
counterclockwise-oriented and\/ $\Bbb C$-plane. Then\/ $\varphi$
{\rm(}see\/~{\rm4.4.8} for the definition\/{\rm)} restricted to\/ $S_1$
is a rotation about\/ $p_1$ by the angle\/
$-2\Area\Delta(p_1,p_2,p_3)$.}

\medskip

{\bf Proof.} Let $f\notin\overline{\Bbb B}$ denote the polar point to
the complex geodesic $C$ containing $\Delta(p_1,p_2,p_3)$. Obviously,
$\varphi f=-f$. By (4.4.10), $\varphi g_1=\varepsilon g_1$. Taking into
account that the points $f,p_1,g_1$ form an orthogonal basis and that
$\varphi\in\SU$, we obtain $\varphi p_1=-\overline\varepsilon p_1$.
Assuming that $\langle f,f\rangle=1$ and $\langle p_1,p_1\rangle=-1$,
we write every point in $S_1\cap\overline{\Bbb B}$ in the form
$zf+p_1$, $|z|\le1$. In terms of $z$, the isometry $\varphi$ acts on
$S_1\cap\overline{\Bbb B}$ as the multiplication by $\varepsilon$. So,
the restriction of $\varphi$ to $S_1\cap\overline{\Bbb B}$ is the
rotation about $p_1$ by the angle $\arg\varepsilon\in[-\pi,\pi]$. It
remains to show that $2\Area\Delta(p_1,p_2,p_3)=-\arg\varepsilon$.

We apply to the representatives $\pi[g_1]g_2,\pi[g_2]g_1,\pi[g_3]g_2$
of $p_1,p_2,p_3$ the known formula
$$\Area\Delta(p_1,p_2,p_3)=\textstyle\frac12\arg\big(-\langle
p_1,p_2\rangle\langle p_2,p_3\rangle\langle p_3,p_1\rangle\big)$$
for the area of the oriented triangle $\Delta(p_1,p_2,p_3)$ (see [AGr1]
or [Gol, p.~25]; by Corollary 4.1.18, the~metric we use differs by the
factor of $4$ from the one in [Gol]) and obtain
$$2\Area\Delta(p_1,p_2,p_3)=\arg\big(-\langle
g_2-t_{12}g_1,g_1-t_{12}g_2\rangle\langle
g_1-t_{12}g_2,g_2-t_{23}g_3\rangle\langle
g_2-t_{23}g_3,g_2-t_{12}g_1\rangle\big)=$$
$$=\arg\big(-t_{12}(t_{12}^2-1)t_{23}(t_{12}t_{23}-t_{31}\overline
\varepsilon)(1-t_{12}^2-t_{23}^2+t_{12}t_{23}t_{31}\varepsilon)\big)=$$
$$=\arg\big((t_{12}t_{23}-t_{31}\overline\varepsilon)(t_{12}^2+t_{23}^2
-t_{12}t_{23}t_{31}\varepsilon-1)\big).$$
It follows from
$1+2t_{12}t_{23}t_{31}\varepsilon_0=t_{12}^2+t_{23}^2+t_{31}^2$ that
$$t_{12}^2+t_{23}^2-t_{12}t_{23}t_{31}\varepsilon-1=
2t_{12}t_{23}t_{31}\varepsilon_0-t_{12}t_{23}t_{31}\varepsilon-t_{31}^2
=t_{31}(t_{12}t_{23}\overline\varepsilon-t_{31})=
t_{31}(t_{12}t_{23}-t_{31}\varepsilon)\overline\varepsilon,$$
which implies the result.\hfill$_\blacksquare$

\medskip

\centerline{\bf4.5.~K\"ahler primitive}

\medskip

The hermitian form (4.1.2) defines a {\it K\"ahler form\/} $\omega$ on
$\Bbb P\setminus\partial_\infty\Bbb B$ by the rule
$\omega(v_p,w_p):=\Im\langle v_p,w_p\rangle$, where
$v_p,w_p\in\T_p\Bbb P$ (see (4.1.3) for the definition of $v_p$). In
this subsection, we obtain a primitive $P_c$ for $\omega$ that depends
on the choice of a {\it base\/} point $c$ and study how this primitive
changes when we alter the base point.

\smallskip

{\bf4.5.1.}~Let $c\in\Bbb P$, let
$p\in\Bbb P\setminus(\Bbb Pc^\perp\cup\partial_\infty\Bbb B)$, and let
$v_p\in\T_p\Bbb P$. Define
$$P_c(v_p):=-\Im\Big(\frac{\langle p,p\rangle\langle
v,c\rangle}{2\langle p,c\rangle}\Big).$$
Note that $P_c(v_p)$ does not depend on the choice of representatives
$c,p,v\in V$ that give the same $v_p$. Obviously, $P_c(v_p)$ depends
smoothly on $v_p\in\T_p\Bbb P$, i.e., it defines a $1$-form on
$\Bbb P\setminus(\Bbb Pc^\perp\cup\partial_\infty\Bbb B)$.

Let $c_1,c_2\in\Bbb P$ be nonorthogonal. Define the multi-valued
function
$$f_{c_1,c_2}(p):=\frac12\Arg\Big(\frac{\langle c_1,p\rangle\langle
p,c_2\rangle}{\langle c_1,c_2\rangle}\Big),$$
where $p\in\Bbb P\setminus(\Bbb Pc_1^\perp\cup\Bbb Pc_2^\perp)$. Note
that $f_{c_1,c_2}(p)$ does not depend on the choice of representatives
$c_1,c_2,p\in V$.

\medskip

{\bf Proposition 4.5.2.} {\sl For nonorthogonal\/ $c_1,c_2\in\Bbb P$,
we have $P_{c_1}-P_{c_2}=df_{c_1,c_2}$ on\/
$\Bbb P\setminus(\partial_\infty\Bbb B\cup\Bbb Pc_1^\perp\cup\Bbb
Pc_2^\perp)$.}

\medskip

{\bf Proof.} Choose representatives $c_1,c_2\in V$ such that
$\langle c_1,c_2\rangle\in\Bbb R$. Take
$p\in\Bbb P\setminus(\partial_\infty\Bbb B\cup\Bbb Pc_1^\perp\cup\Bbb
Pc_2^\perp)$
and $v_p\in\T_p\Bbb P$. Then
$$df_{c_1,c_2}(v_p)=\frac
d{d\varepsilon}\Big|_{\varepsilon=0}f_{c_1,c_2}
\big(p+\varepsilon\langle p,p\rangle v\big)=\frac12\frac
d{d\varepsilon}\Big|_{\varepsilon=0}\Arg\Big(\big\langle
c_1,p+\varepsilon\langle p,p\rangle v\big\rangle\big\langle
p+\varepsilon\langle p,p\rangle v,c_2\big\rangle\Big)=$$
$$=\frac12\Im\frac
d{d\varepsilon}\Big|_{\varepsilon=0}\Big(\ln\big\langle
c_1,p+\varepsilon\langle p,p\rangle v\big\rangle+\ln\big\langle
p+\varepsilon\langle p,p\rangle v,c_2\big\rangle\Big)=$$
$$=\frac12\Im\Big(\frac{\langle p,p\rangle\langle c_1,v\rangle}{\langle
c_1,p\rangle}+\frac{\langle p,p\rangle\langle v,c_2\rangle}{\langle
p,c_2\rangle}\Big)=
\big(P_{c_1}-P_{c_2}\big)(v_p).\eqno{_\blacksquare}$$

{\bf4.5.3.}~In order to show that $P_c$ is a primitive for the form
$\omega$, we introduce an auxiliary vector field. Let
$p\in\Bbb P\setminus\partial_\infty\Bbb B$ and let $v_p\in\T_p\Bbb P$.
The smooth vector field
$$\Tn(v_p)(x):=\langle-,x\rangle\frac{\langle
p,p\rangle\pi[x]v}{\langle p,x\rangle}$$
is defined on
$\Bbb P\setminus(\Bbb Pp^\perp\cup\partial_\infty\Bbb B)$. This field
extends the vector $v_p$, i.e., $\Tn(v_p)(p)=v_p$. Note that
$\Tn(v_p)(x)$ does not depend on the choice of representatives
$p,v,x\in V$ that give the same $v_p$.

\medskip

{\bf Lemma 4.5.4.} {\sl Let\/
$p\in\Bbb P\setminus\partial_\infty\Bbb B$ and let\/
$v_p,w_p\in\T_p\Bbb P$. Then\/ $\big[\Tn(v_p),\Tn(w_p)\big](p)=0$.}

\medskip

{\bf Proof.} Denote by $\hat f$ the lift to $V$ of a $C^2$-function $f$
defined in a neighbourhood of $p$. So, $\hat f(ux)=\hat f(x)$ for
$0\ne u\in\Bbb C$. By definition,
$$\Tn(w_p)f(x)=\frac d{d\varepsilon}\Big|_{\varepsilon=0}\hat
f\Big(x+\varepsilon\langle x,x\rangle\frac{\langle
p,p\rangle\pi[x]w}{\langle p,x\rangle}\Big).$$
Therefore,
$$\Tn(v_p)\big(\Tn(w_p)f\big)(p)=$$
$$=\frac d{d\delta}\Big|_{\delta=0}\bigg(\frac
d{d\varepsilon}\Big|_{\varepsilon=0}\hat f\Big(p+\delta\langle
p,p\rangle v+\varepsilon\big\langle p+\delta\langle p,p\rangle
v,p+\delta\langle p,p\rangle v\big\rangle\frac{\langle
p,p\rangle\pi\big[p+\delta\langle p,p\rangle v\big]w}{\big\langle
p,p+\delta\langle p,p\rangle v\big\rangle}\Big)\bigg)=$$
$$=\frac d{d\delta}\Big|_{\delta=0}\bigg(\frac
d{d\varepsilon}\Big|_{\varepsilon=0}\hat f\Big(p+\delta\langle
p,p\rangle v+\varepsilon\big(\langle p,p\rangle+\delta^2\langle
p,p\rangle^2\langle v,v\rangle\big)w-\varepsilon\delta\langle
p,p\rangle\langle w,v\rangle\big(p+\delta\langle p,p\rangle
v\big)\Big)\bigg)=$$
$$=\frac d{d\delta}\Big|_{\delta=0}\bigg(\frac
d{d\varepsilon}\Big|_{\varepsilon=0}\hat
f\Big(\big(1-\varepsilon\delta\langle p,p\rangle\langle
w,v\rangle\big)\big(p+\delta\langle p,p\rangle
v\big)+\varepsilon\big(\langle p,p\rangle+\delta^2\langle
p,p\rangle^2\langle v,v\rangle\big)w\Big)\bigg)=$$
$$\frac d{d\delta}\Big|_{\delta=0}\bigg(\frac
d{d\varepsilon}\Big|_{\varepsilon=0}\hat f\Big(p+\delta\langle
p,p\rangle v+\varepsilon\langle p,p\rangle\frac{1+\delta^2\langle
p,p\rangle\langle v,v\rangle}{1-\varepsilon\delta\langle
p,p\rangle\langle w,v\rangle}w\Big)\bigg)=\frac
d{d\delta}\Big|_{\delta=0}\Big(\frac
d{d\varepsilon}\Big|_{\varepsilon=0}\hat f\big(p+\delta\langle
p,p\rangle v+\varepsilon\langle p,p\rangle w\big)\Big).$$
By symmetry,
$$\Tn(w_p)\big(\Tn(v_p)f\big)(p)=\frac
d{d\delta}\Big|_{\delta=0}\Big(\frac
d{d\varepsilon}\Big|_{\varepsilon=0}\hat f(p+\delta\langle p,p\rangle
w+\varepsilon\langle p,p\rangle v)\Big).\eqno{_\blacksquare}$$

{\bf Lemma 4.5.5.} {\sl Let\/ $c\in\Bbb P$, let\/
$p\in\Bbb P\setminus(\Bbb Pc^\perp\cup\partial_\infty\Bbb B)$, and
let\/ $v_p,w_p\in\T_p\Bbb P$. Then
$$v_pP_c\big(\Tn(w_p)\big)=\frac{\langle p,p\rangle}2\Im\Big(\langle
w,v\rangle+\frac{\langle p,p\rangle\langle v,c\rangle\langle
w,c\rangle}{\langle p,c\rangle^2}\Big).$$}

{\bf Proof.} By definition,
$$P_c\big(\Tn(w_p)\big)(x)=-\Im\frac{\langle x,x\rangle\langle
p,p\rangle\big\langle\pi[x]w,c\big\rangle}{2\langle x,c\rangle\langle
p,x\rangle}=$$
$$=-\Im\frac{\langle x,x\rangle\langle p,p\rangle\langle
w,c\rangle-\langle p,p\rangle\langle w,x\rangle\langle
x,c\rangle}{2\langle x,c\rangle\langle p,x\rangle}=\frac{\langle
p,p\rangle}2\Im\Big(\frac{\langle w,x\rangle}{\langle
p,x\rangle}-\frac{\langle x,x\rangle\langle w,c\rangle}{\langle
p,x\rangle\langle x,c\rangle}\Big).$$
Therefore,
$$v_pP_c\big(\Tn(w_p)\big)=\frac{\langle p,p\rangle}2\frac
d{d\varepsilon}\Big|_{\varepsilon=0}\Im\Big(\frac{\big\langle
w,p+\varepsilon\langle p,p\rangle v\big\rangle}{\big\langle
p,p+\varepsilon\langle p,p\rangle v\big\rangle}-\frac{\big\langle
p+\varepsilon\langle p,p\rangle v,p+\varepsilon\langle p,p\rangle
v\big\rangle\langle w,c\rangle}{\big\langle p,p+\varepsilon\langle
p,p\rangle v\big\rangle\big\langle p+\varepsilon\langle p,p\rangle
v,c\big\rangle}\Big)=$$
$$=\frac{\langle p,p\rangle}2\frac
d{d\varepsilon}\Big|_{\varepsilon=0}\Im\Big(\varepsilon\langle
w,v\rangle-\frac{1+\varepsilon^2\langle p,p\rangle\langle
v,v\rangle}{\langle p,c\rangle+\varepsilon\langle p,p\rangle\langle
v,c\rangle}\langle w,c\rangle\Big)=\frac{\langle
p,p\rangle}2\Im\Big(\langle w,v\rangle+\frac{\langle p,p\rangle\langle
v,c\rangle\langle w,c\rangle}{\langle
p,c\rangle^2}\Big).\eqno{_\blacksquare}$$

{\bf Proposition 4.5.6.} {\sl For every\/ $c\in\Bbb P$, we have\/
$dP_c=\omega$ on\/
$\Bbb P\setminus(\partial_\infty\Bbb B\cup\Bbb Pc^\perp)$.}

\medskip

{\bf Proof.} Let
$p\in\Bbb P\setminus(\partial_\infty\Bbb B\cup\Bbb Pc^\perp)$ and let
$v_p,w_p\in\T_p\Bbb P$. Take $X:=\Tn(v_p)$ and $Y:=\Tn(w_p)$. By~Lemma
4.5.4, $[X,Y](p)=0$. Hence, $P_c\big([X,Y]\big)(p)=0$. By the
Maurer-Cartan identity,
$$dP_c(X,Y)(p)=X\big(P_c(Y)\big)(p)-Y\big(P_c(X)\big)(p)-
P_c\big([X,Y]\big)(p).$$
It follows from Lemma 4.5.5 that
$$dP_c(v_p,w_p)=\frac{\langle p,p\rangle}2\Im\Big(\langle
w,v\rangle+\frac{\langle p,p\rangle\langle v,c\rangle\langle
w,c\rangle}{\langle p,c\rangle^2}\Big)-\frac{\langle
p,p\rangle}2\Im\Big(\langle v,w\rangle+\frac{\langle p,p\rangle\langle
w,c\rangle\langle v,c\rangle}{\langle p,c\rangle^2}\Big)=$$
$$=\frac12\Im\big(\langle p,p\rangle\langle w,v\rangle-\langle
p,p\rangle\langle v,w\rangle\big)=\Im\langle
v_p,w_p\rangle.\eqno{_\blacksquare}$$

\bigskip

\centerline{\bf5.~References}

\medskip

[AGr1] S.~Anan$'$in, C.~H.~Grossi, {\it Coordinate-free classic
geometries,} Mosc.~Math.~J., to appear, see also
http://arxiv.org/abs/math/0702714.

[AGr2] S.~Anan$'$in, C.~H.~Grossi, {\it Invariants of finite
configurations in\/ $\Bbb H_\Bbb C^2$ and their geometric
interpretation,} in preparation.

[AGr3] S.~Anan$'$in, C.~H.~Grossi, {\it Yet another Poincar\'e
polyhedron theorem,} Proc.~Edinburgh Math. Soc.~{\bf54} (2011),
297--308, see also http://arxiv.org/abs/0812.4161.

[BSi] F.~Bonahon, L.~Siebenmann, {\it The classification of Seifert
fibered 3-orbifolds,} in {\it Low Dimensional Topology,} edited by
R.~Fenn, London Mathematical Society LNS {\bf95}, Cambridge University
Press, New York, 1985, 258 pp.

[GKL] W.~M.~Goldman, M.~Kapovich, B.~Leeb, {\it Complex hyperbolic
manifolds homotopy equivalent to a Riemann surface,}
Comm.~Anal.~Geom.~{\bf9}, no.~1 (2001), 61--95.

[GLT] M.~Gromov, H.~B.~Lawson Jr., W.~Thurston, {\it Hyperbolic
$4$-manifolds and conformally flat $3$-manifolds,} Inst.~Hautes
\'Etudes Sci.~Publ.~Math.~{\bf68} (1988), 27--45.

[Gol] W.~M.~Goldman, {\it Complex hyperbolic geometry,} Oxford
Mathematical Monographs. Oxford Science Publications. The Clarendon
Press, Oxford University Press, New York, 1999, xx+316 pp.

[HSa] J.~Hakim, H.~Sandler, {\it Application of Bruhat decomposition to
complex hyperbolic geometry,} J.~Geom.~Anal.~{\bf10} (2000), 435--453.

[Hsi] Hsieh, Po-Hsun, {\it Cotranchal bisectors in complex hyperbolic
space,} Geometriae Dedicata {\bf97} (2003), 93--98.

[Kap1] M.~Kapovich, {\it Flat conformal structures on $3$-manifolds.
{\rm I.} Uniformization of closed Seifert manifolds,} J.~Differential
Geom.~{\bf38}, no.~1 (1993), 191--215.

[Kap2] M.~Kapovich, {\it On hyperbolic\/ $4$-manifolds fibered over
surfaces,} preprint (1993). Available at
http://www.math.ucdavis.edu/${\sim}$kapovich/eprints.html.

[KoM] A.~I.~Kostrikin, Y.~I.~Manin, {\it Linear algebra and geometry,}
Algebra, Logic and Applications Series. Gordon and Breach Science
Publishers, London, 1989, xvi+309 pp.

[Kre] M.~Krebs, {\it Toledo invariants of $2$-orbifolds and Higgs
bundles on elliptic surfaces,} Michigan Math.~J.~{\bf56}, no.~1 (2008),
3--27.

[Kui] N.~H.~Kuiper, {\it Hyperbolic $4$-manifolds and tessellations,}
Inst.~Hautes \'Etudes Sci.~Publ.~Math.~{\bf68} (1988), 47--76.

[Luo] F.~Luo, {\it Constructing conformally flat structures on some
Seifert fibred $3$-manifolds,} Math.~Ann. {\bf294}, no.~3 (1992),
449--458.

[Pra] A.~Pratoussevitch, {\it Traces in complex hyperbolic triangle
groups,} Geometriae Dedicata {\bf111} (2005), 159--185.

[San] H.~Sandler, {\it Distance formulas in complex hyperbolic space,}
Forum Math.~{\bf8}, no.~1 (1996), 93--106.

[Sch] S.~Schleimer, {\it Waldhausen's theorem,} Geom.~Topol.~Monographs
{\bf12} (2007), 299--317.

[Tol] D.~Toledo, {\it Representations of surface groups in complex
hyperbolic space,} J.~Differential Geom. {\bf29}, no.~1 (1989),
125--133.

[Xia] E.~Z.~Xia, {\it The moduli of flat $\PU(2,1)$ structures on
Riemann surfaces,} Pacific J.~Math.~{\bf195}, no.~1 (2000), 231--256.

\enddocument